\documentclass[12pt]{amsart}
\usepackage{amssymb}
\usepackage{amsmath}
\usepackage{mathabx}
\usepackage{amsthm}
\usepackage{amstext}
\usepackage{amsopn}
\usepackage{mathrsfs}
\usepackage{latexsym}
\DeclareMathAlphabet{\mathpzc}{OT1}{pzc}{m}{it}
\allowdisplaybreaks
\newtheorem{Definition}{Definition}[section]
\newtheorem{Proposition}{Proposition}[section]
\newtheorem{Lemma}{Lemma}[section]
\newtheorem{Theorem}{Theorem}[section]
\newtheorem{Corollary}{Corollary}[section]
\newtheorem{Remark}{Remark}[section]
\newtheorem{Example}{Example}[section]

\begin{document}
\bibliographystyle{plain}
\title[Limit distributions of random matrices]
{Limit distributions of random matrices}
\author[R. Lenczewski]{Romuald Lenczewski}
\address{Romuald Lenczewski, \newline
Instytut Matematyki i Informatyki, Politechnika Wroc\l{}awska, \newline
Wybrze\.{z}e Wyspia\'{n}skiego 27, 50-370 Wroc{\l}aw, Poland  \vspace{10pt}}
\email{Romuald.Lenczewski@pwr.wroc.pl}
\begin{abstract}
We study limit distributions of independent random matrices as
well as limit joint distributions of their blocks under normalized partial traces 
composed with classical expectation. In particular, we are concerned with the ensemble 
of symmetric blocks of independent Hermitian random matrices which are asymptotically 
free, asymptotically free from diagonal deterministic matrices, and whose norms are uniformly bounded.
This class contains symmetric blocks of unitarily invariant Hermitian random matrices whose asymptotic distributions are
compactly supported probability measures on the real line. Our approach is based on the concept of matricial freeness which is a generalization 
of freeness in free probability. We show that the associated matricially free Gaussian operators provide a unified 
framework for studying the limit distributions of sums and products of independent rectangular random matrices, 
including non-Hermitian Gaussian matrices and matrices of Wishart type. 
\end{abstract}

\maketitle
\section{Introduction and main results}
One of the most important features of free probability is its close relation to random matrices. 
It has been shown by Voiculescu [33] that Hermitian random matrices with independent Gaussian entries 
are asymptotically free. This result has been generalized by Dykema to non-Gaussian random matrices [12]
and has been widely used by many authors in their studies of asymptotic distributions of random matrices. 
It shows that there is a concept of noncommutative independence, called {\it freeness}, which is fundamental 
to the study of large random matrices and puts the classical result of Wigner [36] on the semicircle law
as the limit distribution of certain symmetric random matrices in an entirely new perspective.

In particular, if we are given an ensemble of independent Hermitian $n\times n$ random matrices 
$$
\{Y(u,n):u\in \mathpzc{U}\}
$$ 
whose entries are suitably normalized and independent complex Gaussian random variables for each natural $n$, then 
$$
\lim_{n\rightarrow \infty}\tau(n)(Y(u_1,n)\ldots Y(u_m,n))=
\Phi(\omega(u_1)\ldots \omega(u_m))
$$ 
for any $u_1, \ldots , u_m\in \mathpzc{U}$, where $\{\omega(u):u\in \mathpzc{U}\}$
is a semicircular family of {\it free Gaussian operators} living in the free Fock space with the vacuum state $\Phi$
and $\tau(n)$ is the normalized trace composed with classical expectation called the {\it trace} in the sequel.
This realization of the limit distribution gives a fundamental relation between random matrices and
operator algebras. 

The basic original random matrix model studied by Voiculescu corresponds to independent
complex Gaussian variables, where the entries $Y_{i,j}(u,n)$ of each matrix $Y(u,n)$ satisfy the Hermiticity condition, 
have mean zero, and the variances of real-valued diagonal
Gaussian variables $Y_{j,j}(u,n)$ are equal to $1/n$, whereas those of the real and imaginary parts of the off-diagonal 
(complex-valued) Gaussian variables $Y_{i,j}(u,n)$ are equal to $1/2n$. 
If we relax the assumption on equal variances, the scalar-valued free probability is no longer sufficient to describe
the asymptotics of Gaussian random matrices. One approach is to use the operator-valued free probability, as in the work of 
Shlakhtyenko [32], who studied the asymptotics of Gaussian random band matrices and proved that they are asymptotically free with amalgamation
over some commutative algebra. This approach was further developed by Benaych-Georges [6] who treated the case when one of the 
asymptotic dimensions vanishes.

Our approach is based on the decomposition of independent Hermitian Gaussian random matrices $Y(u,n)$ 
with block-identical variances of $|Y_{i,j}(u,n)|$ into symmetric blocks $T_{p,q}(u,n)$, where $u\in \mathpzc{U}$ and $n\in \mathbb{N}$, namely
$$
Y(u,n)=\sum_{1\leq p \leq q \leq r}T_{p,q}(u,n),
$$
where symmetric blocks can be written in the form
$$
T_{p,q}(u,n)=\left\{
\begin{array}{ll}
D_{q}Y(u,n)D_{q}& {\rm if}\;\; p=q\\
D_{p}Y(u,n)D_{q}+D_{q}Y(u,n)D_{p}& {\rm if}\;\;p< q
\end{array}
\right.
$$
where $\{D_1, \ldots , D_r\}$ is the family of $n\times n$ diagonal matrices that 
forms a decomposition of the $n\times n$ identity matrix corresponding to 
the partition of $[n]:=\{1, \ldots , n\}$ into disjoint nonempty intervals
$N_1, \ldots , N_r$ and we assume that
$$
\lim_{n\rightarrow \infty}\frac{n_{q}}{n}= d_{q}\geq 0
$$
for any $q$, where $n_q$ denotes the cardinality of $N_q$. Of course,  
the quantities $N_q,n_q,D_q$ depend on $n$, but this is suppressed in our notation.

At the same time, we decompose each semicircular Gaussian operator in terms of the corresponding $r\times r$ array of 
{\it matricially free Gaussian operators} [23], namely
$$
\omega(u)=\sum_{p,q=1}^{r}\omega_{p,q}(u),
$$
living in the matricially free Fock space of tracial type, similar to the matricially free Fock space introduced in [23], 
in which we distinguish a family of states $\{\Psi_1, \ldots , \Psi_r\}$ associated
with a family of vacuum vectors. These states are used to build the array $(\Psi_{p,q})$ by setting $\Psi_{p,q}=\Psi_q$.

Nevertheless, in order to reproduce the limit distributions of symmetric blocks in the case when the variances 
are only block-identical rather than identical, which is obtained by a rescaling, we also need to rescale 
the matricially free Gaussian operators in the above expressions. The corresponding arrays of distributions 
$$
[\sigma(u)]=
\left(
\begin{array}{llll}
\sigma_{1,1}(u) & \kappa_{1,2}(u) & \ldots & \kappa_{1,r}(u)\\
\kappa_{2,1}(u) & \sigma_{2,2}(u) & \ldots & \kappa_{2,r}(u)\\
. & . & \ddots & . \\
\kappa_{r,1}(u) & \kappa_{r,2}(u) & \ldots & \sigma_{r,r}(u)
\end{array}
\right),
$$
consisting of semicircle laws $\sigma_{q,q}(u)$ and Bernoulli laws $\kappa_{p,q}(u)$, replace the semicircle laws
$\sigma(u)$ and play the role of {\it matricial semicircle laws}. 
In turn, the rescaled matricially free Gaussians $(\omega_{p,q}(u))$ 
replace free Gaussians $\omega(u)$ and thus the corresponding 
sums become more general objects called {\it Gaussian pseudomatrices}.

In any case, it is the family of {\it symmetrized Gaussian operators} defined by
$$
\widehat{\omega}_{p,q}(u)=\left\{
\begin{array}{ll}
\omega_{q,q}(u)& {\rm if}\; p=q\\
\omega_{p,q}(u)+\omega_{q,p}(u)& {\rm if}\; p<q
\end{array}\right.
$$
which gives the operatorial realizations of the limit joint distributions of the ensemble of 
symmetric blocks of the family $\{Y(u,n):u\in \mathpzc{U}\}$ of Hermitian Gaussian random matrices (HGRM) in 
the case when the complex random variables in each symmetric block are identically distributed. 

In this setting, we obtain an operatorial realization of the limit joint distributions
of symmetric blocks of independent HGRM under {\it partial traces} 
$\tau_{q}(n)$, by which we understand normalized traces over the sets of basis vectors indexed by $N_q$, 
respectively, composed with the classical expectation. It takes the form 
$$
\lim_{n\rightarrow \infty}\tau_{q}(n)(T_{p_1,q_1}(u_1,n)\ldots T_{p_m,q_m}(u_m,n))=
\Psi_{q}(\widehat{\omega}_{p_1,q_1}(u_1)\ldots \widehat{\omega}_{p_m,q_m}(u_m)),
$$ 
where  $u_1, \ldots , u_m\in \mathpzc{U}$ and the remaining indices belong to $[r]$. 
Next, from the result for symmetric blocks of HGRM we can derive the limit joint distribution of
symmetric blocks of independent (non-Hermitian) Gaussian random matrices (GRM) as in the case of 
GRM themselves given in [33]. 

This result can be viewed as a block refinement of that used by Voiculescu in his fundamental asymptotic freeness result.
Using the convex linear combination 
$$
\Psi=\sum_{q=1}^{r}d_{q}\Psi_q,
$$
we easily obtain a similar formula for $\tau(n)$ 
considered by most authors in their studies of asymptotics of random matrices. It is obvious that asymptotic freeness is a special case of 
our asymptotic matricial freeness and it corresponds to the case when the variances of all $|Y_{i,j}(u,n)|$ 
are identical. However, we also show that by considering partial traces we can produce random 
matrix models for boolean independence, monotone independence and s-freeness. It is not a coincidence since all these 
notions of independence arise in the context of suitable decompositions of free random variables as shown in [21,22].

Let us remark that similar distributions for one HGRM were derived in our previous paper [23]. 
Here, we generalize this result in several important directions:
\begin{enumerate}
\item
we study symmetric blocks of families of independent random matrices rather than of one random matrix,  
\item
we consider the class of Hermitian random matrices which are asymptotically free, 
asymptotically free from $\{D_1, \ldots, D_r\}$, and whose norms are uniformly bounded,
\item
we include the case when some elements of the dimension matrix
$$
D={\rm diag}(d_1,d_2, \ldots , d_r),
$$
called {\it asymptotic dimensions}, are equal to zero,
\item
we observe that similar results are obtained if we take symmetric blocks of random matrices of dimension $d(n)\times d(n)$ for $d(n)\rightarrow \infty$
as $n\rightarrow \infty$.
\end{enumerate}

The class of Hermitian random matrices mentioned above includes unitarily inva\-riant matrices converging in distribution to 
compactly supported probability measures on the real line which are known to be asymptotically free and asymptotically free from 
diagonal deterministic matrices [18]. Joint distributions of their blocks were studied by Benaych-Georges [6] by means 
of freeness with amalgamation under the additional assumption that the singular values of $D_pY(u,n)D_q$ are deterministic.

It is especially interesting to study off-diagonal symmetric blocks consisting 
of two rectangular blocks whose dimensions do not grow proportionately, with the larger dimension 
proportional to $n$ and the smaller dimension growing slower than $n$. Such blocks will be called {\it unbalanced} in contrast
to the {\it balanced} ones, in which both dimensions are proportional to $n$.
In the formulas for joint distributions under the trace $\tau(n)$, the contributions from partial traces associated with 
vanishing asymptotic dimensions disappear in the limit. However, certain limit joint distributions involving unbalanced symmetric blocks 
under the partial traces themselves become non-trivial and interesting.

It should be pointed out that one can derive the limit joint distributions of all 
balanced blocks, using ordinary freeness (with scalar-valued states)
since the family $\{Y(u,n):u\in \mathpzc{U}\}$ is asymptotically free from $\{D_1, \ldots , D_r\}$ under $\tau(n)$ and
$$
\tau_{q}(n)(X)=\frac{n}{n_{q}}\tau(n)(D_{q}XD_{q})
$$
for any $n\times n$ random matrix $X$ and any $q\in [r]$.
However, when we have unbalanced blocks, this does not seem possible. In particular, 
if $d_q=0$, the joint distribution of symmetric blocks under $\tau_{q}(n)$ is non-trivial, whereas
its contribution to that under $\tau(n)$ is zero. Nevertheless, this `singular' case seems to be of importance since it
leads to new random matrix models. For instance, it gives other types of asymptotic independence as those 
mentioned above as well as a new random matrix model for free Meixner laws and the associated 
asymptotic conditional freeness [24].

We also show that the study of non-Hermitian Gaussian random matrices (sometimes called the Ginibre Ensemble) 
can be reduced to the Hermitian case as in Voiculescu's paper. Further, by considering off-diagonal blocks of our matrices, we can recover 
limit joint distributions of independent Wishart matrices [37]
$$
W(n)=B(n)B(n)^*
$$
where $B(n)$ is a complex Gaussian matrix, or a random matrix of similar type, 
where each $B(n)$ is a product of independent Gaussian random matrices
$$
B(n)=Y(u_1,n)Y(u_2,n)\ldots Y(u_k,n),
$$
provided their dimensions are such that the products are well-defined. These results
have led to the formula for the moments of the multiplicative free convolution of 
Marchenko-Pastur laws of arbitrary shape paramaters, expressed in terms of mutlivariate 
Fuss-Narayana polynomials [26].

In all these situations, it suffices to embed the considered matrices
in the algebra of symmetric blocks in an appropriate way. In particular, we obtain
the random matrix model for the noncommutative Bessel laws of Banica {\it et al} [4] in this fashion.
Non-Hermitian Wishart matrices, which are also of interest [19], can be treated in a similar way.
We also expect that a continuous generalization of our block method will allow us to treat 
symmetric triangular random matrices studied by Basu {\it et al} [5] and the triangular ones 
studied by Dykema and Haagerup [13].

The paper is organized as follows. Section 2 is devoted to 
the concept of matricial freeness, the corresponding arrays of Gaussian 
operators and matricial semicircle laws. Symmetrized counterparts of these objects
are discussed in Section 3. In Section 4, we describe the combinatorics of the mixed 
moments of Gaussian operators and their symmetrized counterparts. Realizations of 
canonical noncommutative random variables in terms of matricially free creation and annihilation operators is given in
Section 5. The first main result of this paper is contained in Section 6, where we derive the limit joint 
distributions of symmetric blocks of a large class of independent Hermitian random matrices under partial traces. 
In Section 7, we construct random matrix models for boolean, monotone and s-free independences.  
Non-Hermitian Gaussian random matrices are treated in Section 8. Section 9 is devoted to Wishart matrices, 
matrices of Wishart type and products of independent Gaussian random matrices.

Finally, we take this opportunity to issue an erratum to the definition
of the symmetrically matricially free array of units given in [23]. The new 
definition entails certain changes which are discussed in the Appendix.

\section{Matricial semicircle laws}

Let us recall the basic notions related to the concept of matricial freeness [22].
Let ${\mathcal A}$ be a unital *-algebra with an array $({\mathcal A}_{i,j})$ 
of non-unital *-subalgebras of ${\mathcal A}$ and let $(\varphi_{i,j})$ be an array 
of states on ${\mathcal A}$. Here, and in the sequel, we shall skip $\mathpzc{J}$ 
in the notations involving arrays and we shall tacitly assume that 
$(i,j)\in \mathpzc{J}$, where
\begin{equation*}
\mathpzc{J}\subseteq [r]\times [r],
\end{equation*}
and $r$ is a natural number. Further, we assume that each 
${\mathcal A}_{i,j}$ has an {\it internal unit} $1_{i,j}$, by which we understand
a projection for which $a1_{i,j}=1_{i,j}a=a$ for any $a\in {\mathcal A}_{i,j}$, and
that the unital subalgebra ${\mathcal I}$ of ${\mathcal A}$ generated by all internal 
units is commutative.

For each natural number $m$, let us distinguish the subset of the $m$-fold 
Cartesian product $\mathpzc{J}\times \mathpzc{J} \times \ldots \times \mathpzc{J}$
of the form
\begin{equation*}
\mathpzc{J}_{\;m}=\{((i_1,j_1),  \ldots , (i_m,j_m)): \;(i_1,j_1)\neq \ldots \neq (i_m,j_m)\},
\end{equation*}
and its subset
\begin{equation*}
\mathpzc{K}_{\;m}=\{((i_1,i_2), (i_2,i_3)\ldots , (i_{m},j_m)):\;(i_1,i_2)\neq (i_2,i_3)\neq \ldots \neq (i_m,j_m)\}.
\end{equation*}
In other words, the neighboring pairs of indices in the set $\mathpzc{K}_{\;m}$ are not only different (as in free products), but
are related to each other as in matrix multiplication. 
Objects labelled by diagonal or off-diagonal pairs, respectively, will be called {\it diagonal} and {\it off-diagonal}.

If $\varphi_1, \varphi_2, \ldots , \varphi_r$ are states on ${\mathcal A}$, we form an array of states 
$(\varphi_{i,j})$ as follows:
\begin{equation*}
\varphi_{i,j}=\varphi_{j}
\end{equation*}
and then we will say that $(\varphi_{i,j})$ is {\it defined by the family $(\varphi_j)$}.
In particular, when $\mathpzc{J}=[r]\times [r]$, this array takes the form
\begin{equation*}
(\varphi_{i,j})=
\left(
\begin{array}{llll} 
\varphi_1 & \varphi_2 & \ldots & \varphi_r\\
\varphi_1 & \varphi_2 & \ldots  & \varphi_r\\
.   & .& \ddots & .\\
\varphi_1 & \varphi_2 & \ldots  & \varphi_r
\end{array}
\right).
\end{equation*}
The states $\varphi_1, \varphi_2, \ldots \varphi_r$ can be defined as states which are {\it conjugate} 
to a distinguished state $\varphi$ on ${\mathcal A}$, where conjugation is implemented by certain elements
of the diagonal subalgebras, namely
\begin{equation*}
\varphi_{j}(a)=\varphi(b_{j}^{*}ab_{j})
\end{equation*}
for any $a\in {\mathcal A}$, where $b_{j}\in {\mathcal A}_{j,j}\cap {\rm Ker}\varphi$ is such that 
$\varphi(b_{j}^{*}b_{j})=1$ for any $j\in [r]$. 

In general, other arrays can also be used in the definition of matricial freeness,
but in this paper we will only use those defined by a family of $r$ states, say $(\varphi_j)$, and thus the definition 
of matricial freeness is adapted to this situation. Since the original definition uses the diagonal states [22], 
which here coincide with $\varphi_1, \ldots , \varphi_r$, we shall use the latter in all definitions.

\begin{Definition}
{\rm 
We say that $(1_{i,j})$ is a {\it matricially free array of units} associated with 
$({\mathcal A}_{i,j})$ and $(\varphi_{i,j})$ if for any state $\varphi_j$ it holds that  
\begin{enumerate}
\item[(a)]
$\varphi_j(b_1ab_2)=\varphi_j(b_1)\varphi_j(a)\varphi_j(b_2)$ 
for any $a\in {\mathcal A}$ and $b_1,b_2\in {\mathcal I}$,
\item[(b)]
$\varphi_{j}(1_{k,l})=\delta_{j,l}$ for any $i,j,k,l$,
\item[(c)]
if $a_{r}\in {\mathcal A}_{i_r,j_r}\cap {\rm Ker}\varphi_{i_r,j_r}$, where $1<r\leq m$, then
$$
\varphi_j(a1_{i_1,j_1}a_{2}\ldots a_m)=
\left\{
\begin{array}{ll}
\varphi_j(aa_{2} \ldots a_m) & {\rm if}\;\;((i_1,j_1), \ldots , (i_m,j_m))\in \mathpzc{K}_{\;m}\\
0 & {\rm otherwise}
\end{array}
\right..
$$
where $a\in {\mathcal A}$ is arbitrary and 
$((i_1,j_1), \ldots , (i_m,j_m))\in \mathpzc{J}_{\;m}$.
\end{enumerate}}
\end{Definition}

\begin{Definition}
{\rm 
We say that *-subalgebras $({\mathcal A}_{i,j})$ are 
{\it matricially free} with respect to $(\varphi_{i,j})$ if the array of internal units 
$(1_{i,j})$ is the associated matricially free array of units and
\begin{equation*}
\varphi_j(a_{1}a_{2}\ldots a_{n})=0\;\;\; {\rm whenever }\;\;a_{k}\in {\mathcal A}_{i_k,j_k}\cap {\rm Ker}\varphi_{i_k,j_k}
\end{equation*}
for any state $\varphi_j$, where $((i_1,j_1), \ldots , (i_m,j_m))\in \mathpzc{J}_{\;m}$.
The matricially free array of variables $(a_{i,j})$ in a unital *-algebra ${\mathcal A}$ 
is defined in a natural way.}
\end{Definition}

The most important example of matrically free operators are the matricial 
generalizations of free Gaussian operators defined on a suitable Hilbert space 
of Fock type.

\begin{Definition}
{\rm 
Let $\{\mathcal{H}_{p,q}(u):(p,q)\in \mathpzc{J}, u\in \mathpzc{U}\}$
be a family of arrays of Hilbert spaces and let $\Omega=\{\Omega_1, \ldots, \Omega_r\}$ be a family of unit vectors.
By the {\it matricially free Fock space of tracial type} 
we understand the direct sum of Hilbert spaces
\begin{equation*}
{\mathcal M}= \bigoplus_{q=1}^{r} {\mathcal M}_{q},
\end{equation*}
where each summand is of the form
\begin{equation*}
{\mathcal M}_{q}={\mathbb C}\Omega_{q}\oplus \bigoplus_{m=1}^{\infty}
\bigoplus_{\stackrel{p_1,\ldots , p_m\in [r]}
{\scriptscriptstyle u_1, \ldots , u_n\in \mathpzc{U}}
}
{\mathcal H}_{p_1,p_2}(u_{1})\otimes {\mathcal H}_{p_2,p_3}(u_2)
\otimes \ldots \otimes 
{\mathcal H}_{p_m,q}(u_{m})
\end{equation*}
where the summation extends over those indices which give tensor products of 
Hilbert spaces associated with pairs which belong to $\mathpzc{J}$. The direct sum is 
endowed with the canonical inner product.}
\end{Definition}

Each space ${\mathcal M}_{q}$ is similar to the matricially free-boolean Fock space introduced in [23] except 
that its tensor products do not have to end with diagonal free Fock spaces and that we use
a family of arrays of Hilbert spaces indexed by $u\in \mathpzc{U}$.
Thus, ${\mathcal M}$ is equipped with a family of (vacuum) unit vectors. The 
associated states will be denoted 
$$
\Psi_{q}(a)=\langle a\Omega_{q}, \Omega_{q}\rangle
$$ 
for any $a\in B({\mathcal M})$ and any $q\in [r]$. This family is used to define an array of states.
Thus, let $(\Psi_{p,q})$ be the array defined by the family $(\Psi_q)$, thus $\Psi_{p,q}=\Psi_q$ 
for any $(p,q)\in \mathpzc{J}$.

\begin{Remark}
{\rm If ${\mathcal H}_{p,q}={\mathbb C}e_{p,q}(u)$ for any $p,q,u$, then the orthonormal basis ${\mathcal B}$ of the 
associated matricially free Fock space ${\mathcal M}$ is given by vacuum vectors $\Omega_{1}, \ldots ,\Omega_r$, and simple tensors of the form
$$
e_{p_1,p_2}(u_1)\otimes e_{p_2,p_3}(u_2)\otimes  \ldots \otimes e_{p_m,q}(u_{m})
$$
where $p_1, \ldots , p_m,q\in [r]$, $u_{1},\ldots , u_m\in \mathpzc{U}$ and it is assumed that 
$(p_1,p_2), \ldots , (p_m,q)\in \mathpzc{J}$.
}
\end{Remark}

\begin{Definition}
{\rm Let $(b_{p,q}(u))$ be an array of non-negative real numbers for any $u\in \mathpzc{U}$.
We associate with it the arrays of 
\begin{enumerate}
\item
{\it matricially free creation operators} $(\wp_{p,q}(u))$ 
by setting
\begin{eqnarray*}
\wp_{p,q}(u)\Omega_{q}&=&
\sqrt{b_{p,q}(u)}e_{p,q}(u)\\
\wp_{p,q}(u)(e_{q,s}(t)\otimes w)&=&
\sqrt{b_{p,q}(u)}(e_{p,q}(u)\otimes e_{q,s}(t)\otimes w)
\end{eqnarray*}
for any $e_{q,s}(t)\otimes w\in {\mathcal B}$, where $(p,q),(q,s)\in \mathpzc{J}$ and $u,t\in \mathpzc{U}$, and 
their action onto the remaining vectors is set to be zero,
\item
{\it matricially free annihilation operators} $(\wp^{*}_{p,q}(u))$ consisting of their adjoints,
\item
{\it matricially free Gaussian operators} $(\omega_{p,q}(u))$ consisting of the sums 
\begin{equation*}
\omega_{p,q}(u)=\wp_{p,q}(u)+\wp^{*}_{p,q}(u).
\end{equation*}
\end{enumerate}
and the matrix $(\omega_{p,q}(u))$ will be denoted by $[\omega(u)]$. The number $b_{p,q}(u)$ will be called 
the {\it covariance} of $\wp_{p,q}(u)$ and the {\it variance} of $\omega_{p,q}(u)$ w.r.t. $\Psi_q$.}
\end{Definition}

\begin{Remark}
{\rm 
Note that the above definition is very similar to that in [23] and for that reason we use similar notations and 
the same terminology. Strictly speaking, however, the above operators live in a larger Fock space, built by the action 
of a {\it collection} of arrays of creation operators onto the family of vacuum vectors $\Omega$
rather than by the action of one array onto one distinguished vacuum vector. Let us also remark that if 
all covariances are equal to one, then the operators $\wp_{p,q}(u)$ are partial isometries which 
generate {\it Toeplitz-Cuntz-Krieger algebras} [25] and thus ther are nice objects to work with. 
Finally, in contrast to [23], we allow the covariances to vanish, which leads to trivial operators.}
\end{Remark}

\begin{Example}
{\rm The moments of each operator $\omega_{p,q}(u)$ will be easily obtained from Proposition 2.3
(cf. [23, Proposition 4.2]). Therefore, let us compute some simple {\it mixed} moments. First, we take only off-diagonal operators:
\begin{eqnarray*}
\Psi_{q}(\omega_{p,q}(u)\omega_{k,p}^{4}(t)\omega_{p,q}(u))&=&
\Psi_{q}(\wp_{p,q}^{*}(u)\wp_{k,p}^{*}(t)\wp_{k,p}(t)\wp_{k,p}^{*}(t)\wp_{k,p}(t)\wp_{p,q}(u))\\
&=&b_{k,p}^{2}(t)b_{p,q}(u)
\end{eqnarray*}
where $p\neq q \neq k \neq p$ and $t,u\in \mathpzc{U}$ are arbitrary.
If some operators are diagonal, we usually obtain
more terms. For instance
\begin{eqnarray*}
\Psi_{q}(\omega_{p,q}(u)\omega_{p,p}^{4}(t)\omega_{p,q}(u))
&=&\Psi_{q}(\wp_{p,q}^{*}(u)\wp_{p,p}^{*}(t)\wp_{p,p}(t)\wp_{p,p}^{*}(t)\wp_{p,p}(t)\wp_{p,q}(u))\\
&+&
\Psi_{q}(\wp_{p,q}^{*}(u)\wp_{p,p}^{*}(t)\wp_{p,p}^{*}(t)\wp_{p,p}(t)\wp_{p,p}(t)\wp_{p,q}(u))
\\
&=&2b_{p,p}^{2}(t)b_{p,q}(u)
\end{eqnarray*}
for any $p\neq q$ and any $u,t\in \mathpzc{U}$. Let us observe that the pairings between the diagonal annihilation and creation 
operators and between the off-diagonal ones are the same as those between annihilation and creation operators 
in the free case and in the boolean case, respectively.}
\end{Example}

Now, we would like to prove matricial freeness of {\it collective Gaussian operators }
of the form
$$
\omega_{p,q}=\wp_{p,q}+\wp_{p,q}^{*},
$$
where
$$
\wp_{p,q}=\sum_{u\in \mathpzc{U}}\wp_{p,q}(u)\;\;{\rm and}\;\;\wp_{p,q}^{*}=\sum_{u\in \mathpzc{U}}\wp_{p,q}^{*}(u),
$$
with respect to the array $(\Psi_{p,q})$.
As in [23, Proposition 4.2], we will prove a more general result on matricial freeness of the array 
of algebras $({\mathcal A}_{p,q})$ generated by collective creation, annihilation and unit operators.

For that purpose, let us define a suitable array of collective units $(1_{p,q})$. The easiest way to do it is to say that 
$1_{p,q}$ is the orthogonal projection onto the subspace onto which the associated collective 
creation or annihilation operators act non-trivially. However, a more explicit definition 
will be helpful.

\begin{Definition}
{\rm 
By the {\it array of collective units} we will understand $(1_{p,q})$, where 
$$
1_{p,q}=r_{p,q}+s_{p,q} \;\;{\rm for}\;\;(p,q)\in \mathpzc{J},
$$
where $s_{p,q}$ is the orthogonal projection onto 
$$
{\mathcal N}_{p,q}:=\bigoplus_{u\in \mathpzc{U}}{\rm Ran}(\wp_{p,q}(u)),
$$ 
and $r_{p,q}$ is the orthogonal projection onto
\begin{eqnarray*}
{\mathcal K}_{q,q}&=&{\mathbb C}\Omega_{q}\oplus \bigoplus_{k\neq q}\mathcal{N}_{q,k},\\
{\mathcal K}_{p,q}&=&{\mathbb C}\Omega_{q}\oplus \bigoplus_{k}\mathcal{N}_{q,k},
\end{eqnarray*}
for any diagonal $(q,q)\in \mathpzc{J}$ and off-diagonal $(p,q)\in \mathpzc{J}$, respectively.
Let us remark that in the latter case $r_{p,q}=1_{q,q}$.}
\end{Definition}

\begin{Proposition}
The following relations hold: 
$$
\wp_{p,q}(u)\wp_{k,l}(t)=0\;\;and\;\;\wp_{p,q}^{*}(u)\wp_{k,l}^{*}(t)=0
$$ 
for $q\neq k$ and any $u,t\in \mathpzc{U}$. Moreover,
$$
\wp_{p,q}^{*}(u)\wp_{p,q}(u)=b_{p,q}(u)1_{p,q}
$$
for any $p,q,u$, and otherwise $\wp_{p,q}^{*}(u)\wp_{k,l}(t)=0$ for any $(p,q,u)\neq (k,l,t)\in \mathpzc{J}\times \mathpzc{U}$.
\end{Proposition}
{\it Proof.}
These relations follow from Proposition 2.1.
\hfill $\blacksquare$\\

One can give a nice tensor product realization of matricially free creation and annihilation operators [25].
Let ${\mathcal A}$ be the unital *-algebra generated by the system of *-free creation operators,
$$
\{\ell(p,q,u): 1\leq p, q\leq [r], u\in \mathpzc{U}\},
$$ 
and let $\varphi$ be the vacuum state on ${\mathcal A}$. In general, we assume that the covariances 
of these creation operators are arbitrary nonnegative numbers, that is 
$$
\ell(p,q,u)^*\ell(p,q,u)=b_{p,q}(u)\geq 0.
$$
The case of trivial covariance is included for convenience. This enables us to 
suppose that ${\mathcal J}=[r]\times [r]$ even if some of the operators in the considered arrays are zero.

The algebra of $r\times r$ matrices with entries in ${\mathcal A}$, namely $M_{r}({\mathcal A})\cong {\mathcal A}\otimes M_{r}({\mathbb C})$,
is equipped with the natural involution 
$$
(a\otimes e(p,q))^{*}=a^{*}\otimes e(q,p)
$$
where $a\in {\mathcal A}$ and $p,q\in [r]$ and $(e(p,q))$ is the $r\times r$ array of matrix units in $M_r({\mathbb C})$. 
Consider the states $\Phi_1, \ldots , \Phi_r$ on $M_{r}({\mathcal A})$ 
of the form
$$
\Phi_q=\varphi\otimes \psi_q
$$
for any $q\in [r]$, where $\psi_q(b)=\langle b e(q), e(q)\rangle$ and
$\{e(j):j\in [r]\}$ is the canonical orthonormal basis in ${\mathbb C}^{r}$. 

The result given below was proved in [25], but we provide the proof for the reader's convenience.
\begin{Lemma}
The *-distributions of the arrays 
$$
\{\ell(p,q,u)\otimes e(p,q):p,q\in [r], u\in \mathpzc{U}\}
$$
with respect to the states $\Phi_q$ agree with the corresponding *-distributions of the 
arrays of matricially free creation operators
$$
\{\wp_{p,q}(u):p,q\in [r], u \in \mathpzc{U}\}
$$
with respect to the states $\Psi_q$, respectively.
\end{Lemma}
{\it Proof.}
For simplicity, we assume that the covariances of all $\ell(p,q,u)$ are equal to one. 
Let ${\mathcal F}({\mathcal H})$ be the free Fock space over the direct sum of Hilbert spaces
$$
{\mathcal H}=\bigoplus_{1\leq p,q\leq r}\bigoplus_{u\in\mathpzc{U}}{\mathcal H}(p,q,u),
$$ 
where ${\mathcal H}(p,q,u)={\mathbb C}e(p,q,u)$ for any $p,q,u$, with each $e(p,q,u)$ being 
a unit vector, and denote by $\Omega$ the vacuum vector in ${\mathcal F}({\mathcal H})$.
Define an isometric embedding 
$$
\tau:{\mathcal M}\rightarrow {\mathcal F}({\mathcal H})\otimes {\mathbb C}^{r}
$$
by the formulas
\begin{eqnarray*}
\tau(\Omega_{q})&=&\Omega\otimes e(q)\\
\tau(e_{p_1,p_2}(u_1)\otimes \ldots \otimes e_{p_n,q}(u_n))&=
&e(p_1,p_2,u_1)\otimes\ldots \otimes e(p_n,q,u_n)\otimes e(p_1),
\end{eqnarray*}
for any $q,p_1, \ldots , p_n$ and $u_1, \ldots , u_n$, where
$\{e(1),\ldots , e(r)\}$ is the canonical basis in ${\mathbb C}^{r}$.
We then have
$$
\tau \wp_{p,q}(u)=(\ell(p,q,u)\otimes e(p,q))\tau 
$$
since
$$
\tau(\wp_{p,q}(u)\Omega_q)=\tau(e_{p,q}(u))=e(p,q,u)\otimes e(p)
$$
$$
=(\ell(p,q,u)\otimes e(p,q))(\Omega\otimes e(q))=
(\ell(p,q,u)\otimes e(p,q))\tau(\Omega_{q})
$$
and
\begin{eqnarray*}
&&\tau(\wp_{p,q}(u)(e_{q,p_2}(u_1)\otimes \ldots \otimes e_{p_n,t}(u_n))\\
&=&\tau(e_{p,q}(u)\otimes e_{q,p_2}(u_1)\otimes \ldots \otimes e_{p_n,t}(u_n))\\
&=&e(p,q,u)\otimes e(q,p_2,u_1)\otimes \ldots 
\otimes e(p_n,t,u_n)\otimes e(p)\\
&=&(\ell(p,q,u)\otimes e(p,q))
(e(q,p_2,u_1)\otimes \ldots 
\otimes e(p_n,t,u_n)\otimes e(q))\\
&=&
(\ell(p,q,u)\otimes e(p,q))
\tau(e_{q,p_2}(u_1)\otimes \ldots \otimes e_{p_n,t}(u_n))
\end{eqnarray*}
for any values of indices and arguments, whereas the actions 
onto the remaining basis vectors gives zero.
This proves that $\wp_{p,q}(u)$ intertwines with $\ell(p,q,u)\otimes e(p,q)$.
Therefore, the *-distributions of $\wp_{p,q}(u)$ under the states $\Psi_k$ 
agree with the corresponding *-distributions of $\ell(p,q,u)\otimes e(p,q)$ under the states $\Phi_k$, respectively,
which finishes the proof. 
\hfill $\blacksquare$\\

Using the above, it is easy to see that the matricially free Gaussian operators give a decomposition of free semicircular operators
even if the variances of the summands are different.
This fact will be used in the results related to asymptotic freeness of Hermitian random matrices.
 
\begin{Proposition}
Let $\{\omega_{p,q}(u):p,q \in [r], u\in \mathpzc{U}\}$ be a family of matricially free Gaussian operators
w.r.t. $(\Psi_{p,q})$, where $\omega_{p,q}(u)$ has variance $d_p\geq 0$ for any $p,q,u$, 
where $d_1+\ldots +d_r=1$.
Then
\begin{enumerate}
\item 
the operator of the form
$$
\omega(u)=\sum_{p,q=1}^{r} \omega_{p,q}(u)
$$
has the standard semicircular distribution w.r.t. $\Psi=\sum_{q}d_q\Psi_q$,
\item
it holds that $\wp(u)^*\wp(v)=\delta_{u,v}1$, where
$$
\wp(u)=\sum_{p,q=1}^{r}\wp_{p,q}(u)\;\;\;and\;\;\;
\wp(u)^{*}=\sum_{p,q=1}^{r}\wp_{p,q}^{*}(u)
$$
for any $u$,
\end{enumerate}
\end{Proposition}
{\it Proof.}
Using Lemma 2.1, we obtain
\begin{eqnarray*}
\wp(u)^*\wp(v)
&=&
(\sum_{p,q}\ell(p,q,u)^*\otimes e(q,p))(\sum_{i,j}\ell(i,j,v)\otimes e(i,j))\\
&=&
\sum_{p,q,j}\ell(p,q,u)^*\ell(p,j,v)\otimes e(q,j)\\
&=&
\delta_{u,v}\sum_{p,q}\ell(p,q,u)^*\ell(p,q,u)\otimes e(q,q)\\
&=&
\delta_{u,v}\sum_{q}\left(\sum_{p}d_{p}\right)(1\otimes e(q,q))\\
&=&
\delta_{u,v}1
\end{eqnarray*}
which proves (2). Then, (1) easily follows since each $\wp_{p,q}(u)^*$ kills each vector $\Omega_j=\Omega\otimes e(j)$. 
This completes the proof.
\hfill $\blacksquare$\\

An important fact in the asymptotics of random matrices is that 
certain families are asymptotically free from deterministic (in particular, diagonal) matrices. 
For that purpose, we would like to find a suitable family of operators on ${\mathcal M}$ which would 
give a limit realization of the family $\{D_1, \ldots , D_r\}$ of diagonal matrices defined in the Introduction.
This family will be denoted by $\{P_1, \ldots, P_r\}$.

\begin{Lemma}
Let $\{\wp(u):u\in \mathpzc{U}\}$ be as in Proposition 2.2 and
let  $P_q=1\otimes e(q,q)$ for any $q\in [r]$. Then the family 
$\{\wp(u):u\in \mathpzc{U}\}$ is *-free from the family $\{P_1, \ldots, P_r\}$ with respect
to the state $\Psi=\sum_{q}d_q\Psi_q$, where $d_1+\ldots + d_r=1$.
\end{Lemma}
{\it Proof.}
We need to show that 
$$
\Psi(Q_{q_0}a_1(u_1)Q_{q_1}a_2(u_2)\ldots Q_{m-1}a_m(u_m)Q_{q_m})=0
$$ 
whenever each  $a_q(u_q)\in {\rm Ker \Psi}$ is a polynomial in $\wp(u_q)$, $q\in [m]$, whereas each $Q_{q_j}$ is 
either equal to one (then we have $u_j\neq u_{j+1}$), or $Q_{q_j}=P_{q_j}-\Psi(P_{q_j})$ (in that case
we have to consider any $u_j,u_{j+1}$). Note that the case when $m=0$ is obviously true.
As usual, we write each polynomial $a_j(u_j)$ in the reduced form, that is we assume that 
it is a linear combination of products of the form 
$$
\wp(u_j)^{p}((\wp(u_j))^{*})^{q}
$$
where $p+q>0$. It suffices to show that any product of matrices 
$$
P_{q}Q_{q_0}a_1(u_1)Q_{q_1}a_2(u_2)\ldots Q_{m-1}a_m(u_m)Q_{q_m}P_{q}
$$ 
(we treat all factors as operator-valued matrices) is a linear combination of terms of the form
$w\otimes e(q,q)$, where $\varphi(w)=0$. This means that any product of the corresponding entries 
$\ell(p_j,q_j,u_j)$ and $\ell(p_j,q_j,u_j)^*$ arising in the computation of the
above moment does not reduce to a non-zero constant. If $u_j\neq u_{j+1}$,
such a reduction is not possible since the neighboring products of operators appearing at the first tensor site 
are free since they correspond
to different elements of $\mathpzc{U}$. In turn, if $u_j=u_{j+1}=u$, then we must have
$Q_{q_j}\neq 1$ in between $a_j(u)$ and $a_{j+1}(u)$. Then, the only terms in polynomials $a_{j}(u)$ and $a_{j+1}(u)$ 
which can produce non-zero constants at the first tensor site are such in which $Q_{q_j}$ is placed between $\wp(u)$* and $\wp(u)$
(the case when $m=1$ is easy and is omitted).
In order to examine this situation, we use the explicit form of $Q_{q_j}=Q_{p}$, where we set $q_j=p$ for simplicity,
$$
Q_{p}={\rm diag}(-d_{p}, \ldots , -d_p, 1-d_{p}, -d_p, \ldots , -d_p),
$$
with $1-d_p$ at the position $p$. An easy matrix multiplication gives 
$$
\wp(u)^{*}Q_p\wp(u)={\rm diag} (c_1, \ldots , c_r)
$$
where 
\begin{eqnarray*}
c_j
&=&
(1-d_p)\ell(p,j,u)^*\ell(p,j,u)-d_p\sum_{i\neq p}\ell(i,j,u)^{*}\ell(i,j,u)\\
&=&
(1-d_p)d_p-d_p\left(\sum_{i\neq p}d_{i}\right)=\;0
\end{eqnarray*}
and thus such terms also give zero contribution.
This proves that $\{\wp(u):u\in \mathpzc{U}\}$ is *-free from $\{P_1, \ldots , P_r\}$, which completes the proof.
\hfill $\blacksquare$\\

\begin{Remark}
{\rm Using Proposition 2.2, we can easily obtain the formulas
$$
\wp_{p,q}(u)=P_p\wp(u)P_q\;\;\;{\rm and}\;\;
\wp_{p,q}(u)^*=P_q \wp(u)^*P_p
$$
which leads to
$$ 
\omega_{p,q}(u)=P_p\wp(u)P_q+P_q\wp(u)^*P_p
$$
for any $p,q,u$, as well as similar formulas for the symmetrized Gaussian operators. If $\wp_{p,q}(u)$ has covariance $b_{p,q}(u)$, then 
we just need to rescale the above expressions. 
This enables us to express mixed moments of matricially free Gaussian operators and of their symmetrized counterparts 
under $\Psi$ in terms of mixed moments of $\omega(u)$ and $\{P_1, \ldots , P_r\}$. The same is true for
any $\Psi_q$ for which $d_{q}\neq 0$. } 
\end{Remark}

The arrays of laws of the operators $\omega_{p,q}(u)$ consist of semicircle and Bernoulli laws. Therefore, let us recall
their definitions. By the semicircle law of radius $2\alpha>0$ we understand the continuous distribution on the 
interval $[-2\alpha,2\alpha]$ with density
$$
d\sigma=\frac{\sqrt{4\alpha^{2}-x^{2}}}{2\pi \alpha^2}dx
$$
and the Cauchy transform
$$
G(z)=\frac{z-\sqrt{z^2-4\alpha^2}}{2\alpha^2},
$$
where the branch of the square root is chosen so that 
$\sqrt{z^2-4\alpha^2}>0$ if $z\in {\mathbb R}$ and $z\in (2\alpha, \infty)$. In turn,
by the Bernoulli law concentrated at $\pm \alpha$ we understand the discrete distribution 
$$
\kappa=\frac{\delta_{-\alpha}+\delta_{\alpha}}{2}
$$
with the Cauchy transform
$$
G(z)=\frac{1}{z-\alpha^2/z}.
$$
\begin{Proposition}
For any $(p,q)\in \mathpzc{J}$, let 
\begin{equation*}
{\mathcal A}_{p,q}={\mathbb C}\langle\wp_{p,q}, \wp^{*}_{p,q}, 1_{p,q}\rangle
\end{equation*}
and let $(\Psi_{p,q})$ be the array of states on $B({\mathcal M})$ defined by 
the states $\Psi_1, \ldots , \Psi_r$ associated with the vacuum vectors $\Omega_{1}, \ldots , \Omega_{r}$, 
respectively.
Then
\begin{enumerate}
\item
the array $({\mathcal A}_{p,q})$ is matricially free with respect to $(\Psi_{p,q})$, where the array of units 
is given by Definition 2.5,
\item
the $\Psi_{q,q}$-distribution of non-trivial ${\omega}_{q,q}(u)$ is the semicircle law of radius $2\sqrt{b_{q,q}(u)}$
for any $(q,q)\in \mathpzc{J}$ and $u\in \mathpzc{U}$,
\item
the $\Psi_{p,q}$-distribution of non-trivial ${\omega}_{p,q}(u)$ is the Bernoulli law concentrated at $\pm\sqrt{b_{p,q}(u)}$
for any off-diagonal $(p,q)\in \mathpzc{J}$ and $u\in \mathpzc{U}$.
\end{enumerate}
\end{Proposition}
{\it Proof.}
The proof reduces to that of [23, Proposition 4.2] since 
we can use Proposition 2.2 to get
$$
\wp_{p,q}= P_p\left(\sum_{u\in \mathpzc{U}}\wp(u)\right)P_q \;\;{\rm and}\;\;\wp_{p,q}^*=P_{q}\left(\sum_{u\in \mathpzc{U}}\wp(u)\right)^{*}P_p
$$
where $\{\wp(u):u\in \mathpzc{U}\}$ is *-free with respect to $\Psi$ and $\wp(u)^{*}\wp(v)=\delta_{u,v}$.
\hfill $\blacksquare$ \\

\begin{Corollary}
If $\mathpzc{J}=[r]\times [r]$, then the array of distributions of non-trivial $\omega_{p,q}(u)$, $u\in \mathpzc{U}$, 
in the states $(\Psi_{p,q})$ takes the matrix form 
$$
[\sigma(u)]=
\left(
\begin{array}{llll}
\sigma_{1,1}(u) & \kappa_{1,2}(u) & \ldots & \kappa_{1,r}(u)\\
\kappa_{2,1}(u) & \sigma_{2,2}(u) & \ldots & \kappa_{2,r}(u)\\
. & . & \ddots & . \\
\kappa_{r,1}(u) & \kappa_{r,2}(u) & \ldots & \sigma_{r,r}(u)
\end{array}
\right)
$$
where $\sigma_{q,q}(u)$ is the semicircle law of radius $2\sqrt{b_{q,q}(u)}$ 
for any $q$ and $\kappa_{p,q}(u)$ is the Bernoulli law concentrated at 
$\pm\sqrt{b_{p,q}(u)}$ for any $p\neq q$.
\end{Corollary}
{\it Proof.}
This is a special case of Proposition 2.3.
\hfill $\blacksquare$\\

Each matrix of the above form plays the role of a {\it matricial semicircle law}
(called {\it standard} if each diagonal law is the semicircle law of 
radius two and each off-diagonal law is the Bernoulli law concentrated 
at $\pm 1$).
In the case when some of the operators $\omega_{p,q}(u)$ are equal to zero, the corresponding 
measures are replaced by $\delta_{0}$. More interestingly, 
the family $\{\sigma(u), u\in \mathpzc{U}\}$ can be treated as a 
family of independent matricial semicircle laws.

\section{Symmetrized matricial semicircle laws}

The matricially free Gaussian operators are the ones which give the operatorial realizations of limit distributions 
for Hermitian Gaussian random matrices and their symmetric blocks in the most general case 
studied in this paper, including the situation when the dimension matrix $D$ is singular (contains zeros on the diagonal).

In particular, in the case when the dimension matrix $D$ has one zero on the diagonal and we evaluate 
the limit distribution of certain rectangular symmetric blocks under the associated partial trace, we obtain
the case studied by Benaych-Georges [6].
We will show that this general situation requires us to use non-symmetrized arrays of Gaussian operators.
Nevertheless, the special case when $D$ has no zeros on the diagonal is worth to be studied separately since in this case the 
limit distribution assumes an especially nice form. Then each $\omega_{p,q}(u)$ is non-trivial and, moreover, 
it always appears together with $\omega_{q,p}(u)$ 
whenever $p\neq q$. This leads to the following definition.

\begin{Definition}
{\rm If $\mathpzc{J}$ is symmetric, then the {\it symmetrized creation operators} 
are ope\-ra\-tors
of the form
\begin{equation*}
\widehat{\wp}_{p,q}(u)=
\left\{
\begin{array}{ll}
\wp_{q,q}(u) & {\rm if}\;p=q\\
\wp_{p,q}(u)+\wp_{q,p}(u)& {\rm if}\;p\neq q
\end{array}
\right.
\end{equation*} 
for any $(p,q)\in \mathpzc{J}$ and  $u\in \mathpzc{U}$. 
Their adjoints $\widehat{\wp}_{p,q}^{*}(u)$ will be called {\it symmetrized annihilation operators}. 
The {\it symmetrized Gaussian operators} are sums of the form 
\begin{equation*}
\widehat{\omega}_{p,q}(u)=\widehat{\wp}_{p,q}(u)+\widehat{\wp}_{p,q}^{*}(u)
\end{equation*} 
for any $(p,q)\in \mathpzc{J}$ and $u\in \mathpzc{U}$. 
Each matrix $(\widehat{\omega}_{p,q}(u))$ will be called a {\it symmetrized Gaussian matrix}.
Of course, the matrices of symmetrized operators are symmetric.}
\end{Definition}

We have shown in [23] that one array of symmetrized Gaussian operators 
gives the asymptotic joint distributions of symmetric random blocks of 
one Gaussian random matrix if the matrix $D$ is non-singular. 
In other words, symmetric blocks of a square Gaussian random matrix
behave as symmetrized Gaussian operators when the size of the matrix goes to infinity. 
We will generalize this result to the case when we have a family of matrices indexed by a finite set 
$\mathpzc{U}$ and we will allow the matrix $D$ to be singular.
The second situation leads to the following definition.

\begin{Definition}
{\rm If $\omega_{p,q}(u)\neq 0$ and $\omega_{q,p}(u)\neq 0$, then $\widehat{\omega}_{p,q}(u)$ will be called {\it balanced}.
If $\omega_{q,p}(u)=0$ and $\omega_{p,q}(u)\neq 0$ or $\omega_{p,q}(u)=0$ and $\omega_{q,p}(u)\neq 0$,
then $\widehat{\omega}_{p,q}(u)$ will be called {\it unbalanced}.
If $\omega_{p,q}(u)=\omega_{q,p}(u)=0$, then $\widehat{\omega}_{p,q}(u)=0$ will be called {\it trivial}.
}
\end{Definition}

\begin{Example}
{\rm In order to see the difference in computations of moments between the matricially free Gaussian operators and their symmetrizations, 
consider the following examples:
\begin{eqnarray*}
\Psi_1(\omega_{1,2}^{4})&=&0\\
\Psi_2(\omega_{1,2}^{4})&=& \Psi_{2}(\wp_{1,2}^{*}\wp_{1,2}^{}\wp_{1,2}^{*}\wp_{1,2}^{})=b_{1,2}^{2}\\
\Psi_1(\widehat{\omega}_{1,2}^{4}) &=&\Psi_{1}(\omega_{2,1}^{4})+\Psi_{1}(\omega_{2,1}^{}\omega_{1,2}^{2}\omega_{2,1}^{})\\
&=&\Psi_{1}(\wp_{2,1}^{*}\wp_{2,1}^{}\wp_{2,1}^{*}\wp_{2,1}^{})+\Psi_{1}(\wp_{2,1}^{*}\wp_{1,2}^{*}\wp_{1,2}^{}\wp_{2,1}^{})=b_{2,1}^{2}+b_{1,2}b_{2,1}\\
\Psi_2(\widehat{\omega}_{1,2}^{4})&=&\Psi_{2}(\omega_{1,2}^{4})+\Psi_{1}(\omega_{1,2}^{}\omega_{2,1}^{2}\omega_{1,2}^{})\\
&=&\Psi_{1}(\wp_{1,2}^{*}\wp_{1,2}^{}\wp_{1,2}^{*}\wp_{1,2}^{})+\Psi_{1}(\wp_{1,2}^{*}\wp_{2,1}^{*}\wp_{2,1}^{}\wp_{1,2}^{})=b_{1,2}^{2}+b_{1,2}b_{2,1}
\end{eqnarray*}
where, for simplicity, we took $\omega_{p,q}(u)=\omega_{p,q}$ for some fixed $u$ in all formulas. If, for instance, 
$\omega_{1,2}=0$, then the third moment reduces to $b_{2,1}^{2}$ and the remaining three moments vanish.
}
\end{Example}

In order to state the relation between the symmetrized operators, let us
define the {\it collective symmetrized units} in terms of the array 
$(1_{p,q})$ of Definition 2.5 as
\begin{equation*}
\widehat{1}_{p,q}=1_{p,q}+1_{q,p}-1_{p,q}1_{q,p}
\end{equation*}
respectively, for any $(p,q)\in \mathpzc{J}$. The symmetrized units are described in more detail 
in the propositions given below.

\begin{Proposition}
Let $\mathpzc{J}$ be symmetric and let ${\mathcal N}_{p,q}$ be defined as in Definition 2.5 for any $(p,q)\in \mathpzc{J}$.
\begin{enumerate}
\item
If $(q,q)\in {\mathpzc{J}}$, then $\widehat{1}_{q,q}$ is the orthogonal projection onto 
$$
{\mathbb C}\Omega_{q}\oplus \bigoplus_{k}\mathcal{N}_{q,k}.
$$
\item
If $(p,q)\in {\mathpzc{J}}$, where $p\neq q$, then $\widehat{1}_{p,q}$ is the orthogonal projection onto 
$$
{\mathbb C}\Omega_{p}\oplus {\mathbb C}\Omega_{q}\oplus \bigoplus_{k}(\mathcal{N}_{p,k}\oplus \mathcal{N}_{q,k}).
$$
\end{enumerate}
\end{Proposition}
{\it Proof.}
The above formulas follow directly from Definition 2.5.
\hfill $\blacksquare$\\

\begin{Proposition}
The collective symmetrized units satisfy the following relations:
\begin{eqnarray*}
\widehat{1}_{p,q}\widehat{1}_{k,l}&=&\left\{\begin{array}{ll}
0 & if\;\;\{p,q\}\cap \{k,l\}=\emptyset\\
\widehat{1}_{q,q} & if\;\;p\neq k\;and \;q=l\\
\widehat{1}_{p,q} & if\;\;(p,q)=(k,l)
\end{array}
\right.
\end{eqnarray*}
\end{Proposition}
{\it Proof.}
These relations are immediate consequences of Proposition 3.1.
\hfill $\blacksquare$\\

\begin{Proposition}
The following relations hold for any $u,t\in \mathpzc{U}$:
$$
\widehat{\wp}_{p,q}(u)\widehat{\wp}_{k,l}(t)=0\;\;and\;\;\widehat{\wp}_{p,q}^{*}(u)\widehat{\wp}_{k,l}^{*}(t)=0
$$ 
when $\{p,q\}\cap \{k,l\}=\emptyset$. Moreover, 
$$
\widehat{\wp}_{p,q}^{*}(u)\widehat{\wp}_{k,l}^{}(t)=0
$$ 
when $\{p,q\}\neq \{k,l\}$ or $u\neq t$. Finally, if the matrix $(b_{p,q}(u))$ is symmetric, then
$$
\widehat{\wp}_{p,q}^{*}(u)\widehat{\wp}_{p,q}^{}(u)=b_{p,q}(u)\widehat{1}_{p,q}
$$
for any $(i,j)\in \mathpzc{J}$. 
\end{Proposition}
{\it Proof.}
These relations follow from Proposition 2.1.
\hfill $\blacksquare$\\

In view of the last relation of Proposition 3.3, let us assume now that the matrices $(b_{i,j}(u))$ are symmetric for all $u\in \mathpzc{U}$.
We would like to define a symmetrized array of algebras related to the array $({\mathcal A}_{i,j})$ of Proposition 2.3. 
For that purpose, we introduce {\it collective symmetrized creation operators} 
\begin{equation*}
\widehat{\wp}_{p,q}=\sum_{u\in \mathpzc{U}}\widehat{\wp}_{p,q}(u)
\end{equation*}
and their adjoints $\widehat{\wp}_{p,q}^{*}$ for any $(p,q)\in\mathpzc{J}$. In turn, the sums
$$
\widehat{\omega}_{p,q}=\widehat{\wp}_{p,q}+\widehat{\wp}_{p,q}^{*}=\sum_{u\in \mathpzc{U}}\widehat{\omega}_{p,q}(u)
$$
will be called the {\it collective symmetrized Gaussian operators}. 

\begin {Proposition}
If $(b_{p,q}(u))$ is symmetric for any $u\in\mathpzc{U}$, then
\begin{equation*}
\widehat{\wp}_{p,q}^{*}\widehat{\wp}_{p,q}^{}=b_{p,q}\widehat{1}_{p,q},\;\;
\widehat{\wp}_{p,q}\widehat{\wp}_{k,l}=0\;\;and\;\;\widehat{\wp}_{p,q}^{*}\widehat{\wp}_{k,l}^{*}=0
\end{equation*}
whenever $\{p,q\}\cap \{k,l\}=\emptyset$, where $b_{p,q}=\sum_{u\in \mathpzc{U}}b_{p,q}(u)$.
\end{Proposition}
{\it Proof.}
The proof follows from Proposition 3.3. For instance,
$$
\widehat{\wp}_{p,q}^{*}\widehat{\wp}_{p,q}^{}=\sum_{u,t\in \mathpzc{U}}\widehat{\wp}_{p,q}^{*}(u)\widehat{\wp}_{p,q}^{}(t)=
\sum_{u\in \mathpzc{U}}b_{p,q}(u)\widehat{1}_{p,q}=b_{p,q}\widehat{1}_{p,q}
$$
The remaining relations are proved in a similar way.
\hfill $\blacksquare$\\

Let us discuss the notion of symmetric matricial freeness from an abstract point of view. 
We assume that the set $\mathpzc{J}$ is symmetric, by which we mean that $(j,i)\in \mathpzc{J}\;\;{\rm whenever} \;\;(i,j)\in \mathpzc{J}$.
Then we consider a symmetric array of subalgebras $({\mathcal A}_{i,j})$ 
of a unital algebra ${\mathcal A}$ and we assume that $(1_{i,j})$ is the associated symmetric array of internal units. 
By ${\mathcal I}$ we denote the unital algebra
generated by the internal units and we assume that it is commutative.
Moreover, the array of states $(\varphi_{i,j})$ on ${\mathcal A}$ remains as in the definition
of matricial freeness. Since the array $({\mathcal A}_{i,j})$ is symmetric,
we thus associate two states, $\varphi_{i,j}$ and $\varphi_{j,i}$, with each 
off-diagonal subalgebra ${\mathcal A}_{i,j}$.

Denote by $\widehat{\mathpzc{J}}$ the set of non-ordered pairs $\{i,j\}$, 
even in the case when $i=j$. Instead of sets $\mathpzc{J}_{\;m}$, we shall use their 
symmetric counterparts, namely subsets of the $m$-fold Cartesian product $\widehat{\mathpzc{J}}\times \widehat{\mathpzc{J}}\times \ldots \times \widehat{\mathpzc{J}}$ of the form
\begin{equation*}
\widehat{\mathpzc{J}}_{\;m}=\{(\{i_1,j_1\}, \ldots , \{i_m,j_m\}):
\{i_1,j_1\}\neq \ldots \neq \{i_m,j_m\}\},
\end{equation*}
and their subsets
\begin{equation*}
\widehat{\mathpzc{K}}_{\;m}=
\{(\{i_1,j_1\}, \ldots , \{i_m,j_m\}):
\{i_k,j_k\}\cap \{i_{k+1},j_{k+1}\} \neq \emptyset\;{\rm for}\; 1\leq k \leq m-1\}
\end{equation*}
where $m\in {\mathbb N}$. These sets comprise tuples of non-ordered pairs, 
in which neighboring pairs are different (as in the case of free products) and are related to each other
as in non-trivial multiplication of symmetric blocks of matrices.

In order to define the notion of {\it symmetric matricial freeness}, we shall first define the array of 
symmetrically matricially free units. This definition will differ from that given in [23]
since we have discovered that the conditions on the units need to be strengthened 
in order that the symmetrized array of Gaussian algebras be symmetrically matricially free as claimed in [23, Proposition 4.2].
In order to formulate stronger conditions on the moments involving symmetrized units, we also choose a slightly different 
formulation which contains relations between these units. This modifiation allows us to state the condition on the moments 
in the simplest possible form, involving only the diagonal units.

\begin{Definition}
{\rm 
Let $({\mathcal A}_{i,j})$ be a symmetric array of subalgebras of ${\mathcal A}$ 
with a symmetric array of internal units $(1_{i,j})$, and let $(\varphi_{i,j})$ be an
array of states on ${\mathcal A}$ defined by the family $(\varphi_{j})$. If, for some 
$a \in {\mathcal A}_{i,j}^{0}$, where $i\neq j$, it holds that
$$
\varphi_{j}(b1_{i,i}\,a)=\varphi_{j}(ba)\;\;{\rm and}\;\;\varphi_{j}(b1_{j,j}\,a)=0
$$
or
$$
\varphi_{j}(b1_{j,j}\,a)=\varphi_{j}(ba)\;\;{\rm and}\;\;\varphi_{q}(b1_{i,i}\,a)=0
$$
for any $b\in {\mathcal A}_{i,j}^{0}$, then we will say that $a$ is {\it odd} or {\it even}, respectively.
The subspaces of ${\mathcal A}_{i,j}^{0}$ spanned by even and odd elements will be called even and odd, respectively.
If each off-diagonal ${\mathcal A}_{i,j}^{0}$ is a direct sum of an odd subspace and an even subspace, the array
$({\mathcal A}_{i,j})$ will be called {\it decomposable}.}
\end{Definition}

\begin{Example}
{\rm The idea of even and odd elements comes from *-algebras generated by symmetrized creation operators. 
If $p\neq q$ and $k$ is odd, it is easy to see that
\begin{eqnarray*}
\Psi_{q}(b\widehat{1}_{p,p}^{}\widehat{\wp}_{p,q}^{\,k})&=&
\Psi_{q}(b\wp_{p,q}^{}\wp_{q,p}^{}\ldots \wp_{p,q}^{})\\
&=& \Psi_{q}(b\widehat{\wp}_{p,q}^{\,k})\\
\Psi_{q}(b\widehat{1}_{q,q}^{}\widehat{\wp}_{p,q}^{\,k})&=&0
\end{eqnarray*}
for any $b\in {\mathcal A}_{p,q}^{0}$ and thus $\widehat{\wp}_{p,q}^{\,k}$ is odd. In turn,
if $k$ is even,
\begin{eqnarray*}
\Psi_{q}(b\widehat{1}_{q,q}^{}\widehat{\wp}_{p,q}^{\,k})&=&
\Psi_{q}(b\wp_{q,p}^{}\wp_{p,q}^{}\ldots \wp_{p,q}^{})\\
&=&
\Psi_{q}(b\widehat{\wp}_{p,q}^{\,k})\\
\Psi_{q}(b\widehat{1}_{p,p}^{}\widehat{\wp}_{p,q}^{\,k})&=&0
\end{eqnarray*}
for any $b\in {\mathcal A}_{p,q}^{0}$
and thus $\widehat{\wp}_{p,q}^{\,k}$ is even. 
The main property of the pair $(\wp_{p,q}, \wp_{q,p})$
used here is that only alternating products of these two operators do not vanish and thus a non-trivial product of odd or even order
maps $\Omega_q$ onto a linear span of simple tensors which begin with $e_{p,q}(u)$, or $e_{q,p}(u)$, where $u\in \mathpzc{U}$, 
respectively.}
\end{Example}

\begin{Definition}
{\rm 
We say that the symmetric array $(1_{i,j})$ is a {\it symmetrically matricially free array of units} 
associated with a symmetric decomposable array $({\mathcal A}_{i,j})$ and the array $(\varphi_{i,j})$ 
if for any state $\varphi_j$ it holds that
\begin{enumerate}
\item
$\varphi_j(u_1au_2)=\varphi_j(u_1)\varphi_j(a)\varphi_j(u_2)$ 
for any $a\in {\mathcal A}$ and $u_1,u_2\in {\mathcal I}$,
\item
$\varphi_q(1_{k,l})=\delta_{q,k}+\delta_{q,l}-\delta_{q,k}\delta_{q,l}$ for any $q,k,l$,
\item
if $a_{r}\in {\mathcal A}_{i_r,j_r}\cap {\rm Ker}\varphi_{i_r,j_r}$, where $1<r\leq m$, then
$$
\varphi_j(a1_{i_1,j_1}a_{2}\ldots a_m)=
\left\{
\begin{array}{ll}
\varphi_j(aa_{2} \ldots a_m) & {\rm if}\;i_2=j_2\;\vee\;(i_2\neq j_2\;\wedge\;a_2\;{\rm odd})\\
0 & {\rm if} \;i_2\neq j_2\;\wedge\;a_2\;{\rm even}
\end{array}
\right.,
$$
where $a\in {\mathcal A}$ is arbitrary and $(\{i_1,j_1\}, \ldots , \{i_m,j_m\})\in \widehat{\mathpzc{K}}_{\;m}$, and
the moment also vanishes when $(\{i_1,j_1\}, \ldots , \{i_m,j_m\})\in \widehat{\mathpzc{J}}_{\;m}\setminus \widehat{\mathpzc{K}}_{\;m}$.
\end{enumerate}}
\end{Definition}

\begin{Definition}
{\rm We say that a symmetric decomposable array  $({\mathcal A}_{i,j})$ is 
{\it symmetrically matricially free} with respect to $({\varphi}_{i,j})$ if
\begin{enumerate}
\item for any $a_{k}\in {\rm Ker}\varphi_{i_k,j_k}\cap {\mathcal A}_{i_k,j_k}$, where $k\in [m]$
and $\{i_1,j_1\}\neq \ldots \neq \{i_m,j_m\}$, and for any state $\varphi_j$ it holds that
\begin{equation*}
\varphi_j(a_1a_2\ldots a_m)=0
\end{equation*}
\item
$(1_{i,j})$ is the associated symmetrically matricially free array of units.
\end{enumerate}}
\end{Definition} 

The array of variables $(a_{i,j})$ in a unital algebra ${\mathcal A}$ will be called 
{\it symmetrically matricially free} with respect to $(\varphi_{i,j})$ 
if there exists a symmetrically matricially free array of units 
$(1_{i,j})$ in ${\mathcal A}$ such that the array of algebras $({\mathcal A}_{i,j})$, 
each generated by $a_{i,j}+a_{j,i}$ and $1_{i,j}$, respectively, 
is symmetrically matricially free with respect to $(\varphi_{i,j})$.
The definition of *-symmetrically matricially free arrays of variables is similar 
to that of *-matricially free arrays.

\begin{Proposition}
Assume that $\mathpzc{J}$ and $(b_{p,q}(u))$ are symmetric
and $b_{p,q}(u)>0$ for any $(p,q)\in \mathpzc{J}$ and $u\in \mathpzc{U}$.  Let 
\begin{equation*}
\widehat{\mathcal A}_{p,q}={\mathbb C}\langle\widehat{\wp}_{p,q}, \widehat{\wp}^{*}_{p,q}\rangle
\end{equation*}
and let $(\Psi_{p,q})$ be the array of states on $B({\mathcal M})$ defined by 
the states $\Psi_1, \ldots , \Psi_r$ associated with the vacuum vectors $\Omega_{1}, \ldots , \Omega_{r}$, 
respectively. Then
\begin{enumerate}
\item
the array $(\widehat{\mathcal A}_{p,q})$ is symmetrically matricially free with respect to $(\Psi_{p,q})$, where 
the associated array of units is $(\widehat{1}_{p,q})$,
\item
the $\Psi_{p,q}$-distribution of $\,\widehat{\omega}_{p,q}(u)$ is the semicircle law of radius $2\sqrt{b_{p,q}(u)}$
for any $(p,q) \in {\mathpzc{J}}$ and $u\in \mathpzc{U}$.
\end{enumerate}
\end{Proposition}
{\it Proof.}
The proof of the second statement is the same as in the case when $\mathpzc{U}$ consists of one element [23, Proposition 8.1], so we will
only be concerned with the first one. 
In that connection, let us observe that 
$\widehat{1}_{p,q}\in \widehat{\mathcal A}_{p,q}$ for any $(p,q)\in \mathpzc{J}$ by Proposition 3.4, so that is why
we did not include the symmetrized units in the set of generators.
In order to prove symmetric matricial freeness of the array $(\widehat{\mathcal A}_{p,q})$, one needs to 
prove that the conditions of Definitions 3.4 and 3.5 are satisfied. 
If $p\neq q$, in view of commutation relations of Proposition 3.3,
any element $a\in\widehat{\mathcal A}_{p,q}^{0}:=\widehat{\mathcal A}_{p,q}\cap {\rm Ker}\Psi_{p,q}$
is spanned by products of the form
$$
\widehat{\wp}_{p_1,q_1}(u_1)\ldots \widehat{\wp}_{p_i,q_i}(u_i)
\widehat{\wp}_{k_1,l_1}^{*}(t_1)\ldots \widehat{\wp}_{k_j,l_j}^{*}(t_j),
$$
where $u_1, \ldots , u_i, t_1, \ldots, t_j\in \mathpzc{U}$ and $i+j>0$, with all pairs of matricial indices 
taken from the set $\{(p,q), (q,p)\}$ and alternating within both groups, i.e. 
$$
(p_{r+1},q_{r+1})=(q_r,p_r)\;{\rm and}\;(k_{r+1},l_{r+1})=(l_r,k_r)
$$ 
for any $r$.
Both facts allow us to conclude that any element 
of $\widehat{\mathcal A}_{p,q}^{0}$, when acting onto a simple tensor $w$,
either gives zero or
$$
e_{p_1,q_1}(u_{1})\otimes e_{p_2,q_2}(u_2) \ldots \otimes e_{p_i,q_i}(u_{i})\otimes w
$$
whenever $w$ begins with some $e_{q_i,k}(u)$ for some $k$ and $u$, 
where matricial indices are as above. 
The case when $p=q$ is identical to that treated in the proof of Proposition 2.3.
Therefore, a succesive action of elements $a_1, \ldots , a_n$ 
satisfying the assumptions of Definition 3.4 either 
gives zero or a vector which is orthogonal to ${\mathbb C}\Omega_q$, which implies that
$$
\Psi_{q}(a_{1}\ldots a_{n})=0
$$
for any $q\in \mathpzc{J}$, which proves the first condition of Definition 3.5. 
Let us prove now that $(\widehat{1}_{p,q})$  satisfies the conditions of Definition 3.4. 
The first condition holds since an analogous condition holds for the array $(1_{p,q})$. 
As for the second condition, it suffices to prove it for $p_1=q_1$ by Proposition 3.2. 
Let us first observe that each off-diagonal ${\mathcal A}_{p,q}^{0}$ 
is decomposable in the sense of Definition 3.3. Namely, the odd (even) subspace of ${\mathcal A}_{p,q}^{0}$ is spanned by products of
creation and annihilation operators of the form given above in which $i$ is odd (even).
Suppose first that $p_2\neq q_2$. Then, under the assumptions of Definition 3.3, 
it is easy to see that if $q_1=p_2$ and $a_{2}\in {\mathcal A}_{p,q}^{0}$ is 
odd, then either $a_{2}\ldots a_{m}\Omega_q=0$, or 
$$
a_2\ldots a_m\Omega_q=e_{p_2,q_2}(u_{1})\otimes \ldots \otimes e_{p_i,q_i}(u_i) \otimes w
$$
where $q_i=q_2$, for some $u_1, \ldots , u_i$ and some vector $w$, and
pairs $(q_1,q_2), \ldots , (p_i,q_i)$ are taken from the set 
$\{(p_2,q_2), (q_2,p_2)\}$ and alternate. By Proposition 3.1, the unit $\widehat{1}_{q_1,q_1}$ 
leaves the vector $a_{2}\ldots a_{m}\Omega_q$ invariant.
In turn, if $a_2$ is even, then either $a_{2}\ldots a_{m}\Omega_q=0$, or
$$
a_{2}\ldots a_{m}\Omega_q=e_{q_2,p_2}(u_{1})\otimes \ldots \otimes e_{p_i,q_i}(u_i) \otimes w
$$
where $q_i=q_2$ for some $u_1, \ldots , u_i$ and thus $\widehat{1}_{q_1,q_1}a_{2}\ldots a_{m}\Omega_q=0$.
Finally, if $p_2=q_2$ and $p_2\in \{p_1,q_1\}$, then it is obvious that 
$\widehat{1}_{p_1,q_1}$ leaves $a_2\ldots a_m\Omega_j$ invariant.
This completes the proof of symmetric matricial freeness.
\hfill $\blacksquare$ \\

\begin{Corollary}
Under the assumptions of Proposition 3.5, if $\mathpzc{J}=[r]\times [r]$, then the array of distributions of
$(\widehat{\omega}_{p,q}(u))$, $u\in \mathpzc{U}$, in the states $(\Psi_{p,q})$ is given by 
$$
[\widehat{\sigma}(u)]=
\left(
\begin{array}{llll}
\sigma_{1,1}(u) & \sigma_{1,2}(u) & \ldots & \sigma_{1,r}(u)\\
\sigma_{1,2}(u) & \sigma_{2,2}(u) & \ldots & \sigma_{2,r}(u)\\
. & . & \ddots & . \\
\sigma_{1,r}(u) & \sigma_{2,r}(u) & \ldots & \sigma_{r,r}(u)
\end{array}
\right)
$$
where $\sigma_{p,q}(u)$ is the semicircle law of radius $2\sqrt{b_{p,q}(u)}$ 
for any $(p,q)\in \mathpzc{J}$ and $u\in \mathpzc{U}$. 
\end{Corollary}

Each symmetric matrix of the above form plays the role of a {\it symmetrized matricial semicircle law}.
It clearly differs from the matricial semicircle law of Corollary 2.1, but 
it contains all information needed to compute all mixed moments of the symmetrized Gaussian operators 
in the states $\Psi_1, \ldots , \Psi_r$. For that reason, we use square matrices since we need two 
distributions of each off-diagonal operator $\widehat{\omega}_{p,q}(u)$, in the states $\Psi_p$ and $\Psi_q$. 

Even if the matrices $B(u)=(b_{p,q}(u))$ are not symmetric, a result similar to Proposition 3.5 can 
be proved. For that purpose, we shall need the probability measure $\vartheta$ on ${\mathbb R}$ corresponding to the two-periodic 
continued fraction $(a,b,a,\ldots )$. Its Cauchy transform is 
$$
G(z)=\frac{z^2+b-a-\sqrt{(z^2-b-a)^{2}-4ab}}{2zb}.
$$
It has the absolutely continuous part 
$$
d\vartheta=\frac{\sqrt{4ab-(x^2-a-b)^2}}{2\pi bx}dx
$$
supported on $|\sqrt{a}-\sqrt{b}|\leq |\,x|\leq \sqrt{a}+\sqrt{b}$
and an atom of mass $1/2-a/2b$ at $x=0$ if $a\neq b$ (see, for instance, [21, Example 9.2]). 
In particular, if $a=b$, then $\vartheta$ is the semicircle law of radius $2a$.

\begin{Proposition}
Assume that $\mathpzc{J}$ is symmetric and that 
$b_{p,q}(u)>0$ for any $(p,q)\in \mathpzc{J}$ and $u\in \mathpzc{U}$.  Let 
\begin{equation*}
\widehat{\mathcal A}_{p,q}={\mathbb C}\langle\widehat{\wp}_{p,q}, \widehat{\wp}^{*}_{p,q}, r_{p,q}, r_{q,p}\rangle
\end{equation*}
and let $(\Psi_{p,q})$ be the array of states on $B({\mathcal M})$ defined by 
the states $\Psi_1, \ldots , \Psi_r$ associated with the vacuum vectors $\Omega_{1}, \ldots , \Omega_{r}$, 
respectively. Then
\begin{enumerate}
\item
the array $(\widehat{\mathcal A}_{p,q})$ is symmetrically matricially free with respect to $(\Psi_{p,q})$, where 
the associated array of units is $(\widehat{1}_{p,q})$,
\item
the $\Psi_{p,q}$-distribution of $\,\widehat{\omega}_{p,q}(u)$ is $\vartheta_{p,q}(u)$ 
for any $(p,q) \in {\mathpzc{J}}$ and $u\in \mathpzc{U}$, where $\vartheta_{p,q}(u)$
is the distribution corresponding to the two-periodic continued fraction with the Jacobi sequence $(b_{p,q}(u), b_{q,p}(u), \ldots )$.
\end{enumerate}
\end{Proposition}
{\it Proof.}
The proof of matricial freeness is similar to that of Proposition 3.5, but more technical due to 
more complicated commutation relations between the off-diagonal collective creation and annihilation operators
$$
\widehat{\wp}_{p,q}^{*}\widehat{\wp}_{p,q}=b_{p,q}r_{p,q}+b_{q,p}r_{q,p}
$$
where $b_{p,q}=\sum_{u\in \mathpzc{U}}b_{p,q}(u)$ and $b_{q,p}=\sum_{u\in \mathpzc{U}}b_{q,p}$
for any $p\neq q$. In turn, we will show in Example 4.4 that the corresponding probability measures on the real line 
correspond to two-periodic continuous fractions with sequences of Jacobi coefficients 
$(b_{q,p}(u), b_{p,q}(u), \ldots)$ and $(b_{p,q}(u), b_{q,p}(u), \ldots )$, respectively.
\hfill 
$\blacksquare$\\

\begin{Corollary}
Under the assumptions of Proposition 3.6, if $\mathpzc{J}=[r]\times [r]$, the array of distributions of
$(\widehat{\omega}_{p,q}(u))$, $u\in \mathpzc{U}$, in the states $(\Psi_{p,q})$ is given by 
$$
[\widehat{\sigma}(u)]=
\left(
\begin{array}{llll}
\sigma_{1,1}(u) & \vartheta_{1,2}(u) & \ldots & \vartheta_{1,r}(u)\\
\vartheta_{2,1}(u) & \sigma_{2,2}(u) & \ldots & \vartheta_{2,r}(u)\\
. & . & \ddots & . \\
\vartheta_{r,1}(u) & \vartheta_{r,2}(u) & \ldots & \sigma_{r,r}(u)
\end{array}
\right)
$$
where and $\vartheta_{p,q}(u)$ is the distribution corresponding to
the two-periodic continued fraction with the Jacobi sequence $(b_{p,q}(u), b_{q,p}(u), \ldots )$
for any $(p,q)\in \mathpzc{J}$ and $u\in \mathpzc{U}$. 
\end{Corollary}

\section{Combinatorics of mixed moments}

It can be seen from the computations in Section 2 that the mixed moments of matricially free Gaussian operators are related to
noncrossing pair partitions. The difference between them and free Gaussians is that one has to use colored 
noncrossing pair partitions, where coloring is adapted to the pairs of matricial indices from $\mathpzc{J}$ 
and to the additional indices from the set $\mathpzc{U}$. It will be convenient to
assume here that $\mathpzc{U}=[t]$, where $t$ is a natural number.

For a given noncrossing pair partition $\pi$,
we denote by $\mathcal{B}(\pi)$, $\mathcal{L}(\pi)$ and $\mathcal{R}(\pi)$ the sets of its blocks, 
their left and right legs, respectively. If $\pi_{i}=\{l(i),r(i)\}$ and $\pi_{j}=\{l(j),r(j)\}$
are blocks of $\pi$ with left legs $l(i)$ and $l(j)$ and right legs $r(i)$ and $r(j)$, respectively, then
$\pi_i$ is {\it inner} with respect to $\pi_j$ if $l(j)<l(i)<r(i)<r(j)$.
In that case $\pi_j$ is {\it outer} with respect to $\pi_i$.
It is the {\it nearest outer block} of $\pi_i$ if there is no block $\pi_k=\{l(k),r(k)\}$
such that $l(j)<l(k)<l(i)<r(i)<r(k)<r(j)$. Since the nearest outer block, if it exists, is unique,
we can write in this case 
$$
\pi_j=o(\pi_i), \;l(j)=o(l(i))\;{\rm and}\;r(j)=o(r(i)).
$$
If $\pi_i$ does not have an outer block, it is called a {\it covering} block.
In that case we set $o(\pi_i)=\pi_0$, where we define $\pi_0=\{0, m+1\}$
and call the {\it imaginary block}.

\begin{Example}
{\rm Let $\sigma\in \mathcal{NC}_{6}^{2}$ be as in Fig. 1.
Its blocks are $\{1,6\}, \{2,3\}, \{4,5\}$ and the imaginary block is $\{0,7\}$.
The left and right legs of $\pi$ are $\mathcal{L}(\sigma)=\{1,2,4\}$ and $\mathcal{R}(\sigma)=\{3,5,6\}$.
The block $\{1,6\}$ is the nearest outer block of both $\{2,3\}$ and $\{4,5\}$ and the imaginary block $\{0,7\}$
is the nearest outer block of $\{1,6\}$.}
\end{Example}

Computations of mixed moments of $(\omega_{p,q}(u))$ in the states $\Psi_q$ are based on the classes
of colored noncrossing pair partitions adapted to ordered tuples of indices. Since, in addition to 
the pair of matricial indices $(p,q)$, we have an additional index $u$ as compared with the case studied previously, 
we distinguish the matricial indices from the non-matricial one and we use an abbreviated notation
for the pair $\mathpzc{v}=(p,q)$. In the random matrix context, the matricial indices will be related 
to blocks of a random matrix, whereas $u$ will label independent matrices. 
If the set $\mathpzc{U}$ consists of one element, we recover the analogous 
definition of [23].

\begin{Definition}
{\rm We will say that $\pi \in \mathcal{NC}_{m}^{2}$ is {\it adapted}
to the tuple $((\mathpzc{v}_1,u_1), \ldots, (\mathpzc{v}_{m},u_m))$, where 
$\mathpzc{v}_{k}=(p_k,q_k)\in [r]\times [r]$ and $u_k\in [t]$ for any $k$, if
\begin{enumerate}
\item[(a)]
$(\mathpzc{v}_{i},u_i)=(\mathpzc{v}_j,u_j)$ whenever $\{i,j\}$ is a block of $\pi$,
\item[(b)]
$q_j=p_{o(j)}$ whenever $\{i,j\}$ is a block of $\pi$ which has an outer block.
\end{enumerate}
The set of such partitions will be denoted by 
$\mathcal{NC}_{m}^{2}((\mathpzc{v}_1,u_1), \ldots , (\mathpzc{v}_{m},u_m))$.
In turn, by $\mathcal{NC}_{m,q}^{2}((\mathpzc{v}_1,u_1), \ldots , (\mathpzc{v}_{m},u_m))$ 
we will denote its subset, for which $q_{k}=q$ whenever $k$ belongs to a covering block.}
\end{Definition}

In order to find combinatorial formulas for the moments of our Gaussian operators, we need
to use colored noncrossing pair partitons. It will suffice to color each $\pi\in \mathcal{NC}_{m}^{2}$, where $m$ is even,
by numbers from the set $[r]$ and label it by numbers from the set $\mathpzc{U}$. 
We will denote by $F_{r}(\pi)$ the set of all mappings $f$ from the set of blocks of $\pi$ 
into $[r]$ called {\it colorings}. In turn, we will denote by $F_{t}(\pi)$ the set of 
all mappings $g$ from the set of blocks of $\pi$ into $[t]$ called {\it labelings}. 

Thus, by a {\it colored labeled noncrossing pair partition} (or, simply, {\it colored noncrossing pair partition}) 
we will understand a triple $(\pi,f,g)$, where $\pi \in \mathcal{NC}_{m}^{2}$, $f\in F_{r}(\pi)$ and $g\in F_{t}(\pi)$. The set 
$$
{\mathcal B}(\pi,f,g)=\{(\pi_1,f,g),  \ldots , (\pi_k,f,g)\}
$$ 
will denote the set of its blocks. We will always assume that also the imaginary block is colored by a number 
from the set $[r]$, but no labeling of the imaginary block is needed). If all blocks are labeled by the same index, 
we will omit both this index and the function $g$ in our notations. 

\begin{Example}
{\rm If $\pi\in\mathcal{NC}_{m,q}^{2}((\mathpzc{v}_{1},u_{1}),\ldots, (\mathpzc{v}_{m},u_m))$, 
then $(\mathpzc{v}_{1},\ldots, \mathpzc{v}_{m})$ defines a unique 
coloring of $\pi$ and $\widehat{\pi}$ in which the block containing $k$ is colored by $p_k$ for any $k$ 
and the imaginary block is colored by $q$. Similarly, $(u_{1},\ldots, u_m)$ 
defines a unique natural labeling of $\pi$. Consider the first colored partition given in 
Fig. 1, where, for simplicity, we assume that $t=1$ and thus we can skip $u$'s.
If we are given the tuple of pairs 
$$
(\mathpzc{v}_1, \mathpzc{v}_2, \mathpzc{v}_3, \mathpzc{v}_4)=((2,1),(2,2),(2,2),(2,1))
$$
then $\pi$ is adapted to this tuple. If $q=1$, then the unique coloring of $\pi$ defined by the
given tuple and the number $q$ is given by $f_1$ since the imaginary block must be colored by $1$ and the
colors of blocks $\{1,4\}$ and $\{2,3\}$ are obtained from the numbers assigned to their left legs, i.e.
$f_1(1)=2$ and $f_1(2)=2$. }
\end{Example}

\begin{Definition}
{\rm Let a real-valued matrix $B(u)=(b_{i,j}(u))\in M_{r}({\mathbb R})$ be given for any $u\in [t]$.
We define a family of real-valued functions $b_q$, where $q\in [r]$, on the set of colored noncrossing pair-partitions by
$$
b_{q}(\pi, f,g)=b_{q}(\pi_{1},f,g)
\ldots b_{q}(\pi_k,f,g)
$$
where $\pi\in \mathcal{NC}_{m}^{2}$, $f\in F_{r}(\pi)$, $g\in F_{t}(\pi)$ 
and $b_q$ is defined on the set of blocks ${\mathcal B}(\pi,f,g)$ as
$$
b_{q}(\pi_k,f,g)=b_{i,j}(u),
$$
whenever block $\pi_k$ is colored by $i$ and labeled by $u$, its nearest outer block is colored by $j$, with the imaginary block 
colored by $q$.}
\end{Definition}

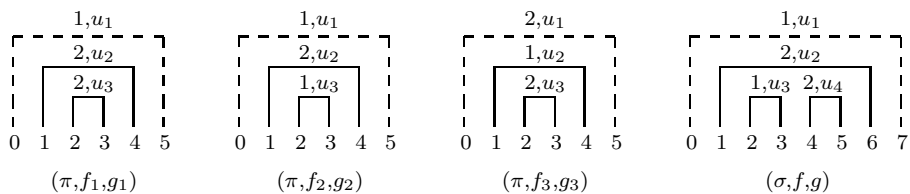
\begin{figure}
\unitlength=1mm
\special{em:linewidth 0.4pt}
\linethickness{0.4pt}
\begin{picture}(160.00,35.00)(16.00,0.00)
\put(39.00,10.00){\line(0,1){8.00}}
\put(43.00,10.00){\line(0,1){4.00}}
\put(47.00,10.00){\line(0,1){4.00}}
\put(51.00,10.00){\line(0,1){8.00}}

\put(35.00,10.00){\line(0,1){1.25}}
\put(35.00,12.50){\line(0,1){1.25}}
\put(35.00,13.00){\line(0,1){1.25}}
\put(35.00,15.50){\line(0,1){1.25}}
\put(35.00,18.00){\line(0,1){1.25}}
\put(35.00,20.75){\line(0,1){1.25}}

\put(55.00,10.00){\line(0,1){1.25}}
\put(55.00,12.50){\line(0,1){1.25}}
\put(55.00,13.00){\line(0,1){1.25}}
\put(55.00,15.50){\line(0,1){1.25}}
\put(55.00,18.00){\line(0,1){1.25}}
\put(55.00,20.60){\line(0,1){1.25}}

\put(34.50,7.00){$\scriptstyle{0}$}
\put(38.50,7.00){$\scriptstyle{1}$}
\put(42.50,7.00){$\scriptstyle{2}$}
\put(46.50,7.00){$\scriptstyle{3}$}
\put(50.50,7.00){$\scriptstyle{4}$}
\put(54.50,7.00){$\scriptstyle{5}$}

\put(43.00,15.00){$\scriptstyle{2,u_3}$}
\put(43.00,19.00){$\scriptstyle{2,u_2}$}
\put(43.00,23.50){$\scriptstyle{1,u_1}$}

\put(35.00,22.00){\line(1,0){1.25}}

\put(37.67,22.00){\line(1,0){1.25}}
\put(40.33,22.00){\line(1,0){1.25}}
\put(42.98,22.00){\line(1,0){1.25}}
\put(45.63,22.00){\line(1,0){1.25}}
\put(48.30,22.00){\line(1,0){1.25}}
\put(50.95,22.00){\line(1,0){1.25}}
\put(53.75,22.00){\line(1,0){1.25}}

\put(39.00,18.00){\line(1,0){12.00}}
\put(43.00,14.00){\line(1,0){4.00}}


\put(69.00,10.00){\line(0,1){8.00}}
\put(73.00,10.00){\line(0,1){4.00}}
\put(77.00,10.00){\line(0,1){4.00}}
\put(81.00,10.00){\line(0,1){8.00}}

\put(65.00,10.00){\line(0,1){1.25}}
\put(65.00,12.50){\line(0,1){1.25}}
\put(65.00,13.00){\line(0,1){1.25}}
\put(65.00,15.50){\line(0,1){1.25}}
\put(65.00,18.00){\line(0,1){1.25}}
\put(65.00,20.75){\line(0,1){1.25}}

\put(85.00,10.00){\line(0,1){1.25}}
\put(85.00,12.50){\line(0,1){1.25}}
\put(85.00,13.00){\line(0,1){1.25}}
\put(85.00,15.50){\line(0,1){1.25}}
\put(85.00,18.00){\line(0,1){1.25}}
\put(85.00,20.60){\line(0,1){1.25}}

\put(64.50,7.00){$\scriptstyle{0}$}
\put(68.50,7.00){$\scriptstyle{1}$}
\put(72.50,7.00){$\scriptstyle{2}$}
\put(76.50,7.00){$\scriptstyle{3}$}
\put(80.50,7.00){$\scriptstyle{4}$}
\put(84.50,7.00){$\scriptstyle{5}$}

\put(73.00,15.00){$\scriptstyle{1,u_3}$}
\put(73.00,19.00){$\scriptstyle{2,u_2}$}
\put(73.00,23.50){$\scriptstyle{1,u_1}$}

\put(65.00,22.00){\line(1,0){1.25}}
\put(67.67,22.00){\line(1,0){1.25}}
\put(70.33,22.00){\line(1,0){1.25}}
\put(72.98,22.00){\line(1,0){1.25}}
\put(75.63,22.00){\line(1,0){1.25}}
\put(78.30,22.00){\line(1,0){1.25}}
\put(80.95,22.00){\line(1,0){1.25}}
\put(83.75,22.00){\line(1,0){1.25}}

\put(69.00,18.00){\line(1,0){12.00}}
\put(73.00,14.00){\line(1,0){4.00}}


\put(99.00,10.00){\line(0,1){8.00}}
\put(103.00,10.00){\line(0,1){4.00}}
\put(107.00,10.00){\line(0,1){4.00}}
\put(111.00,10.00){\line(0,1){8.00}}

\put(95.00,10.00){\line(0,1){1.25}}
\put(95.00,12.50){\line(0,1){1.25}}
\put(95.00,13.00){\line(0,1){1.25}}
\put(95.00,15.50){\line(0,1){1.25}}
\put(95.00,18.00){\line(0,1){1.25}}
\put(95.00,20.75){\line(0,1){1.25}}

\put(115.00,10.00){\line(0,1){1.25}}
\put(115.00,12.50){\line(0,1){1.25}}
\put(115.00,13.00){\line(0,1){1.25}}
\put(115.00,15.50){\line(0,1){1.25}}
\put(115.00,18.00){\line(0,1){1.25}}
\put(115.00,20.60){\line(0,1){1.25}}

\put(94.50,7.00){$\scriptstyle{0}$}
\put(98.50,7.00){$\scriptstyle{1}$}
\put(102.50,7.00){$\scriptstyle{2}$}
\put(106.50,7.00){$\scriptstyle{3}$}
\put(110.50,7.00){$\scriptstyle{4}$}
\put(114.50,7.00){$\scriptstyle{5}$}

\put(103.00,15.00){$\scriptstyle{2,u_3}$}
\put(103.00,19.00){$\scriptstyle{1,u_2}$}
\put(103.00,23.50){$\scriptstyle{2,u_1}$}

\put(95.00,22.00){\line(1,0){1.25}}
\put(97.67,22.00){\line(1,0){1.25}}
\put(100.33,22.00){\line(1,0){1.25}}
\put(102.98,22.00){\line(1,0){1.25}}
\put(105.63,22.00){\line(1,0){1.25}}
\put(108.30,22.00){\line(1,0){1.25}}
\put(110.95,22.00){\line(1,0){1.25}}
\put(113.75,22.00){\line(1,0){1.25}}

\put(99.00,18.00){\line(1,0){12.00}}
\put(103.00,14.00){\line(1,0){4.00}}

\put(129.00,10.00){\line(0,1){8.00}}
\put(133.00,10.00){\line(0,1){4.00}}
\put(137.00,10.00){\line(0,1){4.00}}
\put(141.00,10.00){\line(0,1){4.00}}
\put(145.00,10.00){\line(0,1){4.00}}
\put(149.00,10.00){\line(0,1){8.00}}

\put(125.00,10.00){\line(0,1){1.25}}
\put(125.00,12.50){\line(0,1){1.25}}
\put(125.00,13.00){\line(0,1){1.25}}
\put(125.00,15.50){\line(0,1){1.25}}
\put(125.00,18.00){\line(0,1){1.25}}
\put(125.00,20.75){\line(0,1){1.25}}

\put(153.00,10.00){\line(0,1){1.25}}
\put(153.00,12.50){\line(0,1){1.25}}
\put(153.00,13.00){\line(0,1){1.25}}
\put(153.00,15.50){\line(0,1){1.25}}
\put(153.00,18.00){\line(0,1){1.25}}
\put(153.00,20.60){\line(0,1){1.25}}

\put(124.50,7.00){$\scriptstyle{0}$}
\put(128.50,7.00){$\scriptstyle{1}$}
\put(132.50,7.00){$\scriptstyle{2}$}
\put(136.50,7.00){$\scriptstyle{3}$}
\put(140.50,7.00){$\scriptstyle{4}$}
\put(144.50,7.00){$\scriptstyle{5}$}
\put(148.50,7.00){$\scriptstyle{6}$}
\put(152.50,7.00){$\scriptstyle{7}$}

\put(133.00,15.00){$\scriptstyle{1,u_3}$}
\put(140.00,15.00){$\scriptstyle{2,u_4}$}
\put(137.00,19.00){$\scriptstyle{2,u_2}$}
\put(137.00,23.50){$\scriptstyle{1,u_1}$}

\put(125.00,22.00){\line(1,0){1.25}}
\put(127.67,22.00){\line(1,0){1.25}}
\put(130.33,22.00){\line(1,0){1.25}}
\put(132.98,22.00){\line(1,0){1.25}}
\put(135.63,22.00){\line(1,0){1.25}}
\put(138.30,22.00){\line(1,0){1.25}}
\put(140.95,22.00){\line(1,0){1.25}}
\put(143.60,22.00){\line(1,0){1.25}}
\put(146.25,22.00){\line(1,0){1.25}}
\put(149.00,22.00){\line(1,0){1.25}}
\put(151.70,22.00){\line(1,0){1.30}}

\put(129.00,18.00){\line(1,0){20.00}}
\put(133.00,14.00){\line(1,0){4.00}}
\put(141.00,14.00){\line(1,0){4.00}}
\put(40.00,2.00){$\scriptstyle{(\pi,f_1,g_1)}$}
\put(70.00,2.00){$\scriptstyle{(\pi,f_2,g_2)}$}
\put(100.00,2.00){$\scriptstyle{(\pi,f_3,g_3)}$}
\put(135.00,2.00){$\scriptstyle{(\sigma,f,g)}$}

\end{picture}
\caption{Colored labeled noncrossing pair partitions}
\end{figure}

It should be remarked that in this paper $b_{q}(\pi,f,g)$ may be equal to zero 
even if $\pi \in \mathcal{NC}_{m}^{2}((\mathpzc{v}_1,u_1), \ldots , (\mathpzc{v}_{m},u_m))$
since we assume that the matrices $B(u)$ may contain zeros. 
Let us also recall our convention saying that if $\mathcal{NC}_{m}^{2}$ or its subset is empty, we shall understand that the summation 
over $\pi\in \mathcal{NC}_{m}^{2}$ or over this subset gives zero. In particular, this will always be the case
when $m$ is odd. 

\begin{Proposition}
For any tuple $((\mathpzc{v}_{1},u_1), \ldots , (\mathpzc{v}_m,u_m))$ and $q\in [r]$, $m\in {\mathbb N}$,
where $\mathpzc{v}_{k}=(p_k,q_k)\in [r]\times [r]$ and $u_k\in [t]$ for each $k$, it holds that
$$
\Psi_{q}\left(\omega_{p_1,q_1}(u_1)\ldots \omega_{p_m,q_m}(u_m)\right)
=
\sum_{\pi\in \mathcal{NC}^{2}_{m,q}((\mathpzc{v}_{1},u_1), \ldots , (\mathpzc{v}_m,u_m))}b_q(\pi,f,g)
$$
where $f$ is the coloring of $\pi$ defined by $(\mathpzc{v}_{1}, \ldots , \mathpzc{v}_m;q)$ 
and $g$ is the labeling defined  by $(u_1, \ldots , u_m)$.
\end{Proposition}
{\it Proof.}
The proof is similar to that of [23, Lemma 5.2] and it reduces to showing that
if $\epsilon_1, \ldots , \epsilon_m\in \{1,*\}$ and
$\pi\in \mathcal{NC}_{m,q}^{2}((\mathpzc{v}_{1},u_1), \ldots , (\mathpzc{v}_{m},u_m))$, where $m$ is even, then
$$
\Psi_{q}(\wp_{p_1,q_1}^{\epsilon_{1}}(u_1)\ldots \wp_{p_m,q_m}^{\epsilon_m}(u_m))
=
b_{q}(\pi,f,g)
$$
where $f$ is the coloring of $\pi$ defined by the collection of indices $\{p_{k}, k\in \mathcal{L}(\pi)\}$
associated with the left legs of the blocks of $\pi$ and the index $q$ coloring the imaginary block, and
$g$ is the natural labeling of $\pi$. 
The only difference between the former proof and this one
is that to each block of $(\pi,f,g)$ we assign a matrix element of $B(u)$ for 
suitable $u$ ($u$ is the same for both legs of each block since the partition satisfies 
condition (a) of Definition 4.1. \hfill $\blacksquare$\\

Using Proposition 4.1, we can derive nice combinatorial formulas for the moments of sums of 
collective Gaussian operators
$$
\omega=\sum_{p,q}\omega_{p,q}, 
$$
where
$$
\omega_{p,q}=\sum_{u\in \mathpzc{U}}\omega_{p,q}(u),
$$ 
and for that purpose, denote by $\mathcal{NC}_{m,q}^{\,2}[r]$ the set of all noncrossing 
pair partitions of $[m]$ colored by the set $[r]$, with the imaginary block colored by $q$.
\begin{Lemma}
The moments of $\omega$ in the state $\Psi_q$, where $q\in [r]$, are given by
\begin{eqnarray*}
\Psi_{q}(\omega^m)&=&\sum_{(\pi,f)\in \mathcal{NC}_{m,q}^{\,2}[r]}b_{q}(\pi,f)
\end{eqnarray*}
where $b_{q}(\pi,f)=b_{q}(\pi_1,f)\ldots b_{q}(\pi_s,f)$ for $\pi=\{\pi_1, \ldots, \pi_s\}$ 
and 
$$
b_{q}(\pi_k,f)=\sum_{u}b_{i,j}(u)
$$ 
whenever block $\pi_k$ is colored by $i$ and its nearest outer block is colored by $j$.
\end{Lemma}
{\it Proof.}
Using Proposition 4.1, we can express each summand in the formula
$$
\Psi_{q}(\omega^{m})=\sum_{p_1,q_1,u_1, \ldots, p_m,q_m,u_m}\Psi_{q}(\omega_{p_1,q_1}(u_1)\ldots \omega_{p_m,q_m}(u_m))
$$ 
in terms of $b_{q}(\pi,f,g)$, where $\pi\in \mathcal{NC}_{m,q}^{\,2}((\mathpzc{v}_1,u_1), \ldots , (\mathpzc{v}_{m},u_m))$, with
$\mathpzc{v}_{k}=(p_k,q_k)$ and $f$ is the coloring defined by the tuple 
$(\mathpzc{v}_1, \ldots , \mathpzc{v}_{m})$ whereas $g$ is the labeling defined by $(u_1, \ldots , u_m)$. 
It is clear that if $(\mathpzc{v}_1, \ldots , \mathpzc{v}_{m})$ is fixed and 
$\pi$ is adapted to it, then all labelings to which $\pi$ is adapted are pairwise identical within blocks
which results the addition formula for $b_{q}(\pi_k,f)$.
Now, a different choice of $(\mathpzc{v}_1, \ldots , \mathpzc{v}_{m})$ must lead to a different 
collection of pairs $(\pi,f)$ for given $g$ since even if the same partition appears on the RHS 
of the formula of Proposition 4.1, the coloring $f$ must be different. In fact, 
$\pi$ is uniquely determined by the sequence $\epsilon=(\epsilon_1, \ldots , \epsilon_m)$
which appears in the nonvanishing moment of type
$$
\Psi_q(\wp_{p_1,q_1}^{\epsilon_1}(u_1)\ldots \wp_{p_m,q_m}^{\epsilon_m}(u_m)),
$$
where $\epsilon_k\in \{1,*\}$. Moments of this type are basic 
constituents of $\Psi_{q}(\omega^{m})$. For instance, if $\epsilon=(*,*,1,*,1,1)$, the associated 
unique noncrossing pair partition is $\sigma$ of Fig.1. If we keep $\pi$ and change at least one matricial index in the given tuple to which
$\pi$ is adapted, we either get a tuple to which $\pi$ is not adapted (if we change $q_k$ for some $k\in \mathcal{L}(\pi)$ since this index is determined by the color of its nearest outer block) or we obtain a different coloring $f$ (if we change $p_k$ for some $k\in \mathcal{L}(\pi)$).  
Therefore, all contributions $b_{q}(\pi,f,g)$ of the form given by Proposition 4.1 associated with  
different tuples of indices correspond to different elements of $\mathcal{NC}_{m,q}^{2}[r]$. Therefore, in order to prove the formula for $\Psi_q(\omega^{m})$,
we only need to justify that all colored partitions from $\mathcal{NC}_{m,q}^{2}[r]$ really do contribute to $\Psi_{q}(\omega^{m})$.
Thus, let $(\pi,f)\in \mathcal{NC}_{m,q}^{2}[r]$ be given. There exists a unique $\epsilon=(\epsilon_1, \ldots , \epsilon_m)$ associated with $\pi$.
In turn, $f$ determines $((p_1,q_1), \ldots, (p_m,q_m))$ by an inductive procedure with respect to blocks' depth. Thus,
if $\{k,l\}$ is a covering block, we choose $q_k=q_l=q$ and $p_k=p_l=f(k)=f(l)$. Next, if $\{k,l\}=o(\{i,j\})$, then we choose
$q_i=q_j=p_k$ and $p_i=p_j=f(i)=f(j)$, etc. We proceed in this fashion until we choose all $p_k,q_k$, $k\in [m]$. 
This completes the proof.
\hfill $\blacksquare$\\

Note that $\omega$ coincides with the sum of all symmetrized Gaussian operators 
$$
\widehat{\omega}=\sum_{p\leq q}\widehat{\omega}_{p,q}
$$ 
and thus Lemma 4.1 gives a formula for the 
moments of $\widehat{\omega}$ as well. However, this is not the case for the mixed moments
of symmetrized Gaussian operators $(\widehat{{\omega}}_{p,q}(u))$ in the state $\Psi_q$. 
These are based on the class of colored noncrossing pair partitions adapted to ordered tuples of indices of type 
$(\mathpzc{w},u)$, where $\mathpzc{w}$ is an abbreviated notation for the set $\{p,q\}$.

The definition of this class is given below. Note that it is slightly stronger than that in [23], which is a consequence of
the stronger definition of symmetric matricial freeness discussed in Section 3.

\begin{Definition}
{\rm We say that $\pi \in \mathcal{NC}_{m}^{2}$ is {\it adapted}
to the tuple $((\mathpzc{w}_{1},u_1), \ldots,(\mathpzc{w}_{m},u_m))$, where $\mathpzc{w}_{k}=\{p_k,q_k\}$ 
and $(p_k,q_k,u_k)\in [r]\times [r]\times [t]$ for any $k$, if there exists a tuple
$((\mathpzc{v}_{1},u_1), \ldots,(\mathpzc{v}_{m},u_m))$, where $\mathpzc{v}_k\in \{(p_k,q_k), (q_k,p_k)\}$
for any $k$, to which $\pi$ is adapted. The set of such partitions will be denoted by 
$\mathcal{NC}_{m}^{2}((\mathpzc{w}_1,u_1), \ldots , (\mathpzc{w}_{m},u_m))$.
Its subset consisting of those partitions for which 
$\pi\in \mathcal{NC}_{m,q}^{2}((\mathpzc{v}_1,u_1), \ldots , (\mathpzc{v}_{m},u_m))$ will be denoted 
$\mathcal{NC}_{m,q}^{2}((\mathpzc{w}_1,u_1), \ldots , (\mathpzc{w}_{m},u_m))$.}
\end{Definition}

\begin{Example}
{\rm For simplicity, consider the case when $t=1$ and we can omit $u$'s. If we are given the tuple of sets
$$
(\mathpzc{w}_1, \mathpzc{w}_2, \mathpzc{w}_3, \mathpzc{w}_4)=(\{1,2\},\{1,2\},\{1,2\},\{1,2\}),
$$
then the partition $\pi$ of Fig.1 is adapted to it since there are two tuples,
$$
((2,1),(1,2),(1,2),(2,1))\;\;{\rm and}\;\;((1,2),(2,1),(2,1),(1,2)),
$$
to which $\pi$ is adapted. If $q=1$ or $q=2$, then the associated coloring is given by $f_2$ or  
$f_3$, respectively. Thus, there are
two colored partitions, $(\pi,f_2)$ and $(\pi,f_3)$, associated with the given partition
$\pi$ and the tuple $(\mathpzc{w}_1, \mathpzc{w}_2, \mathpzc{w}_3, \mathpzc{w}_4)$. 
Nevertheless, once we choose the coloring of the
imaginary block, the coloring of $\pi$ is uniquely determined.
In turn, if we are given the tuple of sets
$$
(\mathpzc{w}_1, \mathpzc{w}_2, \mathpzc{w}_3, \mathpzc{w}_4, \mathpzc{w}_{5}, \mathpzc{w}_{6})=
(\{1,2\},\{1,2\},\{1,2\},\{2,2\}, \{2,2\}, \{1,2\})
$$
then $\sigma$ of Fig.1 is adapted to it since it is adapted to
the tuple 
$$
((2,1),(1,2),(1,2),(2,2),(2,2),(2,1)).
$$
Here, there is no other tuple of this type to which $\sigma$ would be adapted
and therefore the only coloring associated with $\pi$ is given by $g$ and
the only color which can be assigned to the imaginary block is $1$.
}
\end{Example}

It can be seen from the above example that if $\pi \in \mathcal{NC}_{m}^{2}((\mathpzc{w}_1,u_1), \ldots , (\mathpzc{w}_{m},u_m))$, 
then there may be more than one colorings of $\pi$ defined by the colorings of the associated
tuples $((\mathpzc{v}_{1},u_1), \ldots , (\mathpzc{v}_{m},u_m))$. Therefore, 
these tuples may produce more than one coloring of $\pi$ defined by the sets of 
matricial indices. However, when we fix $q$ and require 
that $\pi \in \mathcal{NC}_{m,q}^{2}((\mathpzc{v}_1,u_1) \ldots , (\mathpzc{v}_{m},u_m))$, we obtain a unique coloring 
of $\pi$ since we have a unique associated tuple $((\mathpzc{v}_{1},u_1), \ldots , (\mathpzc{v}_{m},u_m))$.

\begin{Proposition}
If $\pi \in \mathcal{NC}_{m,q}^{2}((\mathpzc{w}_1,u_1) \ldots , (\mathpzc{w}_{m},u_m))$, there 
is only one associated tuple $((\mathpzc{v}_{1}, u_1), \ldots , (\mathpzc{v}_{m}, u_{m}))$ for which
$\pi \in \mathcal{NC}_{m,q}^{2}((\mathpzc{v}_1,u_1) \ldots , (\mathpzc{v}_{m},u_m))$. 
\end{Proposition} 
{\it Proof.}
If $\{m-1,m\}$ is a block, then the second index of $\mathpzc{v}_{m-1}=\mathpzc{v}_{m}$ 
must be $q$ and the imaginary block is also the nearest outer block of 
the block containing $m-2$. This allows us to treat the partition $\pi'$ 
obtained from $\pi$ by removing the block $\{m-1,m\}$ in the same way, which gives the inductive step for this (same depth) case.
If $\{m-1,m\}$ is not a block, then the block containing $m$ is the nearest outer block of that containing $m-1$ and
thus the second index of $\mathpzc{v}_{m-1}$ must be equal to the first index of $\mathpzc{v}_{m}$ and the second index of 
$\mathpzc{v}_{m}$ is $q$.
This determines $\mathpzc{v}_{m}$ and $\mathpzc{v}_{m-1}$ and gives the first inductive step for this (growing depth) case. 
Proceeding in this way, we determine $\mathpzc{v}_{1}, \ldots , \mathpzc{v}_{m}$ in a unique way. 
\hfill $\blacksquare$

\begin{Proposition}
For any tuple $((\mathpzc{w}_1,u_1), \ldots , (\mathpzc{w}_m,u_m))$ and $q\in [r]$, where $\mathpzc{w}_{k}=\{p_{k},q_{k}\}$ and
$p_k,q_k\in [r]$, $u_k\in [t]$ for any $k$ and $m\in {\mathbb N}$, it holds that
$$
\Psi_{q}\left(\widehat{\omega}_{p_1,q_1}(u_1)\ldots \widehat{\omega}_{p_m,q_m}(u_m)\right)
=
\sum_{\pi \in \mathcal{NC}_{m,q}^{2}((\mathpzc{w}_1,u_1), \ldots , (\mathpzc{w}_{m},u_m))}
b_{q}(\pi,f,g) 
$$
where $f$ is the coloring of $\pi$ defined by $(\mathpzc{w}_1, \ldots , \mathpzc{w}_m;q)$
and $g$ is the labeling defined by $(u_1, \ldots , u_m)$.
\end{Proposition}
{\it Proof.}
The LHS is a sum of mixed moments of the type computed in Proposition 4.1 with perhaps
some $p_k$'s interchanged with the corresponding $q_k$'s. These moments
are associated with tuples of the form $((\mathpzc{v}_{1}, u_1), \ldots , (\mathpzc{v}_{m},u_m))$.
It follows from the proof of Proposition 4.1 that with each moment of that type we can associate 
the set $\mathcal{NC}_{m,q}^{2}((\mathpzc{v}_1,u_1), \ldots , (\mathpzc{v}_{m},u_m))$ and each such 
moment contributes 
$$
\sum_{\pi\in \mathcal{NC}^{2}_{m,q}((\mathpzc{v}_{1},u_1), \ldots , (\mathpzc{v}_m,u_m))}b_q(\pi,f,g),
$$
where $\mathcal{NC}^{2}_{m,q}((\mathpzc{v}_{1},u_1), \ldots , (\mathpzc{v}_m,u_m))\subseteq
\mathcal{NC}^{2}_{m,q}((\mathpzc{w}_{1},u_1), \ldots , (\mathpzc{w}_m,u_m))$. It is clear from Definition 4.3 that 
$$
\bigcup_{(\mathpzc{v}_1, \ldots , \mathpzc{v}_m)}
\mathcal{NC}^{2}_{m,q}((\mathpzc{v}_{1},u_1), \ldots , (\mathpzc{v}_m,u_m))=
\mathcal{NC}^{2}_{m,q}((\mathpzc{w}_{1},u_1), \ldots , (\mathpzc{w}_m,u_m)),
$$
where the union is taken over the tuples $(\mathpzc{v}_{1}, \ldots, \mathpzc{v}_{m})$ which are 
related to $(\mathpzc{w}_{1}, \ldots , \mathpzc{w}_m)$, i.e. if $\mathpzc{w}_{k}=\{p_k,q_k\}$, then 
$\mathpzc{v}_{k}\in \{(p_k,q_k), (q_k,p_k)\}$. Since the sets on the LHS are disjoint by Proposition 4.2, 
the proof is completed. \hfill $\blacksquare$

\begin{Example}
{\rm Let $p\neq q$ and $u$ (omitted in the notations) be fixed. Then, by Proposition 4.3, 
the moments of the balanced operator $\widehat{\omega}_{p,q}$ take the form
$$
\Psi_{q}(\widehat{\omega}_{p,q}^{m})=\sum_{\pi\in \mathcal{NC}_{m}^{2}}b_{p,q}^{|{\mathcal B}_{o}(\pi)|}b_{q,p}^{|{\mathcal B}_{e}(\pi)|}
$$
where ${\mathcal B}_{o}(\pi)$ and ${\mathcal B}_{e}(\pi)$ denote the blocks of $\pi$ of odd and even depths, respectively,
since all blocks of odd depths contribute $b_{p,q}$ whereas all blocks of even depths contribute $d_{q,p}$.
Therefore, the Cauchy transform of the limit distribution can be represented as the two-periodic continued fraction
with alternating Jacobi coefficients $(b_{p,q}, b_{q,p}, b_{p,q}, \ldots)$. A similar expression is obtained for $\Psi_p$.
In particular, if $(b_{p,q},b_{q,p})=(d_{p},d_q)$, then the moments of $\widehat{\omega}_{p,q}$ in the states $\Psi_q$ and $\Psi_p$ are polynomials in $d_p,d_q$. For instance, 
$$
\Psi_{q}\left(\widehat{\omega}_{p,q}^{6}\right)=d_p^{\,3}+3d_{p}^{\,2}d_q+d_{p}d_{q}^{\,2}
$$
by counting blocks of odd and even depths in the following set of partitions:\\
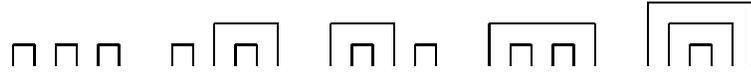
\begin{picture}(200.00,45.00)(-75.00,0.00)
\put(0.00,10.00){\line(0,1){8.00}}
\put(8.00,10.00){\line(0,1){8.00}}
\put(16.00,10.00){\line(0,1){8.00}}
\put(24.00,10.00){\line(0,1){8.00}}
\put(32.00,10.00){\line(0,1){8.00}}
\put(40.00,10.00){\line(0,1){8.00}}

\put(0.00,18.00){\line(1,0){8.00}}
\put(16.00,18.00){\line(1,0){8.00}}
\put(32.00,18.00){\line(1,0){8.00}}
\put(60.00,10.00){\line(0,1){8.00}}
\put(68.00,10.00){\line(0,1){8.00}}
\put(76.00,10.00){\line(0,1){16.00}}
\put(84.00,10.00){\line(0,1){8.00}}
\put(92.00,10.00){\line(0,1){8.00}}
\put(100.00,10.00){\line(0,1){16.00}}

\put(60.00,18.00){\line(1,0){8.00}}
\put(76.00,26.00){\line(1,0){24.00}}
\put(84.00,18.00){\line(1,0){8.00}}

\put(120.00,10.00){\line(0,1){16.00}}
\put(128.00,10.00){\line(0,1){8.00}}
\put(136.00,10.00){\line(0,1){8.00}}
\put(144.00,10.00){\line(0,1){16.00}}
\put(152.00,10.00){\line(0,1){8.00}}
\put(160.00,10.00){\line(0,1){8.00}}

\put(152.00,18.00){\line(1,0){8.00}}
\put(120.00,26.00){\line(1,0){24.00}}
\put(128.00,18.00){\line(1,0){8.00}}

\put(180.00,10.00){\line(0,1){16.00}}
\put(188.00,10.00){\line(0,1){8.00}}
\put(196.00,10.00){\line(0,1){8.00}}
\put(204.00,10.00){\line(0,1){8.00}}
\put(212.00,10.00){\line(0,1){8.00}}
\put(220.00,10.00){\line(0,1){16.00}}

\put(180.00,26.00){\line(1,0){40.00}}
\put(188.00,18.00){\line(1,0){8.00}}
\put(204.00,18.00){\line(1,0){8.00}}

\put(240.00,10.00){\line(0,1){24.00}}
\put(248.00,10.00){\line(0,1){16.00}}
\put(256.00,10.00){\line(0,1){8.00}}
\put(264.00,10.00){\line(0,1){8.00}}
\put(272.00,10.00){\line(0,1){16.00}}
\put(280.00,10.00){\line(0,1){24.00}}

\put(240.00,34.00){\line(1,0){40.00}}
\put(248.00,26.00){\line(1,0){24.00}}

\put(256.00,18.00){\line(1,0){8.00}}

\end{picture}
$\;$\\ 
The moment $\Psi_{p}(\widehat{\omega}_{p,q}^{6})$ is obtained from the above by interchanging
$d_p$ and $d_q$. Thus,
$$
\Psi(\widehat{\omega}_{p,q}^{6})=2d_pd_q^{\,3}+ 6d_p^{\,2}d_q^{2} + 2d_p^{\,3}d_q
$$
since $\Psi=\sum_{j}d_j\Psi_j$ and $\Psi_j(\widehat{\omega}_{p,q}^{6})=0$ for $j\notin \{p,q\}$. 
Moreover, the moment under $\Psi_q$ does not vanish even if $d_p>0$ and $d_q=0$, which corresponds 
to the situation in which $\widehat{\omega}_{p,q}$ is unbalanced. In particular, the first partition in the 
picture given above contributes $d_{p}^{3}$ and is the only one without color $q$.
}
\end{Example}

Finally, we would like to introduce a subset of $\mathcal{NC}_{m,q}^{2}((\mathpzc{w}_1,u_1), \ldots , (\mathpzc{w}_{m},u_m))$ 
consisting of those $\pi$ for which the coloring $f$ satisfies certain additional conditions.
Namely, we would like to distinguish only those colored partitions $(\pi,f)$ whose blocks are not colored by $q$.
Thus, only the imaginary block is colored by $q$, which is not a contradiction since 
the imaginary block is not a block of $\pi$.
This subset will be denoted by $\mathcal{NCI}_{m,q}^{\,2}((\mathpzc{w}_1,u_1), \ldots , (\mathpzc{w}_{m},u_m))$. 

\begin{Example}
{\rm In Fig.1, $(\pi,f_1)\in 
\mathcal{NCI}_{4,1}^{\,2}(\{1,2\}, \{2,2\}, \{2,2\},\{1,2\})$ since $q=1$ and no blocks of $\pi$ are colored by $1$. 
In turn, the remaining colored partitions are not of this type since the colors of their imaginary blocks 
ars assigned to other blocks, too.}
\end{Example}

\begin{Proposition}
For fixed $q\in [r]$, suppose that an operator in any of the arrays $(\widehat{\omega}_{i,j}(u))$ is unbalanced if and only if 
it is of the form $\omega_{p,q}(u)$, where $q \neq p\in [r]$ and $u\in [t]$. Then 
$$
\Psi_{q}\left(\widehat{\omega}_{p_1,q_1}(u_1)\ldots \widehat{\omega}_{p_m,q_m}(u_m)\right)
=
\sum_{\pi \in \mathcal{NCI}_{m,q}^{\,2}((\mathpzc{w}_1,u_1), \ldots , (\mathpzc{w}_{m},u_m))}
b_{q}(\pi,f,g) 
$$
where $f$ is the coloring of $\pi$ defined by $(\mathpzc{w}_1, \ldots , \mathpzc{w}_m;q)$ and
$g$ is the labeling defined by $(u_1, \ldots , u_m)$.
\end{Proposition}
{\it Proof.}
The difference between the considered mixed moment and that of 
Proposition 4.3 is that all operators are balanced except those involving the index $q$.
Namely, if $q\in \{p_k,q_k\}$ for some $k$, then we have $\omega_{p_k,q}(u_k)$ or 
$\omega_{q_k,q}(u_k)$ instead of the symmetrized Gaussian operator $\widehat{\omega}_{p_k,q_k}(u_k)$
for such $k$ and $u_k\in [t]$. Since the first index is these operators is 
different from $q$, we have to eliminate from the set $\mathcal{NC}_{m,q}^{\,2}((\mathpzc{w}_1,u_1) \ldots , (\mathpzc{w}_{m},u_m))$ those partitions 
in which $q$ colors any blocks of $\pi$ since it is always the first index which colors the block. 
This means that we are left only with the contributions from 
$\pi\in \mathcal{NCI}_{m,q}^{\,2}((\mathpzc{w}_1,u_1) \ldots , (\mathpzc{w}_{m},u_m))$, 
which completes the proof.
\hfill $\blacksquare$

\begin{Example}
{\rm Let $p\neq q$ and $u$ (omitted in the notations) be fixed. Then, by Proposition 4.4, 
the moments of the unbalanced operator $\widehat{\omega}_{p,q}=\omega_{p,q}$ take the form
$$
\Psi_{q}(\widehat{\omega}_{p,q}^{m})=b_{p,q}^{m/2}
$$
for even $m$ since the class $\mathcal{NCI}_{m,q}^2(\mathpzc{w}_{1}, \ldots , \mathpzc{w}_{m})$ reduces to 
one interval pair partition with each block colored by $p$ beacause $q$ cannot color any block. 
Of course, if $m$ is odd, we get zero. Moreover,
$$
\Psi_{p}(\widehat{\omega}_{p,q}^{m})=0
$$
since in this case the class $\mathcal{NC}_{m,p}^2(\mathpzc{w}_{1},\ldots, \mathpzc{w}_{m})$ reduces 
to the empty set because $b_{q,p}=0$.}
\end{Example}

The formula of Lemma 4.1 holds irrespective of the number of trivial operators in the arrays $(\omega_{p,q}(u))$ except that
some of the contributions vanish. We can derive a similar formula in the case when the 
only unbalanced operators are of the type studied in Proposition 4.4. 
By $\mathcal{NCI}_{m,q}^{\,2}[r]$ we shall denote the subset of $\mathcal{NC}_{m,q}^{\,2}[r]$ consisting 
of those partitions in which no blocks other than the imaginary block are colored by $q$.

\begin{Lemma}
For fixed $q\in [r]$, suppose that an operator in any of the arrays $(\widehat{\omega}_{i,j}(u))$ is unbalanced if and only if 
it is of the form $\widehat{\omega}_{p,q}(u)$, where $q\neq p\in [r]$ and $u\in [t]$. Then 
\begin{eqnarray*}
\Psi_{q}(\widehat{\omega}^m)&=&\sum_{(\pi,f)\in \mathcal{NCI}_{m,q}^{\,2}[r]}b_{q}(\pi,f)
\end{eqnarray*}
where $b_{q}(\pi,f)=b_{q}(\pi_1,f)\ldots b_{q}(\pi_s,f)$ for $\pi=\{\pi_1, \ldots, \pi_s\}$ and $q\in [r]$, $m\in \mathbb{N}$, 
and 
$$
b_{q}(\pi_k,f)=\sum_{u}b_{i,j}(u)
$$ 
whenever block $\pi_k$ is colored by $i$ and its nearest outer block is colored by $j$.
\end{Lemma}
{\it Proof.}
The proof is similar to that of Lemma 4.1 (Proposition 4.4 is used).
\hfill $\blacksquare$\\

We close this section with a formula for the moment generating functions associated with the moments of
$\omega=\widehat{\omega}$ of Lemma 4.1 (a similar formula can be given for the moments of Lemma 4.2). Let
\begin{eqnarray*}
M(z)&=&{\rm diag}(M_{1}(z), \ldots , M_{r}(z)),
\end{eqnarray*}
where 
$$
M_{q}(z)=\sum_{m=0}^{\infty}\Psi_{q}(\omega^{m})z^{m}
$$
for any $q\in [r]$. The coefficients of the matrix-valued series $M(z)$ obtained in this fashion 
play the role of matricial analogs of Catalan numbers called {\it Catalan matrices}.

Let us introduce the diagonalization mapping
$$
\mathpzc{D}:M_{r}({\mathbb C})\rightarrow M_{r}({\mathbb C})
$$
by the formula
$$
\mathpzc{D}(A)={\rm diag}\left(\sum_{i=1}^{r}a_{i,1}, \ldots , \sum_{i=1}^{r}a_{i,r}\right)
$$
for any $A=(a_{i,j})\in M_{r}({\mathbb C})$.  We denote $B=(b_{p,q})$, where 
$b_{p,q}=\sum_{u\in \mathpzc{U}}b_{p,q}(u)$.

\begin{Proposition}
The matrix-valued generating function $M(z)$ assumes the form
$$
M(z)=\sum_{m=0}^{\infty}\mathcal{C}_{n}z^{2m}
$$
where $(\mathcal{C}_{n})$ is a sequence of diagonal matrices satisfying the recurrence formula
$$
\mathcal{C}_n=\sum_{i+j=n-1}\mathpzc{D}(\mathcal{C}_{i}B\mathcal{C}_{j})
$$
for any natural $n$, where $\mathcal{C}_{0}$ is the $r\times r$ identity matrix. 
The series $M(z)$ converges in the supremum norm for
$|z|<(4r\!\parallel \!B \!\parallel)^{-1}$.
\end{Proposition}
{\it Proof.}
We use the combinatorial expression of Lemma 4.1 to write it in the form of a recurrence
$$
\Psi_{q}(\omega^{2m})=\sum_{p=1}^{r}\sum_{k=0}^{m-1}b_{p,q}\Psi_{p}(\omega^{2k})\Psi_{q}(\omega^{2m-2k-2}),
$$
which can be easily justified. Namely, the class of all colored noncrossing pair partitions 
of $[r]$ is divided into $r$ subclasses in which the block containing $1$, say $\{1,2k+1\}$ is colored by $p\in [r]$. 
This block contributes $b_{p,q}$. Moreover, all colored blocks which lie under this block, playing for them the role 
of the imaginary block, contribute $\Psi_{p}(\omega^{2k})$. In turn, all colored blocks containing numbers $\{2k+2, \ldots , 2m\}$
contribute $\Psi_{q}(\omega^{2m-2k-2})$ since their nearest outer block is still the original imaginary block colored by $q$.
It remains to take the sum over $p$ to obtain the recurrence. The latter leads to the equation
for moment generating functions
$$
M_{q}(z)=1+\sum_{p=1}^{r}b_{p,q}M_{p}(z)M_{q}(z)z^{2},
$$
for any $q\in [r]$, which gives
$$
M(z)=I+{\mathcal D}(M(z)BM(z))z^{2}.
$$
Clearly, $M(z)$ is of the form
$$
M(z)=\sum_{m=0}^{\infty}\mathcal{C}_{m}z^{2m},
$$
where $\mathcal{C}_{m}$ is a constant diagonal matrix for each $m$.
Substituting such series into the above equation
and computing its coefficients, we obtain the desired recurrence formula.
As for the convergence of the series $M(z)$, note that the recurrence formula
is of the same type as for Catalan numbers $C_{m}$ and thus the expression for $\mathcal{C}_{m}$
has $C_{m}$ terms. Using this recurrence, we can estimate the supremum norm of each ${\mathcal C}_{m}$:
$$
\parallel \mathcal{C}_{m}\parallel \leq C_{m}r^{m}\parallel \!B \!\parallel ^{m}
$$
and thus the above series converges for 
$$
|z|<\frac{1}{\sqrt{4r\!\parallel \!B \!\parallel}}
$$
by D'Alembert's test, which completes the proof.
\hfill $\blacksquare$\\

\section{Canonical noncommutative random variables}

Free creation and annihilation operators can be used to construct much more general random variables called 
canonical noncommutative random variables. Let us show that these random variables can be decomposed
into sums of matricial type.
  
Let $\{\ell(u):u\in \mathpzc{U}\}$ be a family of *-free standard creation operators, by which we mean that
$\ell(u)^*\ell(v)=\delta_{u,v}$, on the free Fock space
$$
{\mathcal F}={\mathcal F}(\bigoplus_{u\in \mathpzc{U}}{\mathbb C}e(u)).
$$
Using these operators and their adjoints one can 
construct free random variables $\gamma(u)$ whose all moments coincide with the moments of given probability measures $\mu(u)$.

\begin{Definition}
{\rm Let $\mu(u)$ be a probability measure on the real line whose free cumulants are $(r_{k}(u))_{k\geq 1}$, respectively.
The formal sums 
$$
\gamma(u)=\ell(u)^*+\sum_{k=0}^{\infty}r_{k+1}(u)(\ell(u))^{k}
$$
are called {\it canonical noncommutative random variables}. If $\sum_{k=0}^{\infty}|r_{k+1}(u)|<\infty$, 
then $\gamma(u)$ is a bounded operator on ${\mathcal F}$. }
\end{Definition}
\begin{Remark}
{\rm 
If $(r_{k}(u))$ is an arbitrary sequence of real numbers and $|\mathpzc{U}|=t$, 
we can treat $\gamma(u)$ as an element of the unital *-algebra $\widetilde{\mathcal{E}}_{t}$ of formal sums
$$
\sum_{\stackrel 
{p\geq 0}
{\scriptscriptstyle 0\leq q \leq Q}}
\sum_{\stackrel
{u_1, \ldots , u_p\in [t]}
{\scriptscriptstyle v_1, \ldots , v_p\in [t]}}
c_{u_1,\ldots , u_p, v_1, \ldots , v_q}\ell(u_1)\ldots \ell(u_p)\ell(v_1)^{*}\ldots \ell(v_q)^{*}
$$
equipped with the linear functional $\phi_t$ sending the above element to $c_{\emptyset, \emptyset}$, the coefficient
corresponding to $p=q=0$. }
\end{Remark}

Therefore, by Proposition 2.2, we can also obtain realizations of bounded free canonical noncommutative random variables 
$\gamma(u)$ on ${\mathcal M}$ by setting $\ell(u)=\wp(u)$.
If $\gamma(u)$ is not bounded, then there exists a 
sequence $(\gamma^{(n)}(u))$ of bounded free 
canonical noncommutative random variables such that all moments of orders $\leq n$ of $\gamma^{(n)}(u)$ in the state $\Psi$
agree with the corresponding moments of $\gamma(u)$ in the state $\phi_t$. Consequently, this property will hold true if we
compute mixed moments of the whole family $\{\gamma(u): u\in \mathpzc{U}\}$. 

\begin{Proposition}
If $\sum_{k=0}^{\infty}|r_{k+1}(u)|<\infty$, then the canonical noncommuative random variable $\gamma(u)$ has the decomposition
$$
\gamma(u)=\sum_{p,q=1}^{r}\gamma_{p,q}(u),
$$
for any $r\in \mathbb{N}$, where
\begin{eqnarray*}
\gamma_{p,q}(u)&=&\wp_{p,q}(u)^* +\delta_{p,q}r_{1}(u)P_{q}\\
&+&
\sum_{k=1}^{\infty}r_{k+1}(u)\sum_{q_1, \ldots , q_{k-1}}
\wp_{p,q_1}(u)\wp_{q_1,q_2}(u)\ldots \wp_{q_{k-1},q}(u)
\end{eqnarray*}
for any $p,q,u$, where each $\wp_{p,q}(u)$ has covariance $d_{p}$ and $d_1+\ldots +d_r=1$. 
\end{Proposition}
{\it Proof.}
It follows from Proposition 2.2 that we can realize canonical free creation operators $\ell(u)$ 
as $\wp(u)$, where each $\wp_{p,q}(u)$ has covariance $d_p$. Thus, the above decomposition of canonical random variables 
follows from the following elementary computation:
\begin{eqnarray*}
\gamma(u)&=&\sum_{p,q}P_p(\wp(u)^* +\sum_{k=0}^{\infty}r_{k+1}(u)(\wp(u))^{k}) P_{q}\\
&=&\sum_{p,q}P_q\wp(u)^*P_p +
\delta_{p,q}r_{1}(u)P_{q}\\
&+&\sum_{p,q}\sum_{k=1}^{\infty}r_{k+1}(u)\sum_{q_1, \ldots , q_{k-1}}
\wp_{p,q_1}(u)\wp_{q_1,q_2}(u)\ldots \wp_{q_{k-1},q}(u_k)\\
&=&\sum_{p,q}\gamma_{p,q}(u)
\end{eqnarray*}
which completes the proof.\hfill $\blacksquare$\\

\begin{Remark}
{\rm Let us observe that the off-diagonal operators $\gamma_{p,q}(u)$ are generalizations of matricially free 
Gaussian operators and can also be written in the form
$$
\gamma_{p,q}(u)=P_{q}\wp(u)^*P_{p}+ \sum_{k=0}^{\infty}r_{k+1}(u)P_{p}(\wp(u))^k P_{q}.
$$ 
It will turn out that they describe the asymptotics of unbalanced blocks of Hermitian 
random matrices which are asymptotically free. Their symmetrized couterparts
$$
\widehat{\gamma}_{p,q}(u)=\gamma_{p,q}(u)+\gamma_{q,p}(u)=P_{p}\gamma(u)P_{q}+P_{q}\gamma(u)P_{p}
$$
for $p\neq q$ and 
$$
\widehat{\gamma}_{q,q}(u)=\gamma_{q,q}(u)=P_{q}\gamma(u)P_{q}
$$
for any $q$, 
are generalizations of symmetrized Gaussian operators and they describe the asymptotics of symmetric blocks 
of these matrices. Note that $\gamma_{p,q}(u)\neq P_{p}\gamma(u)P_q$ for $p\neq q$. }
\end{Remark}

The combinatorics of mixed moments of canonical noncommutative random variables is based on 
all noncrossing partitions. Let us observe that with the mixed moment $\Phi(\gamma(u_1)\ldots \gamma(u_m))$
we can associate noncrossing partitions $\pi=\{\pi_1, \ldots , \pi_s\}$ of the set $[m]$
such that $u_i=u$ whenever $i\in \pi_k$. Then, we can associate a power of a free creation operator $\ell(u)^{k}$ with the 
last leg of any block and $\ell(u)^*$ with its remaining legs unless the block is a singleton, to which 
we associate the unit. 
We show below that a similar procedure can be applied to the moments 
$\Psi_{q}\left(\gamma_{p_1,q_1}(u_1)\ldots \gamma_{p_m,q_m}(u_m)\right)$. 
Then, it can be seen that we arrive at a generalization of Definition 4.1.

\begin{Definition}
{\rm We will say that $\pi \in \mathcal{NC}_{m}$ is {\it adapted}
to the tuple $((\mathpzc{v}_1,u_1), \ldots, (\mathpzc{v}_{m},u_m))$, where 
$\mathpzc{v}_{i}=(p_i,q_i)\in [r]\times [r]$ and $u_i\in [t]$ for any $i$, if
for any block $(i(1)<\ldots <i(k))$ of $\pi$ it holds that
\begin{enumerate}
\item[(a)]
$u_{i(1)}=\ldots =u_{i(k)}$,
\item[(b)]
$(p_{i(k-1)},q_{i(1)})=(p_{i(k)},q_{i(k)})$,
\item[(c)]
$(q_{i(1)}, q_{i(2)}, \ldots, q_{i(k-1)})=(p_{o(i(1))}, p_{i(1)}, \ldots,p_{i(k-2)})$ 
whenever the block has an outer block,
\end{enumerate}
where $o(i(1))$ is the smallest number in the nearest outer block of the given block which is 
greater than $i(1)$. Denote by $\mathcal{NC}_{m}((\mathpzc{v}_1,u_1), \ldots , (\mathpzc{v}_{m},u_m))$
the set of such partitions and by $\mathcal{NC}_{m,q}((\mathpzc{v}_1,u_1), \ldots , (\mathpzc{v}_{m},u_m))$
its subset for which $q_{i(k)}=q$ whenever the block $(i(1)<\ldots <i(k))$ is a covering block.
}
\end{Definition}

We also need to generalize Definition 4.2. For that purpose, we need to 
color arbitrary noncrossing partitions. However, one color for each block 
is not enough. If $\pi\in \mathcal{NC}_{m}$, then the associated
colored noncrossing partition will be the triple $(\pi,f,g)$, where 
\begin{enumerate}
\item
$f$ is an $[r]$-valued function which assigns a number to each 
singleton and to each subblock consisting of two consecutive numbers of any block $(i(1)<\ldots <i(k))$,
\item
$g$ is an $[t]$-valued function which assigns an element of $[t]$ to each block.
\end{enumerate}
The set of colorings of $\pi$ by functions $f$ and $g$ will be again denoted by $F_{r}(\pi)$
and $F_{t}(\pi)$, respectively. We will say that block $(i_1<\ldots <i_k)$ is {\it colored} by  $(p_1, \ldots , p_{k-1})$ if
$f$ assigns $p_j$ to the subblock joining $i_j$ and $i_{j+1}$ for any $j\in [k-1]$.
We will say that this block is {\it labeled} by $u$ if $g$ assigns $u$ to this block. 

\begin{Definition}
{\rm 
Let $D=(d_1, \ldots, d_r)$ be the dimension matrix, 
and let $(r_{k}(u))$ be sequences of free cumulants. Let $\pi=\{\pi_1, \ldots , \pi_n\}\in \mathcal{NC}_{m}$ and 
$f\in F_{r}(\pi)$, $g\in F_{t}(\pi)$. If a block $\pi_j=(i_1<\ldots <i_k)$ is colored by $(p_1, \ldots, p_{k-1})$
and labeled by $u$, where $k>1$, we assign to it the number 
$$
b(\pi_j,f,g)=d_{p_1}\ldots d_{p_{k-1}}r_{k}(u),
$$
and we assign $r_{1}(u)$ to any singleton. 
Extending this definition by multiplicativity over blocks, 
$$
b(\pi, f,g)=b(\pi_{1},f, g)\ldots b(\pi_n,f, g),
$$
we obtain real-valued functions on the set of colored noncrossing partitions. 
}
\end{Definition}

\begin{Remark}
{\rm For simplicity, the above $b$ does not depend on the colors $q$ of the imaginary blocks. 
However, we can use more general functions, similar to those of Definition 4.2.
Consider arrays $(c_{p,q}(u))$ of real numbers and assign to $\pi_j$ the number
$$
b_{q}(\pi_j,f,g)=d_{p_1}\ldots d_{p_{k-1}}c_{p_{0},p_1}(u)\ldots c_{p_{k-1},p_{0}}(u)r_{k}(u),
$$
where $p_0$ is the color of its nearest outer subblock and $q$ is the color of the imaginary block (thus,
we assign $r_{1}(u)c_{p_0,p_0}$ to a singleton). 
Let us observe that in the case when $\pi_j$ is a 2-block and $r_2(u)=1$, we have
$b_{q}(\pi_j,f,g)=d_{p}c_{q,p}(u)c_{p,q}(u)=b_{p,q}(u)$ if $(c_{p,q}(u))$ is symmetric, which 
reduces to Definition 4.2.}
\end{Remark}

\begin{Example}
{\rm Consider the noncrossing partition $\pi$ of the set $[8]$ consisting of four blocks:
$\pi_1=\{1,8\},\pi_2=\{2,3,4\},\pi_3=\{5\},\pi_4=\{6,7\}$, colored as shown below and labelled by
the same $u$ (thus $g$ is omitted in the notations).

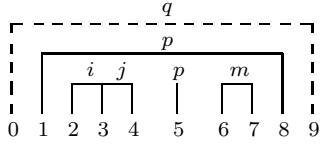
\unitlength=1mm
\special{em:linewidth 0.4pt}
\linethickness{0.4pt}
\begin{picture}(160.00,30.00)(70.00,1.00)

\put(129.00,10.00){\line(0,1){8.00}}
\put(133.00,10.00){\line(0,1){4.00}}
\put(137.00,10.00){\line(0,1){4.00}}
\put(141.00,10.00){\line(0,1){4.00}}
\put(147.00,10.00){\line(0,1){4.00}}
\put(153.00,10.00){\line(0,1){4.00}}
\put(157.00,10.00){\line(0,1){4.00}}
\put(161.00,10.00){\line(0,1){8.00}}

\put(125.00,10.00){\line(0,1){1.25}}
\put(125.00,12.50){\line(0,1){1.25}}
\put(125.00,13.00){\line(0,1){1.25}}
\put(125.00,15.50){\line(0,1){1.25}}
\put(125.00,18.00){\line(0,1){1.25}}
\put(125.00,20.75){\line(0,1){1.25}}

\put(165.00,10.00){\line(0,1){1.25}}
\put(165.00,12.50){\line(0,1){1.25}}
\put(165.00,13.00){\line(0,1){1.25}}
\put(165.00,15.50){\line(0,1){1.25}}
\put(165.00,18.00){\line(0,1){1.25}}
\put(165.00,20.75){\line(0,1){1.25}}

\put(124.50,7.00){$\scriptstyle{0}$}
\put(128.50,7.00){$\scriptstyle{1}$}
\put(132.50,7.00){$\scriptstyle{2}$}
\put(136.50,7.00){$\scriptstyle{3}$}
\put(140.50,7.00){$\scriptstyle{4}$}
\put(146.50,7.00){$\scriptstyle{5}$}
\put(152.50,7.00){$\scriptstyle{6}$}
\put(156.50,7.00){$\scriptstyle{7}$}
\put(160.50,7.00){$\scriptstyle{8}$}
\put(164.50,7.00){$\scriptstyle{9}$}

\put(135.00,15.00){$\scriptstyle{i}$}
\put(139.00,15.00){$\scriptstyle{j}$}
\put(146.50,15.00){$\scriptstyle{p}$}
\put(154.00,15.00){$\scriptstyle{m}$}
\put(145.00,19.00){$\scriptstyle{p}$}
\put(145.00,23.50){$\scriptstyle{q}$}

\put(125.00,22.00){\line(1,0){1.50}}
\put(128.00,22.00){\line(1,0){1.25}}
\put(130.50,22.00){\line(1,0){1.25}}
\put(133.00,22.00){\line(1,0){1.25}}
\put(135.50,22.00){\line(1,0){1.25}}
\put(138.00,22.00){\line(1,0){1.25}}
\put(140.50,22.00){\line(1,0){1.25}}
\put(143.00,22.00){\line(1,0){1.25}}
\put(145.50,22.00){\line(1,0){1.25}}
\put(148.00,22.00){\line(1,0){1.25}}
\put(150.50,22.00){\line(1,0){1.25}}
\put(153.00,22.00){\line(1,0){1.25}}
\put(155.50,22.00){\line(1,0){1.25}}
\put(158.00,22.00){\line(1,0){1.25}}
\put(160.75,22.00){\line(1,0){1.35}}
\put(163.50,22.00){\line(1,0){1.50}}

\put(129.00,18.00){\line(1,0){32.00}}
\put(133.00,14.00){\line(1,0){8.00}}
\put(153.00,14.00){\line(1,0){4.00}}
\end{picture}
\noindent
By Definition 5.3, 
$$
b(\pi_1,f)=d_{p}r_{2},\;\; b(\pi_2,f)=d_{i}d_{j}r_{3}, \;\;b(\pi_3,f)=r_1,\;\;b(\pi_4,f)=d_{m}r_2
$$
where $r_j=r_{j}(u)$ for any $j$. For arbitrary scaling constants, we need to multiply
these contributions by $c_{p,q}c_{q,p}$, $c_{p,i}c_{i,j}c_{j,p}$, $c_{p,p}$ and 
$c_{p,m}c_{m,p}$, respectively. 
}
\end{Example}

In order to find a combinatorial formula for the mixed moments of operators $\gamma_{p,q}(u)$, we need
to color $\pi\in \mathcal{NC}_{m,q}((\mathpzc{v}_1,u_1), \ldots , (\mathpzc{v}_{m},u_m))$
in a suitable way. Namely, we will color each block $(i_1<\ldots <i_k)$ of $\pi$ by $(p_{i_1}, \ldots , p_{i_{k-1}})$
and we will label it by $u$ whenever $u_{i_1}=\ldots=u_{i_{k-1}}=u$. The imaginary block is colored by $q$. 
We will say that this coloring of is defined by $(\mathpzc{v}_{1}, \ldots , \mathpzc{v}_m;q)$.
Similarly, we will say that this labeling of $\pi$ is defined by $(u_1, \ldots , u_m)$.

\begin{Proposition}
For any tuple $((\mathpzc{v}_{1},u_1), \ldots , (\mathpzc{v}_m,u_m))$ and $q\in [r]$, $m\in {\mathbb N}$,
where $\mathpzc{v}_{k}=(p_k,q_k)\in [r]\times [r]$ and $u_k\in [t]$ for each $k$, it holds that
$$
\Psi_{q}\left(\gamma_{p_1,q_1}(u_1)\ldots \gamma_{p_m,q_m}(u_m)\right)
=
\sum_{\pi\in \mathcal{NC}_{m,q}((\mathpzc{v}_{1},u_1), \ldots , (\mathpzc{v}_m,u_m))}b_q(\pi,f,g)
$$
where $f$ is the coloring of $\pi$ defined by $(\mathpzc{v}_{1}, \ldots , \mathpzc{v}_m;q)$ and $g$ is the 
labeling of $\pi$ defined by $(u_1, \ldots , u_m)$.
\end{Proposition}
{\it Proof.}
In the computation of the moment $\Psi_{q}\left(\gamma_{p_1,q_1}(u_1)\ldots \gamma_{p_m,q_m}(u_m)\right)$ we obtain 
the sum of products of operators of three types:
\begin{enumerate}
\item
projection $P_{q_j}$ corresponding to $k=0$ in the definition of $\gamma_{p_j,q_j}(u_j)$ (these appear only if $p_j=q_j$), to which we 
associate a singleton,
\item
generalized creation operator $P_{p_j}\wp^{k}(u_j)P_{q_j}$ for $k>0$, to which we
associate the last leg of some block,
\item
generalized annihilation operator $P_{q_j}\wp(u)^*P_{p_j}$ to which we associate a leg of some block (except the last one),
\end{enumerate}
where each $j\in [m]$ appears exactly once.
The generalized annihilation operators are matricially free annihilation operators, 
$$
P_{q_j}\wp(u_j)^*P_{p_j}=\wp_{p_j,q_j}(u_j)^*,
$$
and the generalized creation operators are matricial products
of matricially free creation operators by Proposition 5.1, namely
$$
P_{p_j}\wp^{k-1}(u_j)P_{q_j}=\sum_{s_1, \ldots , s_{k-2}}\wp_{p_j,s_1}(u_j)\wp_{s_1,s_2}(u_j)\ldots \wp_{s_{k-2},q_j}(u_j),
$$
where $k>1$.
By the definition of matricially free creation and annihilation operators, it can be seen that
a product of operators of the above three types is not zero if and only if
for each generalized creation operator $P_{p_j}\wp^{k-1}P_{q_{j}}$, there exists a
sequence of $k-1$ matricially free annihilation operators, $(\wp_{s_{k-2},q_j}(u_j)^*, \ldots , \wp_{p_j,s_1}(u_j)^*)$
for some $s_{1}, \ldots , s_{k-2}$, associated with an increasing sequence of numbers $(i_1,\ldots ,i_{k-1})$, all smaller than $j$, 
thus it holds that
$$
(p_{i_1},q_{i_1})=(s_{k-2},q_j),  \ldots, (p_{i_{k-1}}, q_{i_{k-1}})=(p_j,s_1)
$$
with $u_{i_1}=\ldots = u_{i_{k-1}}=u_j$, and moreover, $q_{j}=p_{o(i_1)}$, where 
$o(i_1)$ is the smallest index greater than $i_k$ corresponding to an operator 
belonging to another product of this type within which the considered product is nested 
in a noncrossing way (if such a product does not exist, $q_j=q$).
In other words, we must have $(q_{i_1},q_{i_2}, \ldots , q_{i_{k-1}})=(p_{o(i_1)}, p_{i_1}, \ldots , p_{i_{k-2}})$
and also $(p_{i_k},q_{i_k})=(p_{i_{k-1}},q_{i_1})$ since $i_k=j$. These are precisely 
the conditions of Definition 5.1, which proves that the considered mixed moment
is a sum o mixed moments of matricially free creation and annihilation operators and projections 
which correspond to all noncrossing partitions in the sense described above.  
To compute the mixed moment corresponding to $\pi$, it suffices to compute the contribution from products
corresponding to the blocks of $\pi$. If such a block is of the form $(i<i+1<\ldots <k)$ (there is at least one block like this), 
then it contributes
$$
\Psi_p(\wp_{p_i,p}(u)^*\ldots \wp_{p_{k-1},p_{k-2}}(u)^{*}\wp_{p_{k-1},p_{k-2}}(u)\ldots \wp_{p_{i},p}(u))
=d_{p_i}\ldots d_{p_{k-1}}r_{k}(u)
$$
where $p$ is the coloring of its nearest outer block.
By induction with respect to the depth of blocks, a similar contribution
is obtained for the remaining blocks. The contribution from a partition
$\pi\in \mathcal{NC}_{m,q}((\mathpzc{v}_{1},u_1), \ldots , (\mathpzc{v}_m,u_m))$
is the product over blocks of such expressions in agreement with 
Definition 5.3. This completes the proof.
\hfill $\blacksquare$

\begin{Definition}
{\rm We say that $\pi \in \mathcal{NC}_{m}$ is {\it adapted}
to the tuple $((\mathpzc{w}_{1},u_1), \ldots,(\mathpzc{w}_{m},u_m))$, where $\mathpzc{w}_{k}=\{p_k,q_k\}$ 
and $(p_k,q_k,u_k)\in [r]\times [r]\times [t]$ for any $k$, if there exists a tuple
$((\mathpzc{v}_{1},u_1), \ldots,(\mathpzc{v}_{m},u_m))$, where $\mathpzc{v}_k\in \{(p_k,q_k), (q_k,p_k)\}$
for any $k$, to which $\pi$ is adapted. The set of such partitions will be denoted by 
$\mathcal{NC}_{m}((\mathpzc{w}_1,u_1), \ldots , (\mathpzc{w}_{m},u_m))$.
Its subset consisting of those partitions for which 
$\pi\in \mathcal{NC}_{m,q}((\mathpzc{v}_1,u_1), \ldots , (\mathpzc{v}_{m},u_m))$ will be denoted 
$\mathcal{NC}_{m,q}((\mathpzc{w}_1,u_1), \ldots , (\mathpzc{w}_{m},u_m))$. Next,
its subset consisting of those partitions in which $q$ does not color any blocks but the imaginary block
will be denoted $\mathcal{NCI}_{m,q}((\mathpzc{w}_{1},u_1), \ldots , (\mathpzc{w}_m,u_m))$.
Finally, we will say that the coloring of $\pi$ is defined by $(\mathpzc{w}_{1}, \ldots , \mathpzc{w}_m)$ 
if it is inherited from $(\mathpzc{v}_{1}, \ldots,\mathpzc{v}_{m})$. 
}
\end{Definition}

\begin{Proposition}
For any tuple $((\mathpzc{w}_{1},u_1), \ldots , (\mathpzc{w}_m,u_m))$ and $q\in [r]$, $m\in {\mathbb N}$,
where $\mathpzc{w}_{k}=\{p_k,q_k\}$ and $u_k\in [t]$ for each $k$, it holds that
$$
\Psi_q\left(\widehat{\gamma}_{p_1,q_1}(u_1)\ldots \widehat{\gamma}_{p_m,q_m}(u_m)\right)
=
\sum_{\pi\in \mathcal{NC}_{m,q}((\mathpzc{w}_{1},u_1), \ldots , (\mathpzc{w}_m,u_m))}b_q(\pi,f,g)
$$
where $f$ is the coloring of $\pi$ defined by $(\mathpzc{w}_{1}, \ldots , \mathpzc{w}_m;q)$
and $g$ is the labeling of $\pi$ defined by $(u_1, \ldots , u_m)$.
The summation reduces to $\mathcal{NCI}_{m,q}((\mathpzc{w}_{1},u_1), \ldots , (\mathpzc{w}_m,u_m))$
if $d_q=0$.
\end{Proposition}
{\it Proof.}
The first statement is a consequence of Proposition 5.2. The second statement follows from the fact that
if a segment of a block of $\pi$ is colored by $q$, then $d_q$ appears in the formula for $b_{q}(\pi,f)$
and since $d_q=0$, the contribution of such $\pi$ vanishes.
\hfill $\blacksquare$

\begin{Example}
{\rm Let us compute the lowest order mixed moments of operators of type $\gamma_{p,q}=\gamma_{p,q}(u)$ for fixed $u$, 
supposing that all scaling constants are equal to one. We obtain
\begin{eqnarray*}
\Psi_q(\gamma_{p,q})&=&\delta_{p,q}r_{1}\Psi_{q}(P_{q})=\delta_{p,q}r_{1}\\
\Psi_q(\gamma_{p,q}\gamma_{p,q})&=&r_{2}\Psi_{q}(\wp_{p,q}^*\wp_{p,q})+\delta_{p,q}r_{1}^{2}\Psi_{q}(P	_{q})
= r_{2}d_{p}+\delta_{p,q}r_{1}^{2}\\
\Psi_q(\gamma_{s,q}\gamma_{s,p}\gamma_{p,q})&=&r_{3}\Psi_{q}(\wp_{s,q}^*\wp_{p,s}^*\wp_{p,s}\wp_{s,q})\\
&+&
\delta_{s,p}r_{1}r_{2}\Psi_{q}(\wp_{p,q}^{*}P_{p}\wp_{p,q})+
\delta_{p,s}\delta_{p,q}r_{1}r_{2}\Psi_{q}(P_{q}\wp_{q,q}^*\wp_{q,q})\\
&+&
\delta_{p,q}r_{1}r_{2}\Psi_{q}(\wp_{s,q}^*\wp_{s,q}P_{q})+
\delta_{p,q}\delta_{p,s}r_{1}^3\Psi_{q}(P_{q})\\
&=&
r_{3}d_pd_s+\delta_{s,p}r_{1}r_{2}d_p
+\delta_{p,s}\delta_{p,q}r_{1}r_{2}d_q
+\delta_{p,q}r_{1}r_{2}d_s
+\delta_{p,q}\delta_{p,s}r_{1}^3
\end{eqnarray*} 
where we set $r_{k}=r_{k}(u)$ for any $k$.
Similar moments of the symmetrized operators $\widehat{\gamma}_{p,q}$ can be obtained by linearity.
}
\end{Example}

All mixed moments of Proposition 5.3 are polynomials in asymptotic dimensions. 
This case will be interesting in the random matrix context especially if not all asymptotic 
dimensions are equal. Although finding explicit forms of such 
polynomials may be highly non-trivial, a combinatorial formula can be given. 

\begin{Corollary}
Under the assumptions of Proposition 5.3, 
$$
\Psi_q\left(\widehat{\gamma}_{p_1,q_1}(u_1)\ldots \widehat{\gamma}_{p_m,q_m}(u_m)\right)
$$
$$
=
\sum_{\pi\in \mathcal{NC}_{m}((\mathpzc{w}_{1},u_1), \ldots , (\mathpzc{w}_{m},u_m))}
\prod_{\scriptscriptstyle{{\rm blocks}}\;\pi_k}r_{|\pi_k|}(u^{(k)})
\prod_{j\,\in \mathpzc{w}_{1}\cup \ldots \cup \mathpzc{w}_{m}}
d_{j}^{\,|\mathcal{S}_{j}(\pi)|}
$$
where $\mathcal{S}_{j}(\pi)$ is the set of subblocks of blocks of $\pi$ of the form $\{i,i+1\}$ which are colored by $j$, 
and $u^{(k)}=u_i$ for all $i\in \pi_k$. 
\end{Corollary}
{\it Proof.}
Using the formula of Proposition 5.3, we observe that each block $\pi_k$ contributes a cumulant of order 
$|\pi_k|$ and a product of $|\pi_k|-1$ asymptotic dimensions (each subblock $\{i,i+1\}$ colored by $j$ contributes $d_j$).
This finishes the proof.
\hfill $\blacksquare$\\

\section{Hermitian random matrices}

The context for the study of random matrices originated by Voiculescu [33] is 
the following. Let $\mu$ be a probability measure on some measurable space without atoms and let 
$L=\bigcap_{1\leq p <\infty}L^{p}(\mu)$ be endowed with the state
expectation ${\mathbb E}$ given by integration with respect to $\mu$. 
The *-algebra of $n\times n$ random matrices is 
$M_{n}(L)=L\otimes M_{n}({\mathbb C})$ with the state 
$$
\tau(n)={\mathbb E}\otimes {\rm tr}(n)
$$
where ${\rm tr}(n)$ is the normalized trace over the set $\{e(j):j\in [n]\}$ of basis vectors of ${\mathbb C}^{n}$.

In our study, we replace the complete trace $\tau(n)$ by {\it partial traces} 
$$
\tau_{q}(n)={\mathbb E} \otimes {\rm tr}_{q}(n),
$$
where ${\rm tr}_{q}(n)$ is the normalized trace over the set of basis vectors of ${\mathbb C}^{n}$ indexed by $j\in N_{q}$, that is,
the trace divided by $n_{q}=|N_{q}|$.

\begin{Definition}
{\rm The {\it symmetric blocks} of an $n\times n$ random matrix $Y(u,n)$ are matrices of the form 
$$
T_{p,q}(u,n)=\left\{
\begin{array}{ll}
D_{q}Y(u,n)D_{q}& {\rm if}\;\; p=q\\
D_{p}Y(u,n)D_{q}+D_{q}Y(u,n)D_{p}& {\rm if}\;\;p\neq q
\end{array}
\right.
$$
for any $p,q\in [r]$, $u\in \mathpzc{U}$ and $n\in \mathbb{N}$, where
$\{D_1, \ldots , D_r\}$ is the family of $n\times n$ diagonal matrices ($n$ is suppressed in the notation) 
which gives a natural decomposition of the identity correponding to the partition of $[n]$ defined above.
Thus, $(D_{q})_{j,j}=1$ if and only if $j\in N_{q}$ and the remaining entries of $D_{q}$ are zero.}
\end{Definition}

\begin{Definition}
{\rm The sequence of symmetric blocks $(T_{p,q}(u,n))_{n\in \mathbb{N}}$, where $(p,q)\in \mathpzc{J}$, $u\in \mathpzc{U}$ are fixed, 
will be called {\it balanced} if $d_{p}>0$ and $d_q>0$, {\it unbalanced} if $(d_p=0\;{\rm and}\;d_q>0)$ or $(d_p>0\;{\rm and}\;d_q=0)$,
and {\it evanescent} if $d_p=0$ and $d_q=0$ (cf. Definition 3.2).}
\end{Definition}

\begin{Example}
{\rm Consider the sequence of $n\times n$ random matrices consisting of three symmetric blocks: 
two diagonal Hermitian blocks,
$$
T_{1,1}(n)=\left(
\begin{array}{cc}
S_{1,1}(n)& 0\\
0& 0
\end{array}
\right), \;\;\;
T_{2,2}(n)=\left(
\begin{array}{cc}
0& 0\\
0& S_{2,2}(n)
\end{array}
\right),
$$
and the off-diagonal one,
$$
T_{1,2}(n)=\left(
\begin{array}{cc}
0& S_{1,2}(n)\\
S_{1,2}^{*}(n)& 0
\end{array}
\right),
$$
where $u$ is omitted for simplicity. If we suppose that 
$$
\lim_{n\rightarrow \infty}\frac{n_{1}}{n}=d_1=0\;\; {\rm and} \;\;
\lim_{n\rightarrow \infty}\frac{n_{2}}{n}=d_2>0,
$$
then $(T_{1,1}(n))_{n\in \mathbb{N}}$ is evanescent, 
$(T_{2,2}(n))_{n\in \mathbb{N}}$ is balanced and 
$(T_{1,2}(n))_{n\in \mathbb{N}}$ is unbalanced.
It will follow from Theorem 6.1 that the moments of such symmetric blocks under $\tau_q(n)$
tend to the moments of trivial, balanced and unbalanced symmetrized Gaussian operators,
$\widehat{\omega}_{1,1}=0$,
$\widehat{\omega}_{2,2}=\omega_{2,2}$ and
$\widehat{\omega}_{1,2}=\omega_{2,1}$,
respectively, under $\Psi_q$, where $q\in \{1,2\}$.
}
\end{Example}

It is easy to predict that if $D$ is singular, then all mixed moments involving 
evanescent blocks will tend to zero. However, if we compute mixed moments involving 
unbalanced blocks consisting of pairs of mutually adjoint rectangular blocks whose 
one dimension grows too slowly, it is less obvious that one can still obtain nonvanishing limit moments under 
a suitable partial trace.

\begin{Theorem}
Let $\{Y(u,n):u\in \mathpzc{U}, n\in {\mathbb N}\}$ be a family of independent Hermitian random 
matrices whose asymptotic joint distribution under $\tau(n)$ agrees with that of the 
family $\{\gamma(u):u\in \mathpzc{U}\}$ under $\Psi=\sum_{q}d_q\Psi_q$ 
and which is asymptotically free from $\{D_1, \ldots , D_r\}$.
If $\parallel \!Y(u,n)\!\parallel \leq C$ for any $n,u$ and some $C$, then
$$
\lim_{n\rightarrow \infty}\tau_q(n)(T_{p_1,q_1}(u_1,n)\ldots T_{p_m,q_m}(u_m,n))
=
\Psi_q(\widehat{\gamma}_{p_1,q_1}(u_1)\ldots \widehat{\gamma}_{p_m,q_m}(u_m))
$$
as $n\rightarrow \infty$ for any $u_1, \ldots, u_m\in \mathpzc{U}$ and $p_1,q_1, \ldots , p_m,q_m,q\in [r]$.
\end{Theorem}
{\it Proof.}
The mixed moments of the blocks 
$S_{i,j}(u,n)=D_{i}Y(u,n)D_{j}$ under $\tau_q(n)$ can be written as 
mixed moments under $\tau(n)$ of the form
$$
\tau_{q}(n)(S_{p_1,q_1}(u_1)\ldots S_{p_m,q_m}(u_m))
$$
$$
=
\frac{n}{n_{q}}\tau(n)(D_{q}D_{p_1}Y(u_1,n)D_{q_1}\ldots D_{p_m}Y(u_m,n)D_{q_m}D_{q})
$$
for all values of the indices. 

Suppose first that $d_q>0$.
By assumption, the family $\{Y(u,n):u\in \mathpzc{U}\}$ is asymptotically free under 
$\tau(n)$ and free from $\{D_1, \ldots , D_r\}$. Therefore, by Lemma 2.2, the moment on the RHS converges to the corresponding 
mixed moment
$$
\frac{1}{d_{q}}\Psi(P_{q}P_{p_1}\gamma(u_1)P_{q_1}\ldots P_{p_m}\gamma(u_m)P_{q_m}P_{q})
$$
$$
=
\Psi_q(P_{p_1}\gamma(u_1)P_{q_1}\ldots P_{p_m}\gamma(u_m)P_{q_m})
$$
as $n\rightarrow \infty$ provided the asymptotic joint distribution of $\{D_1, \ldots , D_r\}$ under $\tau(n)$ 
agrees with that of $\{P_1, \ldots , P_r\}$ under $\Psi$. That is easily checked 
since $P_p\perp P_q$ for $p\neq q$ and 
$$
\tau(n)(D_q)=\frac{n_q}{n}\rightarrow d_{q}=\Psi(P_q).
$$
The corresponding formula for symmetric blocks easily follows. It is worth to remark that, in general, 
the moments of blocks $S_{p,q}(u,n)$ do not tend to the corresponding moments of operators $\gamma_{p,q}(u)$.

Now, fix $q$ and suppose that $d_q=0$. Then, the asymptotic mixed moments under $\Psi_q$ cannot be obtained directly 
from free probability since all mixed moments under $\Psi$ of the type given above 
tend to zero and since $n/n_q\rightarrow \infty$, we obtain an indeterminate expression. However, let us
observe that for the family of moments considered above for various $d_{q}>0$, the limit
$$
\lim_{d_q\rightarrow 0}\Psi_q(P_{p_1}\gamma(u_1)P_{q_1}\ldots P_{p_m}\gamma(u_m)P_{q_m})
$$
always exists and is finite since the combinatorial expression of Proposition 5.2 is a polynomial in asymptotic dimensions, 
including $d_q$. Moreover, each operator $\widehat{\gamma}_{p,q}(u)$ tends (in the sense of convergence of moments) 
to the operator $\gamma_{p,q}(u)$ for any $p$. In fact, this is the main reason why we decomposed the off-diagonal operators $\widehat{\gamma}_{p,q}(u)$ as in Remark 5.2.
This is because an expression of the form $\wp_{q,p}(u)^*\wp_{q,p}(u)$ (or, equivalently, $P_{p}\wp(u)^*P_{q}\wp(u)P_{p}$) 
under $\Psi$ produces $d_q$. Therefore, the expression
$$
\lim_{d_q\rightarrow 0}\Psi_q(\widehat{\gamma}_{p_1,q_1}(u_1)\ldots \widehat{\gamma}_{p_m,q_m}(u_m))
$$
is still well defined except that we understand that $\widehat{\gamma}_{p,q}(u)=\gamma_{p,q}(u)$ for any $p,u$.
Translating that limit into the framework of random matrices, we can consider a family of 
independent Hermitian random matrices $\{Y(u,n,d_q):u\in \mathpzc{U}, n\in \mathbb{N}\}$ for any $0< d_q<\delta$, 
where $n_q/n\rightarrow d_q$. Then we obtain
$$
\lim_{d_q\rightarrow 0}\lim_{n\rightarrow \infty}\tau_q(n)(T_{p_1,q_1}(u_1,n,d_q)\ldots T_{p_m,q_m}(u_m,n,d_q))
$$
$$
=
\Psi_q(\widehat{\gamma}_{p_1,q_1}(u_1)\ldots \widehat{\gamma}_{p_m,q_m}(u_m))
$$
as $n\rightarrow \infty$, where $\widehat{\gamma}_{p,q}(u)=\widehat{\gamma}_{q,p}(u)=\gamma_{p,q}(u)$ for any $p,u$
and the remaining symmetrized operators $\widehat{\gamma}_{i,j}(u)$, $i,j\neq q$, 
do not reduce to the non-symmetrized ones. 
Therefore, we obtain the asymptotics of the desired form if we take the iterated limit. 

Let us justify the fact that the same limit is obtained when 
$n_q/n\rightarrow 0$ as $n\rightarrow \infty$.  
Let  $\{D_1, \ldots , D_r\}$ and $\{D_{1}', \ldots , D_{r}'\}$ be two different 
$n$-dependent sequences of families of diagonal matrices giving two different decompositions of the $n\times n$ identity matrices, respectively.
Let $D$ an $D'$, respectively, be the corresponding matrices of asymptotic dimensions.
We would like to show that for any $\epsilon>0$ there exists $\delta>0$ such that if 
$\max_{1\leq p\leq r}|d_{p}-d_{p}'|<\delta$, then for large $n$ it holds that
$$
|{\rm tr}(n)(D_{q_0}'Y_1D_{q_{1}}'\ldots D_{q_{m-1}}'Y_mD_{q_{m}}'-
D_{q_0}Y_1D_{q_1}\ldots D_{q_{m-1}}Y_mD_{q_{m}})|< \epsilon
$$
for any $q_0,\ldots , q_m$, where we write $Y(u_j,n)=Y_j$ for brevity.
Denote $D_j'-D_j=\Delta_{j}$ for any $j$. Then we have
$$
{\rm tr}(n)(D_{q_0}'Y_1D_{q_{1}}'\ldots D_{q_{m-1}}'Y_mD_{q_{m}}'-
D_{q_0}Y_1D_{q_{1}}\ldots D_{q_{m-1}}Y_mD_{q_{m}})
$$
$$
=\sum_{j=0}^{m}{\rm tr}(n)(D_{q_0}Y_1D_{q_{1}}\ldots Y_{j}\Delta_{q_{j}}Y_{j+1}D_{q_{j+1}}'\ldots
D_{q_{m-1}}'Y_mD_{q_{m}}').
$$
Now, we use the Schatten $p$-norms $\parallel \!A\!\parallel_{p}=\sqrt[p]{{\rm tr}(n)(|A|^p))}$ for $p\geq 1$
to get some estimates. It holds that
$|{\rm tr}(n)(A)|\leq \parallel \!A \!\parallel _1\leq \parallel \! A \!\parallel _{p}\leq \parallel \!A\!\parallel$
for any $p\geq 1$. This inequality, together with a repeated use of the H\"{o}lder inequality $\parallel \!AB \!\parallel_{r}\leq \parallel \!A \!\parallel_{p}\parallel \!B \!\parallel_{q}$, where $r^{-1}=p^{-1}+q^{-1}$, gives an upper bound for the absolute value of the RHS 
of the form
$$
(m+1) \cdot \max_{0\leq j \leq m}\parallel \!\Delta_{q_j}\!\parallel_{2m+1} \cdot \max_{1\leq j \leq m}\parallel \!Y_{j}\!\parallel_{2m+1}^{m} 
$$
since the $p$-th norm of each $D_{q}$ is bounded by $1$ for any $p$. Now, since
$$
\parallel\!\Delta_{q}\!\parallel_{2m+1}=\left|\frac{n_{q'}-n_{q}}{n}\right|^{1/(2m+1)}
$$
each $\parallel\!\Delta_{q_j}\!\parallel_{2m+1}$ can be made smaller than $\epsilon/M$ for large $n$ and small $\delta$, where
$$
M=(m+1)\cdot \max_{1\leq j \leq m}\parallel \!Y_{j}\!\parallel_{2m+1}^{m}\leq (m+1)C^{m}
$$ 
and thus our assertion follows. This means that the convergence of the considered mixed moments as $n\rightarrow \infty$ 
is uniform within the considered family of random matrices for given $\delta$ and we can take the iterated limit 
to get the limit moments for $d_{q}=0$. This completes our proof.
\hfill $\blacksquare$

\begin{Corollary}
Under the assumptions of Theorem 6.1, 
$$
\lim_{n\rightarrow \infty}
\tau(n)(T_{p_1,q_1}(u_1,n)\ldots T_{p_m,q_m}(u_m,n))
=
\Psi(\widehat{\gamma}_{p_1,q_1}(u_1)\ldots\widehat{\gamma}_{p_m,q_m}(u_m))\\ 
$$
for any $(p_1,q_1),\ldots , (p_m,q_m)\in \mathpzc{J}$ and $u_1, \ldots u_m\in \mathpzc{U}$.
\end{Corollary}
{\it Proof.}
We have
$$
\tau(n)=\frac{1}{n}\sum_{q=1}^{r}n_q\tau_{q}(n),
$$
and since $n_{q}/n\rightarrow d_{q}$ as $n\rightarrow \infty$ for any $q\in [r]$, the assertion follows from Theorem 6.1.
\hfill $\blacksquare$\\

\begin{Remark}
{\rm Let us make a few observations.
\begin{enumerate} 
\item
The above results hold for the symmetric blocks of unitarily invariant Hermitian random matrices
whose distributions converge to compactly supported measures on the real line becasue they are known to 
be asymptotically free and asymptotically free from the deterministic diagonal matrices [18, Theorem 4.3.5]
and their norms are uniformly bounded.
\item
If we rescale all variables in each block $T_{p,q}(u)$ by a constant $c_{p,q}(u)$
(each matrix of such scaling constants is assumed to be symmetric), we obtain the limit mixed moments of symmetric blocks
from Theorem 6.1 by linearity. It suffices to multiply each $\widehat{\gamma}_{p,q}(u)$ by $c_{p,q}(u)$.
\item
An important example is that of Hermitian Gaussian random matrices, where
block variances are $v_{p,q}(u)=c_{p,q}(u)^2$. The limit moments are expressed in terms of operators 
$\omega_{p,q}(u)$ with variances in the state $\Psi_q$ equal to 
$b_{p,q}(u)=d_pv_{p,q}(u)$, respectively.
\end{enumerate}
}
\end{Remark}

\begin{Definition}
{\rm By a {\it Hermitian Gaussian random matrix} (HGRM) we shall understand a
Hermitian $n$-dimensional matrix $Y(n)=(Y_{i,j}(n))$ of complex-valued random variables such that
\begin{enumerate}
\item
the variables
$$
\{{\rm Re}Y_{i,j}(n),\;{\rm Im}Y_{i,j}(n):\,1\leq i\leq j\leq n\}
$$ 
are independent Gaussian random variables,
\item
${\mathbb E}(Y_{i,j}(n))=0$ for any $i,j,n$,
\item
${\mathbb E}(|Y_{j,j}(n)|^{2})=v_{q,q}/n$ and
${\mathbb E}(|{\rm Re}Y_{i,j}(n)|^{2})={\mathbb E}(|{\rm Im}Y_{i,j}(n)|^{2})=v_{p,q}/2n$
for any $(i,j)\in N_p\times N_q$ and $n$, where $p\neq q$ and $V:=(v_{p,q})\in M_{r}({\mathbb R})$ consists of positive entries
and is symmetric.
\end{enumerate}}
\end{Definition}

\begin{Definition}
{\rm We will say that the family of HGRM of Definition 6.3 is {\it independent} if the variables
$$
\{{\rm Re}Y_{i,j}(u,n),{\rm Im}Y_{i,j}(u,n):\,1\leq i\leq j\leq n,\;u\in \mathpzc{U}\}
$$ 
are independent and we denote by $V(u)=(v_{p,q}(u))$ the corresponding variance matrices.}
\end{Definition}

In analogy to the case of one HGRM studied in [23], the asymptotic joint distributions of 
symmetric blocks of a family of HGRM are expressed in terms of the matrices
$$
B(u)=DV(u),
$$
for any $u\in \mathpzc{U}$. Let us observe that if all block variances are equal to one, then
$$
B(u)=\left(
\begin{array}{llll}
d_{1} & d_{1}  & \ldots & d_{1}\\
d_{2} & d_{2}  & \ldots & d_{2}\\
. & . & \ddots &.\\

d_{r} & d_{r}  & \ldots & d_{r}
\end{array}
\right),
$$
for any $u\in \mathpzc{U}$. Even this simple situation is interesting when we investigate joint distributions of blocks. 

\begin{Corollary}
If $\{Y(u,n):u\in \mathpzc{U}, n\in {\mathbb N}\}$ is a family of independent HGRM, then 
$$
\lim_{n\rightarrow \infty}\tau_q(n)(T_{p_1,q_1}(u_1,n)\ldots T_{p_m,q_m}(u_m,n))
=
\Psi_q(\widehat{\omega}_{p_1,q_1}(u_1)\ldots \widehat{\omega}_{p_m,q_m}(u_m))
$$
as $n\rightarrow \infty$ for any $u_1, \ldots, u_m\in \mathpzc{U}$ and $p_1,q_1, \ldots , p_m,q_m,q\in [r]$.
\end{Corollary}
{\it Proof.}
If $r_{k}(u)=\delta_{k,2}$ for any $u\in\mathpzc{U}$, we obtain $\widehat{\gamma}_{p,q}(u)=\widehat{\omega}_{p,q}(u)$ 
for any $p,q,u$. Then, we use Theorem 6.1 together with asymptotic freeness of independent HGRM under $\tau(n)$
to prove our claim (the case of arbitrary variances is also covered since we can rescale
the blocks and the operators). \hfill $\blacksquare$\\

The family of polynomials of Corollary 5.1 gives the limit mixed moments of Theorem 6.1 
in the case when all scaling constant are equal to one.
If $D$ is singular, then certain limit mixed moments of Theorem 6.1, and thus also the corresponding 
polynomials, may vanish. We close this section with giving the necessary and sufficient conditions for
them to be non-trivial. This result shows that if $d_q$ vanishes, the only interesting situation 
is when only the imaginary block is colored by $q$.

\begin{Corollary}
Under the assumptions of Theorem 6.1, if
$$
\lim_{n \rightarrow \infty}\tau_{q}(n)(T_{p_1,q_1}(u_1,n)\ldots T_{p_m,q_m}(u_m,n))\neq 0
$$
then the following conditions hold:
\begin{enumerate}
\item
there are no evanescent blocks in the above moment,
\item 
all unbalanced blocks are of the form $T_{p,q}(u,n)$ for some $p\neq q$,
\item
$\mathcal{NCI}_{m,q}((\mathpzc{w}_{1},u_1), \ldots , (\mathpzc{w}_{m},u_m))\neq \emptyset$.
\end{enumerate}
\end{Corollary}
{\it Proof.}
If $T_{p_i,q_i}(u,n_i)$ is evanescent for some $i\in [m]$, then the corresponding $\gamma_{p_i,q_i}(u_i)$ 
is trivial and thus $M_{m}=0$. Next, if $T_{p_i,q_i}(u,n_i)$ is unbalanced for some $i\in [m]$, 
where $q\notin \{p_i,q_i\}$, then, assuming (without loss of generality) that the 
corresponding symmetrized operator is of the form $\widehat{\gamma}_{p_i,q_i}(u_i)=\gamma_{p_i,q_i}(u_i)$
(and thus $d_{q_i}=0$) and taking the largest such $i$, we observe that it must act on $\Omega_{q_i}$ or on 
a vector of the form $e_{q_{i},t}(u)\otimes w$ in order to give a non-zero result. However, since $q_i\neq q$,
we must have $\Omega_{q_i}\neq \Omega_q$ and, moreover, the index $q_{i}$ is `inner' and therefore must deliver 
$d_{q_{i}}$ in association with some annihilation operator of the form $\wp_{q_{i},t}(u)$.
Since $d_{q_{i}}=0$, the moment must vanish.  Thus, each unbalanced block has to be of the form 
$T_{p,q}(u,n)$ for some $p\in [r]$ and $u\in \mathpzc{U}$. Finally, if 
$\mathcal{NCI}_{m,q}((\mathpzc{w}_{1},u_1), \ldots , (\mathpzc{w}_{m},u_m))= \emptyset$, 
then $M_{m}=0$ by Proposition 5.3. This completes the proof.
\hfill $\blacksquare$\\

\section{Asymptotic monotone independence and s-freeness}

In this section we show how to use the asymptotics of symmetric blocks under partial traces to obtain 
random matrix models for booolean, monotone and s-free independences. 
The framework of matricial freeness random variables is very effective here since the matricial 
nature of our variables allows us to choose only those with suitable indices to construct random variables which are independent 
in the required sense. Then one can identify these variables with blocks of random matrices. 
Partial traces and unbalanced blocks are particularly useful. 

We will restrict our attention to the case of two independent asymptotically free 
Hermitian random matrices decomposed in terms of three symmetric blocks as 
$$
Y(u,n)=\sum_{1\leq i\leq j\leq 2}T_{i,j}(u,n), 
$$
where $u\in \{1,2\}$. Choosing appropriate symmetric blocks of these matrices, we will construct pairs of
independent Hermitian random matrices $\{X(1,n), X(2,n)\}$ which are asymptotically 
boolean independent, monotone independent or s-free under one of the partial traces. Thus, the distributions of 
the sums
$$
X(n)=X(1,n)+X(2,n)
$$
will be asymptotically boolean, monotone or s-free convolutions of the 
asymptotic distributions of the summands. The boolean case is rather straightforward, 
but the other two are slightly more involved. 

In view of Theorem 6.1, if we are able to construct 
families of random variables which are independent in some sense from a family of arrays of canonical noncommutative random variables, 
a corresponding random matrix model is close at hand.
The definitions of boolean independence, monotone independence and s-freeness can be found in [22]. 
Let us remark that the notion of s-freeness is related to the subordination property for the 
free additive convolution [7,34].
Moreover, we recall that in the last two cases the order in which the variables appear is relevant. 

\begin{Lemma}
Let $\mathpzc{U}=\{1,2\}$, $r=2$ and $r_{1}(u)=0$ for $u\in \{1,2\}$. Then
\begin{enumerate}
\item the pair $\{\gamma_{1,2}(1), \gamma_{1,2}(2)\}$ is boolean independent with respect to $\Psi_2$,
\item the pair $\{\gamma_{1,2}(1), \gamma_{1,1}(2)+\gamma_{1,2}(2)\}$ is monotone independent with respect to $\Psi_2$,
\item the pair $\{\gamma_{1,1}(1)+\gamma_{1,2}(1), \gamma_{1,1}(2)\}$ is s-free with respect to $(\Psi_2,\Psi_1)$,
\end{enumerate}
where the operators $\gamma_{p,q}(u)$ are given by Proposition 5.1 and the covariance of $\wp_{p,q}(u)$ is $d_p$
for any $p,q,u$, where $d_1=1$ and $d_2=0$.
\end{Lemma}
{\it Proof.}
Using Proposition 5.1 and the assumptions of this lemma, we can write the operators 
under consideration in the form
\begin{eqnarray*}
\gamma_{1,1}(u)&=&\wp_{1,1}(u)^{*}+\sum_{k=1}^{\infty}r_{k+1}(u)\wp_{1,1}^{k}(u)\\
\gamma_{1,2}(u)&=&\wp_{1,2}(u)^{*}+\sum_{k=1}^{\infty}r_{k+1}(u)\wp_{1,1}^{k-1}(u)\wp_{1,2}(u)
\end{eqnarray*}
for $u\in \{1,2\}$, where we use the fact that $\wp_{2,1}(u)=\wp_{2,2}(u)=0$ since $d_2=0$.
To obtain (1), it suffices to observe that $\Omega_2$ is the 
only basis vector that can appear in the range of a polynomial in $\gamma_{1,2}(1)$ ($\gamma_{1,2}(2)$) onto which $\gamma_{1,2}(2)$ ($\gamma_{1,2}(1)$)
acts non-trivially. The same is true if we interchange both operators. To prove (2), denote 
$$
\gamma(1)=\gamma_{1,2}(1)\;\;{\rm and}\;\;
\gamma(2)=\gamma_{1,1}(2)+\gamma_{1,2}(2).
$$
We need to show that the pair $\{\gamma(1),\gamma(2)\}$ (in that order) is monotone independent w.r.t. $\Psi_2$, i.e.
$$
\Psi_2(w_1a_1b_1a_2w_2)=\Psi_2(b_1)\Psi_2(w_1a_1a_2w_2)
$$
for any $a_1,a_2\in {\mathbb C}[\gamma(1),1_1]$, $b_1\in {\mathbb C}[\gamma(2),1_2]$, where $1_1=1_{1,2}$ and $1_{2}=1$  and $w_1,w_2$ are arbitrary
elements of ${\mathbb C}\langle \gamma(1),\gamma(2),1_1,1_2\rangle$.
It suffices to consider the action of $\gamma(1)$ and $\gamma(2)$ onto their invariant subspace in ${\mathcal M}$
of the form
$$
{\mathcal M}'={\mathbb C}\Omega_2\oplus {\mathcal F}({\mathcal H}_{1,1}(2))\otimes {\mathcal F}({\mathcal H}_{1,1}(1))\otimes 
({\mathcal H}_{1,2}(2) \oplus{\mathcal H}_{1,2}(1))
$$
where ${\mathcal F}({\mathcal H})$ is the free Fock space over ${\mathcal H}$ with the vacuum vector $\Omega$ and
${\mathcal H}_{p,q}(u)={\mathbb C}_{p,q}(u)$ for any $p,q,u$ (we identify $\Omega\otimes e_{1,2}(u)$ with $e_{1,2}(u)$). 
Now, the range of any polynomial in $\gamma(1)$ is 
contained in 
$$
{\mathcal M}_{1}'={\mathbb C}\Omega_2\oplus {\mathcal F}({\mathcal H}_{1,1}(1))\otimes {\mathcal H}_{1,2}(1)
$$ 
since the only vector in ${\mathcal M}'$ onto which $\wp_{1,2}(1)$ acts nontrivially is $\Omega_2$, giving $e_{1,2}(1)$,
and then the action of powers of $\wp_{1,1}(1)$ gives the remaining basis vectors of ${\mathcal M}_{1}'$. 
Of course, $1_1$ preserves any vector from ${\mathcal M}_{1}'$. Therefore, it suffices to compute the action of
any polynomial in $\gamma(2)$ onto any vector from ${\mathcal M}_{1}'$. Now, looking at the formulas for $\gamma_{1,1}(2)$
and $\gamma_{1,2}(2)$, we observe that the action of powers of $\gamma(2)$ 
onto any basis vector of ${\mathcal M}_{1}'$ is the same as the action of the canonical noncommutative random variable 
onto the vacuum vector in the free Fock space. Thus, we have
$$
(\gamma(2))^{k}v=\alpha_{k}v+{\rm higher \;order\;terms}
$$
for any $k\geq 1$ and any basis vector $v$ of ${\mathcal M}'$. In particular,
\begin{eqnarray*}
(\gamma(2))^{k}\Omega_2&=&\alpha_{k}\Omega_2\;{\rm mod}\;({\mathcal M}'\ominus {\mathbb C}\Omega_2)\\
(\gamma(2))^{k}e_{1,2}(1)&=&\alpha_{k}e_{1,2}(1)\;{\rm mod}\;({\mathcal M}'\ominus ({\mathbb C}\Omega_2\oplus {\mathcal H}_{1,2}(1))),
\end{eqnarray*}
etc. Thus $\Psi_2((\gamma(2))^{k})=\alpha_k$ and, moreover, since all higher order terms are
in  ${\rm Ker}(\gamma(1))$, we can pull out $\alpha_{k}$ and the required condition for monotone independence holds for positive powers of $\gamma(2)$.
It is easy to see that it also holds if $b_1=1_2$, which completes the proof of (2). To prove (3), consider 
$$
\delta(1)=\gamma_{1,1}(1)+\gamma_{1,2}(1), \;\;\delta(2)=\gamma_{1,1}(2)
$$
and their invariant subspace in ${\mathcal M}$ of the form
$$
{\mathcal M}''={\mathbb C}\Omega_2\oplus {\mathcal F}({\mathcal H}_{1,1}(1)\oplus {\mathcal H}_{1,1}(2))\otimes {\mathcal H}_{1,2}(1).
$$
Observe that ${\mathcal M}''$ is isomorphic to 
$$
{\mathcal F}_{1,2}={\mathbb C}\Omega\oplus {\mathcal F}_{1}^{0}\oplus ({\mathcal F}_{2}^{0}\otimes {\mathcal F}_{1}^{0})\oplus
({\mathcal F}_{1}^{0}\otimes {\mathcal F}_{2}^{0}\otimes {\mathcal F}_{1}^{0})\oplus \ldots
$$
where ${\mathcal F}_{j}^{0}={\mathcal F}({\mathbb C}e(j))\ominus {\mathbb C}\xi_{j}$ and the isomorphism $\tau$ is given by
$\tau(\Omega_{2})=\Omega$ and 
$$
\tau(e_{1,1}(j_1)\otimes \ldots \otimes e_{1,1}(j_{k-1})\otimes e_{1,2}(1))=e(j_1)\otimes \ldots e(j_{k-1})\otimes e(1)
$$
for any $j_1,\ldots , j_{k-1}\in \{1,2\}$ and $k\in {\mathbb N}$. The space ${\mathcal F}_{1,2}$ is the s-free product of two free Fock spaces,
$({\mathcal F}_{1},\xi_1)$ and $({\mathcal F}_{2}, \xi_2)$, the subspace of their free product 
$({\mathcal F}_{1},\xi_1)*({\mathcal F}_{2}, \xi_2)$ (for the definition of the s-free product of Hilbert spaces and of the 
s-free convolution describing the free subordination property, see [21]). Now, observe that the action of $\delta(1)$ onto 
${\mathcal M}^{''}$ can be identified with the action of the canonical noncommutative random variable $\gamma_1$ (built 
from $\ell_1$ and its adjoint) restricted to ${\mathcal F}_{1,2}$. Similarly, the action of $\delta(2)$ onto ${\mathcal M}^{''}$ can be identified with the action of $\gamma_2$ (built from $\ell_2$ and its adjoint) restricted to ${\mathcal F}_{1,2}$. This proves that the pair $\{\delta(1),\delta(2)\}$ (in that order) is s-free with respect to the pair of states $(\Psi_2,\Psi_1)$, which gives (3).
\hfill $\blacksquare$\\

We are ready to state certain results concerning asymptotic properties of 
independent Hermitian random matrices $X(1,n), X(2,n)$ which are built from at most three symmetric blocks
of two independent matrices $Y(1,n),Y(2,n)$ which satisfy the assumptions of Theorem 6.1.
We assume that one asymptotic dimension is equal to zero, thus 
one block is balanced, one is unbalanced and one is evanescent. 

\begin{Theorem}
Under the assumptions of Theorem 6.1, let $\mathpzc{U}=\{1,2\}$, where $r_{1}(u)=0$ for $u\in \{1,2\}$, 
$r=2$, $d_2=0$ and $d_1=1$.
Then
\begin{enumerate}
\item
the pair $\{T_{1,2}(1,n), T_{1,2}(2,n)\}$ is asymptotically boolean independent with respect to $\tau_{2}(n)$,
\item
the pair $\{T_{1,2}(1,n),Y(2,n)\}$ is asymptotically monotone independent with respect to $\tau_2(n)$,
\item
the pair $\{Y(1,n), T_{1,1}(2,n)\}$ is asymptotically s-free with respect to $(\tau_2(n), \tau_1(n))$.
\end{enumerate}
\end{Theorem}
{\it Proof.}
Since $d_1=1$ and $d_2=0$, blocks $T_{1,1}(u,n), T_{1,2}(u,n)$ and $T_{2,2}(u,n)$ are balanced, unbalanced and evanescent, respectively. 
Therefore, by Theorem 6.1, we have convergence as $n\rightarrow \infty$ 
\begin{eqnarray*}
T_{1,2}(u,n)&\rightarrow &\gamma_{1,2}(u)\\
T_{1,1}(u,n)&\rightarrow &\gamma_{1,1}(u)\\
T_{2,2}(u,n)&\rightarrow &0
\end{eqnarray*}
in the sense of moments under $\tau_2(n)$ for $u\in \{1,2\}$, 
where the limit moments are computed in the state $\Psi_2$. By Lemma 7.1, the proof
is completed. \hfill $\blacksquare$\\

\section{Non-Hermitian Gaussian random matrices}

We would like to generalize the results of Section 6 to the ensemble of 
independent non-Hermitian random matrices. For simplicity, we restrict our attention to
the ensemble of non-Hermitian Gaussian matrices. This ensemble (or, the family of its entries) 
is sometimes called the Ginibre Ensemble (the family of symmetric blocks of the Ginibre Ensemble
can be called the Ginibre Symmetric Block Ensemble). We keep the same settings for blocks as in Section 6.

We will show that the Ginibre Symmetric Block Ensemble converges in *-moments as $n\rightarrow \infty$ 
to the ensemble of non-self adjoint operators 
$$
\{\eta_{p,q}(u):\, (p,q)\in \mathpzc{J},  u\in \mathpzc{U}\}
$$
where 
$$
\eta_{p,q}(u)=\widehat{\wp}_{p,q}(2u-1)+\widehat{\wp}^{*}_{p,q}(2u)
$$
for $(p,q)\in \mathpzc{J}$, $u\in \mathpzc{U}=[t]$ and $(\widehat{\wp}_{p,q}(u))$, 
where $u\in [2t]$, are arrays of symmetrized creation operators (the set $\mathpzc{J}$ is assumed to be symmetric).
In fact, it suffices to consider the case when $\mathpzc{J}=[r]\times [r]$ and restrict to 
a symmetric subset if needed.

In order to define $t$ arrays $(\eta_{p,q}(u))$, where $u\in [t]$, we need $2t$ arrays 
$(\widehat{\wp}_{p,q}(u))$, where $u\in [2t]$. Moreover, we shall assume that 
$$
b_{p,q}(2u-1)=b_{p,q}(2u)=d_{p}v_{p,q}(u)
$$ 
for any fixed $p,q,u$. Our definition parallels that in the free case, where operators
$$
\eta(u)=\ell(2u-1)+\ell^{*}(2u)
$$
for $u\in [t]$, which are unitarily equivalent to circular operators, 
are introduced [31, Theorem 3.3]. Here, $\{\ell(1),\ldots , \ell(2t)\}$ is a family of free creation operators.

We can now generalize the result of Corollary 6.2 to the Ginibre Symmetric Block Ensemble.
The important difference is that we do not assume that the matrices are Hermitian and thus each matrix 
contains $2n^2$ independent Gaussian random variables. 
However, in order to reduce this model to the Hermitian case, we need to assume that the variance matrices are symmetric. 
The theorem given below is a block refinement of that proved by Voiculescu for the Ginibre Ensemble
[33, Theorem 3.3] except that we assume the Gaussian variables to be independent block-identically distributed ({\it i.b.i.d.})
instead of i.i.d.

\begin{Theorem}
Let $(Y(u,n))$ be the family of complex random matrices, where $n\in \mathbb{N}$, $u\in \mathpzc{U}=[t]$, and 
$$
\{{\rm Re}Y_{i,j}(u,n), {\rm Im}Y_{i,j}(u,n): i,j\in [n], u\in \mathpzc{U}\}
$$
is a family of independent Gaussian random variables with zero mean and variances 
$$
{\mathbb E}(({\rm Re}Y_{i,j}(u,n))^{2})= {\mathbb E}(({\rm Im}Y_{i,j}(u,n))^{2})=\frac{v_{p,q}(u)}{2n}
$$ 
for any $(i,j)\in N_{p,q}$ and $u\in \mathpzc{U}$, where the matrices $(v_{p,q}(u))$ are symmetric. Then 
$$
\lim_{n\rightarrow \infty}\tau_{q}(n)(T_{p_1,q_1}^{\,\epsilon_1}(u_1,n)\ldots T_{p_m,q_m}^{\,\epsilon_m}(u_m,n))
=
\Psi_{q}(\eta_{p_1,q_1}^{\epsilon_1}(u_1)\ldots \eta_{p_m,q_m}^{\epsilon_m}(u_m))
$$
for any $\epsilon_1, \ldots , \epsilon_m\in \{1,*\}$, $u_1, \ldots, u_m\in \mathpzc{U}$ and
$q, p_1,q_1, \ldots , p_m,q_m\in [r]$.
\end{Theorem}
{\it Proof.}
The proof is similar to that in [33, Theorem 3.3].  
Write each matrix in the form
$$
Y(u,n)=2^{-1/2}\left(X(u,n)+i\widetilde{X}(u,n)\right),
$$
where
\begin{eqnarray*}
X(u,n)&=&2^{-1/2}\left(Y(u,n)+Y^{*}(u,n)\right)\\
\widetilde{X}(u,n)&=&i2^{-1/2}\left(Y^{*}(u,n)-Y(u,n)\right)
\end{eqnarray*}
are Hermitian for any $n\in \mathbb{N}$ and $u\in \mathpzc{U}$. 
The symmetric blocks of these matrices, denoted $(U_{p,q}(u,n))$ and $(\widetilde{U}_{p,q}(u,n))$, respectively,
will give the asymptotics of symmetric blocks $(T_{p,q}(u,n))$ and their adjoints since
$$
T_{p,q}(u,n)=2^{-1/2}(U_{p,q}(u,n)+i\widetilde{U}_{p,q}^{*}(u,n))
$$
for any $p,q,n,u$. They are built from variables $X_{i,j}(u,n)$ and $\widetilde{X}_{i,j}(u,n)$, where
\begin{eqnarray*}
{\rm Re}X_{i,j}(u,n)&=&2^{-1/2}\left({\rm Re}Y_{i,j}(u,n)+{\rm Re}Y_{j,i}(u,n)\right)\\
{\rm Im}X_{i,j}(u,n)&=&2^{-1/2}\left({\rm Im}
Y_{i,j}(u,n)-{\rm Im}Y_{j,i}(u,n)\right)\\
{\rm Re}\widetilde{X}_{i,j}(u,n)&=&2^{-1/2}\left({\rm Im}Y_{i,j}(u,n)+{\rm Im}Y_{j,i}(u,n)\right)\\
{\rm Im}\widetilde{X}_{i,j}(u,n)&=&2^{-1/2}\left({\rm Re}Y_{j,i}(u,n)-{\rm Re}Y_{i,j}(u,n)\right),
\end{eqnarray*}
for any $(i,j)\in N_{p,q}$ and any $u\in \mathpzc{U}$.
These satisfy the assumptions of Corollary 6.2. In particular, they are independent due to the fact that 
for fixed $i\neq j $ and $u,n$, the pairs $\{{\rm Re}Y_{i,j}(u,n), {\rm Re}Y_{j,i}(u,n)\}$ and 
$\{{\rm Im}Y_{i,j}(u,n),{\rm Im}Y_{j,i}(s,n)\}$ are identically distributed since the variance matrices 
$(v_{p,q}(u))$ are symmetric by assumption (each variable $X_{i,j}(u)$ indexed by $(i,j)\in N_{p,q}$ has
variance $v_{p,q}(u)/n$). Denote 
$$ 
U_{p,q}(u,n)=Z_{p,q}(2u-1,n)\;\;{\rm and}\;\;\widetilde{U}_{p,q}(u,n)=Z_{p,q}(2u,n)
$$
for any $p,q\in [r]$ and $u\in \mathpzc{U}$.  
Using Corollary 6.2, we can express 
the asymptotic mixed moments in the blocks $(Z_{p,q}(u,n))$ in terms of mixed moments 
in symmetrized Gaussian operators $(\widehat{\omega}_{p,q}(u))$, where $u\in [2t]$, namely
$$
\lim_{n\rightarrow \infty}
\tau_{q}(n)(Z_{p_1,q_1}(u_1,n)\ldots Z_{p_m,q_m}(u_m,n))
=
\Psi_{q}(\widehat{\omega}_{p_1,q_1}(u_1)\ldots\widehat{\omega}_{p_m,q_m}(u_m))\\ 
$$
for any $(p_1,q_1),\ldots , (p_m,q_m)\in \mathpzc{J}$, $q\in [r]$, and $u_1, \ldots u_m\in [2t]$. 
Each of these arrays has semicircle distributions on the diagonal and Bernoulli distributions 
elsewhere. These arrays are independent in the sense discussed in Section 3. Now, the linear map
$$
\theta: M_{r}^{2t}({\mathbb C})\rightarrow M_{r}^{2t}({\mathbb C})
$$
such that 
\begin{eqnarray*}
\theta(e_{p,q}(2u-1))&=&2^{-1/2}(e_{p,q}(2u-1)+ie_{p,q}(2u))\\
\theta(e_{p,q}(2u))&=&2^{-1/2}(e_{p,q}(2u-1)-ie_{p,q}(2u))
\end{eqnarray*}
for any $p,q\in [r]$ and $u\in [t]$, induces a unitary map
$$
{\mathcal M}(\theta):\; \mathcal{M}(\widehat{\mathcal H})\rightarrow \mathcal{M}(\widehat{\mathcal H})
$$
where $\mathcal{M}(\widehat{\mathcal{H}})$ is the matricially free Fock space of tracial type over 
the array $\widehat{\mathcal{H}}=({\mathcal H}_{p,q})$ of Hilbert spaces
$$
{\mathcal H}_{p,q}=\bigoplus_{u=1}^{2t}{\mathbb C}_{p,q}(u),
$$
where $p,q\in [r]$, for which it holds that
\begin{eqnarray*}
{\mathcal M}(\theta)\Omega_q&=&\Omega_q,\\
{\mathcal M}(\theta)\wp_{p,q}(2u-1){\mathcal M}(\theta)^{*}&=&2^{-1/2}(\wp_{p,q}(2u-1)+i\wp_{p,q}(2u)),\\
{\mathcal M}(\theta)\wp_{p,q}(2u){\mathcal M}(\theta)^{*}&=&2^{-1/2}(\wp_{p,q}(2u-1)-i\wp_{p,q}(2u))
\end{eqnarray*}
for any $p,q\in [r]$ and $u\in [t]$. Consequently,
$$
{\mathcal M}(\theta)\eta_{p,q}(u){\mathcal M}(\theta)^{*}=
2^{-1/2}(\widehat{\omega}_{p,q}(2u-1)+i\widehat{\omega}_{p,q}(2u))
$$
for any $p,q\in [r]$ and $u\in [t]$. Therefore, each mixed moment in operators (which replace
the circular operators)
$$
\widehat{\omega}_{p,q}(2u-1)+i\widehat{\omega}_{p,q}(2u)
$$
can be expressed as a corresponding moment in the operators $\eta_{p,q}(u)$, 
where $p,q\in [r]$ and $u\in [t]$, respectively. 
This completes the proof. \hfill $\blacksquare$\\

\begin{Remark}
{\rm Let us observe that the family of arrays of operators
$$
\{[\zeta(u)]:\,u\in \mathpzc{U}\}
$$
where each entry of $[\zeta(u)]$ is of the form
$$
\zeta_{p,q}(u)=\widehat{\omega}_{p,q}(2u-1)+i\widehat{\omega}_{p,q}(2u),
$$
for any $(p,q)\in \mathpzc{J}$, replaces the family of circular operators of Voiculescu [33]
related to the circular law [3,8,14]. 
In particular, it is easy to see that each operator $\zeta_{p,q}(u)$ is a circular operator 
since it is of the form $a+ib$, where $a$ and $b$ are free with respect to the corresponding 
state $\Psi_q$ and have semicircle distributions under $\Psi_q$. If $p=q$, the distributions of $a$ and $b$ are
identical, whereas if $p\neq q$, they are identical if and only if $d_p=d_q$.
}
\end{Remark}

Theorem 8.1 can be applied to the study of Wishart matrices and more general products of complex 
rectangular random Gaussian matrices. Let us state a combinatorial formula for the mixed moments
of the operators from the arrays $(\eta_{p,q}(u))$, where $u\in \mathpzc{U}$. First, however, 
let us present an example with explicit computations.
\begin{Example}
{\rm For simplicity, assume that $t=1$ and denote $\eta_{p,q}=\eta_{p,q}(1)$ for any $p,q\in [r]$. Of course,
$$
\eta_{p,q}=\widehat{\wp}_{p,q}(1)+\widehat{\wp}^{*}_{p,q}(2)\;\;{\rm and}\;\;
\eta_{p,q}^{*}=\widehat{\wp}_{p,q}^{*}(1)+\widehat{\wp}_{p,q}(2)
$$
For any $q\in [r]$, we have
\begin{eqnarray*}
\Psi_q((\eta_{p,q}\eta_{p,q}^{*})^{3})
&=&\Psi_{q}(\wp_{p,q}^{*}(2)\wp_{p,q}(2)\wp_{p,q}^{*}(2)\wp_{p,q}(2)\wp_{p,q}^{*}(2)\wp_{p,q}(2))\\
&+&
\Psi_{q}(\wp_{p,q}^{*}(2)\wp_{p,q}(2)\wp_{p,q}^{*}(2)\wp_{q,p}^{*}(1)\wp_{q,p}(1)\wp_{p,q}(2))\\
&+&
\Psi_{q}(\wp_{p,q}^{*}(2)\wp_{q,p}^{*}(1)\wp_{q,p}(1)\wp_{p,q}(2)\wp_{p,q}^{*}(2)\wp_{p,q}(2))\\
&+&
\Psi_{q}(\wp_{p,q}^{*}(2)\wp_{q,p}^{*}(1)\wp_{q,p}(1)\wp_{q,p}^{*}(1)\wp_{q,p}(1)\wp_{p,q}(2))\\
&+&
\Psi_{q}(\wp_{p,q}^{*}(2)\wp_{q,p}^{*}(1)\wp_{p,q}^{*}(2)\wp_{p,q}(2)\wp_{q,p}(1)\wp_{p,q}(2))\\
&=&
b_{p,q}^{3}(2)+3b_{p,q}^{2}(2)b_{q,p}(1)+b_{p,q}(2)b_{q,p}^{2}(1)
\end{eqnarray*}
and thus the summands correspond to noncrossing colored partitions shown in Example 4.5.
In the random matrix setting of Theorem 8.1, if all block variances are set to one, 
the expressions on the right-hand sides reduce to polynomials in $d_p,d_q$ of the form
$$
d_p, \;\;\;d_{p}^{\,2}+d_{p}d_{q}, \;\;\; d_{p}^{\,3}+3d_{p}^{\,2}d_{q}+d_{p}d_{q}^{\,2}
$$
respectively, since $b_{p,q}(k)=d_{p}$ and $b_{q,p}(k)=d_{q}$ for $k\in \{1,2\}$ (cf., for instance, [16]).
}
\end{Example}

The above example shows the connection between mixed moments of the considered operators 
and noncrossing colored partitions. The main feature of the combinatorics for the non-Hermitian case
is that in order to get a pairing between two operators, one of them has to be starred and one unstarred. This leads to the following 
definition and consequently, to an analog of Proposition 4.3.

\begin{Definition}
{\rm We say that $\pi \in \mathcal{NC}_{m}^{2}$ is {\it adapted}
to $((\mathpzc{w}_{1},u_1, \epsilon_1), \ldots,(\mathpzc{w}_{m},u_m, \epsilon_m))$, where 
$\mathpzc{w}_{k}=\{p_k,q_k\}$ and $(p_k,q_k,u_k)\in [r]\times [r]\times [t]$ and $\epsilon_k\in \{1,*\}$
for any $k$, if there exists a tuple $((\mathpzc{v}_{1},u_1), \ldots,(\mathpzc{v}_{m},u_m))$, where 
$\mathpzc{v}_k\in \{(p_k,q_k), (q_k,p_k)\}$ for any $k$, to which $\pi$ is adapted and $\epsilon_j\neq \epsilon_k$
whenever $\{j,k\}$ is a block. The set of such partitions will be denoted by 
$\mathcal{NC}_{m}^{2}((\mathpzc{w}_1,u_1, \epsilon_1), \ldots , (\mathpzc{w}_{m},u_m, \epsilon_m))$.
Its subset consisting of those partitions which are in 
$\mathcal{NC}_{m,q}^{2}(\mathpzc{w}_1,u_1), \ldots , (\mathpzc{w}_{m},u_m))$ will be denoted 
$\mathcal{NC}_{m,q}^{2}((\mathpzc{w}_1,u_1, \epsilon_1), \ldots , (\mathpzc{w}_{m},u_m, \epsilon_m))$.}
\end{Definition}

\begin{Proposition}
For any tuple $((\mathpzc{w}_{1},u_1, \epsilon_1), \ldots , (\mathpzc{w}_m,u_m, \epsilon_m))$ and $q\in [r]$, $m\in {\mathbb N}$,
where $\mathpzc{w}_{k}=\{p_k,q_k\}$ and $p_k,q_k\in [r]$, $u_k\in [t]$ and $\epsilon_k\in \{1,*\}$ for each $k$, it holds that
$$
\Psi_{q}(\eta_{p_1,q_1}^{\epsilon_1}(u_1)\ldots \eta_{p_m,q_m}^{\epsilon_m}(u_m))
=\sum_{\pi\in \mathcal{NC}_{m,q}^{2}((\mathpzc{w}_{1},u_1,\epsilon_1), \ldots , (\mathpzc{w}_m,u_m,\epsilon_m))}
b_{q}(\pi,f,g)
$$
where $f$ is the coloring of $\pi$ defined by $(\mathpzc{w}_{1}, \ldots , \mathpzc{w}_m;q)$ 
and $g$ is the labeling of $\pi$ defined by $(u_1, \ldots , u_m)$.
\end{Proposition}
{\it Proof.}
The above mixed moment in the operators $\eta_{p_k,q_k}^{\epsilon_{k}}(u_k)$, where $k\in [m]$, can be expressed as the mixed 
moment of the operators 
$$
2^{-1/2}(\widehat{\omega}_{p_k,q_k}(2u_k-1)\pm i\widehat{\omega}_{p_k,q_k}(2u_k)),
$$
where $k\in [m]$, and thus can be written in terms of the mixed moments of the operators $(\widehat{\omega}_{p_k,q_k}(s_k))$, where 
$s_k\in \{2u_{k}-1, 2u_{k}\}$ for any $k\in [m]$.
These, in turn, are given by Proposition 4.3. However, the associated noncrossing pair 
partitions must be adapted to $((\mathpzc{w}_{1}, u_1), \ldots , (\mathpzc{w}_{m}, u_m))$ 
since $2u_j-1=2u_k-1$ or $2u_j=2u_k$ for some $j\neq k$ implies that $u_j=u_k$.
The contributions from such partitions are given by numbers $b_{q}(\pi,f)$. 
In the context of complex Gaussian operators of the type shown above, where the variances
associated with $2u_k-1$ and $2u_k$ coincide for any fixed $k$, these numbers are obtained when computing
the moments given by the products of the above complex linear combinations (one with the plus sign and the other with the minus sign), 
namely 
$$
2^{-1}(\widehat{\omega}_{p_k,q_k}(2u_k-1)\widehat{\omega}_{p_k,q_k}(2u_k-1)+\widehat{\omega}_{p_k,q_k}(2u_k)\widehat{\omega}_{p_k,q_k}(2u_k))
$$  
under the state $\Psi_{p_k}$ or $\Psi_{q_k}$, which gives $b_{q_k,p_k}=d_{q_k}v_{q_k,p_k}(u_k)$ or $b_{p_k,q_k}=d_{p_k}v_{p_k,q_k}(u_k)$, respectively. 
This implies that the associated partition must satisfy
the additional condition $\epsilon_j=\epsilon_k$ whenever $\{j,k\}$ is a block. This completes the proof.
\hfill $\blacksquare$\\

\section{Wishart matrices and related products}
We can also study products of Gaussian random matrices using the framework of matricial freeness.
Since the main objects in our model are blocks of Gaussian random matrices, there is no major 
difference in treating rectangular or square matrices. 

The main point of this section is to show that the framework of matricial freeness is general enough to include
square or rectangular Gaussian random matrices as well as their products in one unified scheme. 
However, it can also be seen that our formalism is much more general and should give a number of interesting examples, 
of which sums and products of independent block-identically distributed Gaussian random matrices are just 
the most natural ones. 

Using matricially free probability, we can reproduce certain results concerning 
products $B(n)$ of rectangular Gaussian random matrices and the asymptotic distributions of $B(n)B^{*}(n)$ 
under the trace composed with classical expectation as $n\rightarrow \infty$. 
In particular, if $B(n)$ is just one complex Gaussian random matrix, the matrix 
$$
W(n)=B(n)B^{*}(n)
$$ 
is the {\it complex Wishart matrix} [37]. 
The limit distribution of a sequence of such matrices is the Marchenko-Pastur distribution (which also plays the role
of the free Poisson distribution) with moments given by Catalan numbers. 

The original result on the asymptotic distributions of the Wishart matrices is due to Marchenko and Pastur [27], 
but many different proofs have been given (see, for instance [15,27,35]). 
In the case when $B(n)$ is a power of a square random matrix, then it has recently been shown by Alexeev 
{\it et al} [1] that the limit moments are given by Fuss-Catalan numbers 
(products of independent random matrices have also been studied recently [2,9,10]).
The distributions defined by the Fuss-Catalan numbers were explicitly determined by Penson and $\dot{\rm Z}$yczkowski [31].
A random matrix model based on a family of distributions called free Bessel laws constructed from the Marchenko-Pastur distribution
by means of free convolutions was given by Banica {\it et al} [4].
In turn, asymptotic freeness of independent Wishart matrices 
was proved by Voiculescu [33] (see also [11,18]). Let us also mention that explicit formulas for 
moments of Wishart matrices were given, for instance, by Hanlon {\it et al} [17] and by Haagerup and Thorbjornsen [16]. 

In our framework, it suffices to use one symmetric off-diagonal (rectangular or square, balanced or unbalanced) block to study the asymptotic
distribution of the Wishart matrix. More generally, using matricially free Gaussian operators,
we can also study the asymptotic distributions of $W(n)$ in the case when the matrices $B(n)$ are sums of independent rectangular 
Gaussian random matrices as in [6], namely
$$
B(n)=Y(u_1,n)+Y(u_2,n)+\ldots +Y(u_m,n),
$$
assuming they have the same dimensions, as well as in the case when $B(n)$ is a product of independent rectangular
Gaussian random matrices
$$
B(n)=Y(u_1,n)Y(u_2,n)\ldots Y(u_k,n),
$$
assuming their dimensions are such that the products are well-defined. Using Theorem 8.1,
one can study non-Hermitian Wishart matrices which are also of interest [19]. 
In both cases, it suffices to take a suitable sequence of the off-diagonal symmetric blocks of a family of independent Gaussian 
random matrices (the first case) or even one Gaussian random matrix matrix (the second case).
Since we know how to compute the asymptotic joint distributions of symmetric blocks, we can immediately obtain results 
concerning asymptotic distributions of such products. 

In our treatment of products, it will be convenient to use sets $\mathcal{NC}_{m}^{2}(W_k)$ of noncrossing $W_k$-pairings, where $W_k$ is a word of the form
$$
W_{k}=(12 \ldots p p^{*}\ldots 2^{*}1^{*})^{k}
$$
where $m=kp$, by which we understand the set of noncrossing pair partitions of the set $[m]$ in which all blocks 
are associated with pairs of letters of the form $\{j,j^{*}\}$.
If $\pi\in \mathcal{NC}_{m}^{2}(W_k)$, we denote by $\mathcal{R}_{j}(\pi)$ and $\mathcal{R}_{j}^{*}(\pi)$ the sets of the 
right legs of $\pi$ which are associated with $j$ and $j^{*}$, respectively. Of course,
$$
{\mathcal R}(\pi)=\bigcup_{j=1}^{p}{\mathcal R}_{j}(\pi)\cup {\mathcal R}_{j}^{*}(\pi)
$$
for any $\pi\in \mathcal{NC}_{m}^{2}(W_k)$. 

\begin{Definition}
{\rm Define homogenous polynomials in variables $d_1,d_2, \ldots , d_{p+1}$ of the form
$$
P_{k}(d_1,d_2, \ldots , d_{p+1})=\sum_{\pi \in \mathcal{NC}_{kp}^{2}(W_{k})}d_{1}^{\,r_1(\pi)}d_{2}^{\,r_2(\pi)}\ldots d_{p+1}^{\,r_{p+1}(\pi)}
$$
for any $k,p\in \mathbb{N}$, where
$$
r_j(\pi)=|\mathcal{R}_{j}(\pi)|+|\mathcal{R}_{j-1}^{*}(\pi)|
$$
for any $\pi\in \mathcal{NC}_{kp}^{2}(W)$ and $1\leq j\leq p+1$, where we set $\mathcal{R}_{0}^{*}(\pi)=\emptyset$ and 
$\mathcal{R}_{p+1}(\pi)=\emptyset$. }
\end{Definition}

Our random matrix model for products of independent Gaussian random matrices will be based 
on products of a chain of symmetric random blocks 
$$
T_{j,j+1}(n)=T_{j,j+1}(u,n),\;{\rm where}\; j\in [p]\; {\rm and}\;n\in \mathbb{N},
$$
with $u\in \mathpzc{U}$ fixed and thus omitted in our notations (taking different $u$'s does not change
the computations since these blocks ore independent anyway).
The blocks are embedded in a Gaussian random matrix $Y(n)$ with 
asymptotic dimensions $d_1,d_2, \ldots , d_{p+1}$ which will be the variables of our polynomials.
The symmetric blocks $T_{j,j+1}(n)$ are not assumed to be Hermitian, thus we shall use
Theorem 8.1 when discussing their singular limit distributions. Using products of 
Hermitian random blocks is also possible, especially in the simplest case of Wishart matrices.
However, it is more convenient to use products of non-Hermitian symmetric random blocks
to give a combinatorial description of the limit moments in the case of singular limit distributions
of higher-order products.

In the proof given below, we use ideas on the enumeration of noncrossing pairings given in [20] 
(related to the results in [28]).

\begin{Theorem}
Under the assumptions of Theorem 8.1, suppose that $v_{j,j+1}=1$ for $j\in [p]$, where $p$ is a fixed 
natural number, and let
$$
B(n)=T_{1,2}(n)T_{2,3}(n)\ldots T_{p,p+1}(n)
$$
for any natural $n$. Then, for any nonnegative integer $k$,
$$
\lim_{n \rightarrow \infty}\tau_{1}(n)\left(\left(B(n)B^{*}(n)\right)^{k}\right)=P_{k}(d_1,d_2, \ldots , d_{p+1})
$$
where $d_1,d_2, \ldots , d_{p+1}$ are asymptotic dimensions.
\end{Theorem}
{\it Proof.}
By Theorem 8.1, we obtain
$$
\lim_{n\rightarrow \infty}\tau_{1}(n)((B(n)B^{*}(n))^{k})=\Psi_{1}((\eta\eta^{*})^{k})
$$
for any natural $k$, where
$$
\eta=\eta_{1,2}\eta_{2,3}\ldots \eta_{p,p+1}
$$
for any given $p$. In view of Proposition 8.1, we can express the right-hand side in 
terms of noncrossing partitions
$$
\Psi_{1}((\eta\eta^{*})^{k})=
\sum_{\pi\in \mathcal{NC}_{m,1}^{2}((\mathpzc{w}_{1},\epsilon_1), \ldots , (\mathpzc{w}_{m},\epsilon_{m}))}
b_{1}(\pi,f)
$$
where $m=kp\,$ and 
\begin{eqnarray*}
(\mathpzc{w}_{1+2lp}, \ldots, \mathpzc{w}_{p+2lp})&=&(\{1,2\}, \ldots , \{p,p+1\})\\
(\mathpzc{w}_{p+1+2lp},\ldots , \mathpzc{w}_{2p+2lp})&=&(\{p,p+1\}, \ldots , \{1,2\})
\end{eqnarray*}
for any $0\leq l \leq k-1$. In the notation for noncrossing pair partitions adapted to $((\mathpzc{w}_{1},\epsilon_1), \ldots , (\mathpzc{w}_{m},\epsilon_{m}))$ we omit $u_1, \ldots , u_{m}$ since in the considered case we have one array 
$(\eta_{p,q})$. Since only operators with the same matricial indices can form a pairing,
$\mathcal{NC}_{m,1}^{2}((\mathpzc{w}_{1},\epsilon_1), \ldots , (\mathpzc{w}_{m},\epsilon_{m}))$ 
can be put in one-to-one correspondence with the 
set $\mathcal{NC}_{kp}^{2}(W_k)$. If now $\pi\in \mathcal{NC}_{kp}^{2}(W_k)$, 
then we will identify $b_{1}(\pi,f)$ with the corresponding expression 
of Proposition 8.1. Recall that
$$
b_{1}(\pi,f)=\prod_{i=1}^{kp}b_{1}(\pi_i,f)
$$
where $\{\pi_1, \ldots , \pi_{kp}\}$ are blocks of $\pi$ and where the coloring $f$ is uniquely determined by $\pi$ 
and by the coloring $q=1$ of the imaginary block.
Since the block variances are set to one, we can express each 
$b_{1}(\pi,f)$ as a monomial of order $m=kp$ in asymptotic dimensions of the form
$$
b_{1}(\pi,f)=d_{1}^{\,r_1(\pi)}d_{2}^{\,r_2(\pi)}\ldots d_{p+1}^{\,r_{p+1}(\pi)}
$$
where the exponents are natural numbers which depend on $\pi$ and for which 
$$
r_1(\pi)+r_2(\pi)+\ldots + r_{p+1}(\pi)=kp.
$$
Let us find an explicit combinatorial formula for these numbers. 

When computing $b_{1}(\pi_i,f)$
for various blocks $\pi_i$, we can see that asymptotic dimensions are assigned to the right legs
associated with letters $j$ and $j^{*}$ according to the rules
$$
j\rightarrow d_j\;\;{\rm and}\;\;j^{*}\rightarrow d_{j+1}
$$
for any $j\in [p]$. If $p=1$, this is immediate since the block containing the last letter $1^*$ in $W_k$ has to be colored by $2$
since $1$ is reserved to match the color of the imaginary block. Suppose that $p>1$.
If the right leg of some block $\pi_i$ is associated to letter $j$ in the word $W_k$ (and its left leg to $j^*$), 
then from the pair $\{j,j+1\}$ of possible colorings of this block we must choose $j$ since the inner block $\pi_s$ associated 
with the letter $j-1$ immediately preceding the considered letter $j$ can be colored by $j-1$ or $j$ (thus $j+1$ is not possible)
to make the partition adapted to the considered tuple. This argument is repeated for $j=p, \ldots , 1$ (in that order), which shows that 
each block containing a letter $j$ is colored by $j$ and therefore we should assign to it the dimension $d_j$.
A similar argument shows that each block containing a letter $j^*$ must be colored by $j+1$
and therefore we should assign to it $d_{j+1}$.

This leads to the formula
$$
r_j(\pi)=|\mathcal{R}_{j}(\pi)|+|\mathcal{R}_{j-1}^{*}(\pi)|\;\;{\rm for}\;j\in [p+1],
$$
where $\pi\in \mathcal{NC}_{m}^{2}(W_k)$, and thus
$$
\Psi_{1}((\eta\eta^{*})^{k})=\sum_{\pi \in \mathcal{NC}_{kp}^{2}(W_{k})}d_{1}^{\,r_2(\pi)}d_{2}^{\,r_2(\pi)}\ldots d_{p+1}^{\,r_{p+1}(\pi)},
$$
which completes the proof.
\hfill $\blacksquare$\\

The special case of $p=1$ corresponds to Wishart matrices and the Marchenko-Pastur distribution [27] with
shape parameter $t> 0$, namely
$$
\varrho_{t}={\rm max}\{1-t,0\}\delta_{0}+\frac{\sqrt{(x-a)(b-x)}}{2\pi x}1\!\!1_{[a,b]}(x)dx
$$
where $a=(1-\sqrt{t})^{2}$ and $b=(1+\sqrt{t})^{2}$ (see Corollary 9.1). At the same time, $\varrho_t$ 
is the free Poisson law in free probability.
It has been shown by Oravecz and Petz [29] that its moments are the 
{\it Narayana polynomials}
$$
N_{k}(t)=\sum_{j=1}^{k}N(k,j)t^{j}
$$ 
for any $k\in \mathbb{N}$, with coefficients
$$
N(k,j)=\frac{1}{j}{k-1 \choose j-1}{k \choose j-1}
$$ 
called {\it Narayana numbers}. These numbers are obtained in several different enumerations. 
For instance, $N(k,j)$ is equal to the number of Catalan paths of lenght $2k$ with $j$ peaks. 
However, we will use another enumeration related to right legs of blocks of the 
associated noncrossing partitions.
\begin{Proposition}
The Narayana number $N(k,j)$ is equal to the number of those noncrossing pair partitions $\pi\in \mathcal{NC}_{2k}^{2}$ 
which have $j$ even numbers in the set $\mathcal{R}(\pi)$.
\end{Proposition}
{\it Proof.}
Our proof is based on the enumeration derived by Osborn [30], in which $N(k,j)$ is the number of Catalan paths of 
lenght $2k$ with $2j$ edges lying in odd bands, where an odd band is a set of the form
${\mathbb R}\times [i,i+1]$ in the plane ${\mathbb R}^{2}$, where $i$ is even, in the standard setting, in which 
a Catalan path begins at $(0,0)$ and ends at $(2k,0)$. 
Now, let us associate with each Catalan path a noncrossing pair partition $\pi\in \mathcal{NC}_{2k}^{2}$
in the canonical way. Then, if $\{l,r\}$ is a block of $\pi$ with $l<r$, where $r$ is even, 
then the edges corresponding to $l$ and $r$ lie in the odd band. Similarly, if $r$ is odd, then these edges lie in the even band.
Consequently, $N(k,j)$ is equal to the number of noncrossing pair partitions of $[2k]$ which have $j$ right legs which are even 
and $k-j$ right legs which are odd. This completes the proof.
\hfill $\blacksquare$

\begin{Example}
{\rm In the path shown below all edges labelled by $a$ lie in the odd bands and all edges labelled by $b$
lie in the even bands.\\
\unitlength=1mm
\special{em.linewidth 0.5pt}
\linethickness{0.5pt}
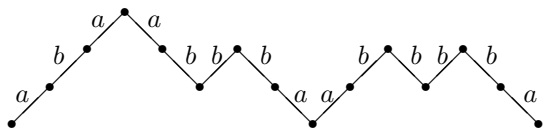
\begin{picture}(140.00,25.00)(-25.00,7.00)

\put(10.00,10.00){\line(1,1){5.00}}
\put(15.00,15.00){\line(1,1){5.00}}
\put(20.00,20.00){\line(1,1){5.00}}
\put(25.00,25.00){\line(1,-1){5.00}}
\put(30.00,20.00){\line(1,-1){5.00}}
\put(35.00,15.00){\line(1,1){5.00}}
\put(40.00,20.00){\line(1,-1){5.00}}
\put(45.00,15.00){\line(1,-1){5.00}}
\put(50.00,10.00){\line(1,1){5.00}}
\put(55.00,15.00){\line(1,1){5.00}}
\put(60.00,20.00){\line(1,-1){5.00}}
\put(65.00,15.00){\line(1,1){5.00}}
\put(70.00,20.00){\line(1,-1){5.00}}
\put(75.00,15.00){\line(1,-1){5.00}}

\put(10.00,10.00){\circle*{1.00}}
\put(15.00,15.00){\circle*{1.00}}
\put(20.00,20.00){\circle*{1.00}}
\put(25.00,25.00){\circle*{1.00}}
\put(30.00,20.00){\circle*{1.00}}
\put(35.00,15.00){\circle*{1.00}}
\put(40.00,20.00){\circle*{1.00}}
\put(45.00,15.00){\circle*{1.00}}
\put(50.00,10.00){\circle*{1.00}}
\put(55.00,15.00){\circle*{1.00}}
\put(60.00,20.00){\circle*{1.00}}
\put(65.00,15.00){\circle*{1.00}}
\put(70.00,20.00){\circle*{1.00}}
\put(75.00,15.00){\circle*{1.00}}
\put(80.00,10.00){\circle*{1.00}}

\put(10.50,13.00){\footnotesize $a$}
\put(15.50,18.00){\footnotesize $b$}
\put(20.50,23.00){\footnotesize $a$}
\put(28.00,23.00){\footnotesize $a$}
\put(33.00,18.00){\footnotesize $b$}
\put(36.50,18.00){\footnotesize $b$}
\put(43.00,18.00){\footnotesize $b$}
\put(47.50,13.00){\footnotesize $a$}
\put(51.00,13.00){\footnotesize $a$}
\put(56.00,18.00){\footnotesize $b$}
\put(63.00,18.00){\footnotesize $b$}
\put(66.50,18.00){\footnotesize $b$}
\put(73.00,18.00){\footnotesize $b$}
\put(78.00,13.00){\footnotesize $a$}
\end{picture}}\\
\end{Example}
\noindent
There are 6 edges lying in the odd bands and 8 edges lying in the even bands, 
thus there are $j=3$ even numbers in the set ${\mathcal R}(\pi)=\{4,5,7,8,11,13,14\}$.
\begin{Corollary}
If $p=1$ and $d_1>0$, then under the assumptions of Theorem 9.1 we obtain
$$
\lim_{n \rightarrow \infty}\tau_{1}(n)\left(\left(B(n)B^{*}(n)\right)^{k}\right)=d_{1}^{\,k}N_{k}(d_{2}/d_{1}),
$$
and thus the limit distribution is the $d_{1}$-dilation of the Marchenko-Pastur 
distribution $\varrho_t$ with the shape parameter $t=d_{2}/d_1$, denoted $\varrho_{d_2,\,d_1}$.
\end{Corollary}
{\it Proof.}
The letters $1$ and $1^{*}$ correspond to odd and even numbers, respectively, if we 
use the notation as in the proof of Theorem 9.1.
In turn, the number of noncrossing pair partitions which have $j$ right legs labelled by 
$1^{*}$ is given by $N(k,j)$ by Proposition 9.1 and thus
$$
P_{k}(d_1,d_2)=\sum_{j=1}^{k}N(k,j)d_{2}^{j}d_{1}^{k-j}=d_{1}^{k}N_{k}(d_{2}/d_{1})
$$
which proves the first assertion. The second one is clear since the 
transformation on moments $m_{k}\rightarrow \lambda^{k}m_{k}$ leads to the dilation of the corresponding measures.
\hfill $\blacksquare$\\
\begin{Remark}
{\rm The usual settings for Wishart matrices are slightly different since different normalizations are used [10,19]. For instance, if 
$B(N)\in GRM(m(N),N,1/N)$ is a Gaussian random matrix of dimension $m(N)\times N$ and $1/N$ is the variance
of each entry, where $N\in \mathbb {N}$, it is assumed that
$$
t=\lim_{N\rightarrow \infty}\frac{m(N)}{N}
$$
and then one computes the limit distribution of $B^{*}(N)B(N)$ under normalized trace ${\rm tr}(N)$ composed with classical expectation.
One obtains Narayana polynomials in $t$ as the limit moments, with 
the corresponding Marchenko-Pastur distribution with shape parameter $t$. In our model, we embed 
$B(n)$ and $B^{*}(n)$ in a larger square matrix of dimension $n$, in which $B(n)=S_{1,2}(n)$ and 
$B^{*}(n)=S_{1,2}^{*}(n)=S_{2,1}(n)$ are off-diagonal blocks and we compute the limit moments of $B(n)B^{*}(n)$ under $\tau_{1}(n)$.
In order to directly compare these two approaches, we set $N=n_1(n)$ and $m(N)=n_2(n)$ to get
$t=d_2/d_1$ and, since the variance in our approach is $1/n$ and in Marchenko-Pastur's theorem the variance is $1/N$, the $k$th 
moment in the limit moment obtained in Marchenko-Pastur's theorem must be multiplied by $d_{1}^{\,k}$ 
to give our asymptotic product polynomial.}
\end{Remark}

\begin{Example}
{\rm If $p=1$, the lowest order polynomials of Definition 9.1 are 
\begin{eqnarray*}
P_{1}(d_1,d_2)&=&d_2\\
P_{2}(d_1,d_2)&=&d_{2}d_{1}+d_{2}^{2}\\
P_{3}(d_1,d_2)&=&d_{2}d_{1}^{2}+3d_{2}^{2}d_{1}+d_{2}^{3}\\
P_{4}(d_1,d_2)&=&d_2d_1^3+6d_2^2d_1^2+6d_{2}^{3}d_{1}+d_{2}^{4}
\end{eqnarray*}
and can be obtained directly from Corollary 9.1 since the corresponding Narayana polynomials are 
$N_{1}(t)=t, N_{2}(t)=t+t^2, N_{3}(t)=t+3t^{2}+t^{3}$, $N_{4}(t)=t+6t^2+6t^3+t^4$, respectively.}
\end{Example}

\begin{Corollary}
If $d_1=d_2=\ldots =d_{p+1}=d$, then 
$$
P_{k}(d_1,d_2, \ldots , d_{p+1})=d^{kp}F(p,k)
$$
where $F(p,k)=\frac{1}{pk+1}{pk+k \choose k}$ are Fuss-Catalan numbers.
\end{Corollary}
{\it Proof.}
If $d_1=\ldots =d_{p+1}=d$, then $P_{k}(d_1,d_2,\ldots , d_{p+1})$ is equal to $d^{kp}$ multiplied by the 
number of noncrossing pair partitions adapted to the word $W_{k}=12\ldots pp^*\ldots 2^*1^*$ by Theorem 9.1. 
It is well known (see, for instance [20]) that the latter is equal to the Fuss-Catalan number $F(p,k)$.
\hfill $\blacksquare$\\
\begin{Example}
{\rm If $p=2$, we obtain by Theorem 9.1 the limit moments
$$
\Psi_{1}((\eta\eta^{*})^{k})=P_{k}(d_{1},d_{2}, d_{3}), \;\;{\rm where}\; \eta=\eta_{1,2}\eta_{2,3}.
$$
For instance, contributions to $P_{2}(d_1,d_2,d_3)$ are
\begin{eqnarray*}
\Psi_1(\wp_{2,1}^{*}\wp_{3,2}^{*}\wp_{3,2}\wp_{2,1}\wp_{2,1}^{*}\wp_{3,2}^{*}\wp_{3,2}\wp_{2,1})&=&d_{2}^{\,2}d_{3}^{\,2}\\
\Psi_1(\wp_{2,1}^{*}\wp_{3,2}^{*}\wp_{3,2}\wp_{1,2}^{*}\wp_{1,2}\wp_{3,2}^{*}\wp_{3,2}\wp_{2,1})&=&d_{1}d_{2}d_{3}^{\,2}\\
\Psi_1(\wp_{2,1}^{*}\wp_{3,2}^{*}\wp_{2,3}^{*}\wp_{1,2}^{*}\wp_{1,2}\wp_{2,3}\wp_{3,2}\wp_{2,1})&=&d_1d_{2}^{\,2}d_{3}
\end{eqnarray*}
with the corresponding partitions of the associated word $S=122^{*}1^{*}122^{*}1^{*}$ 

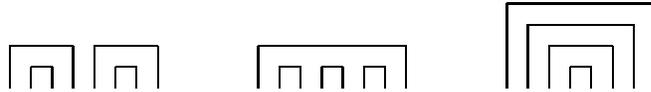
\begin{picture}(250.00,60.00)(00.00,0.00)
\put(76.00,10.00){\line(0,1){16.00}}
\put(84.00,10.00){\line(0,1){8.00}}
\put(92.00,10.00){\line(0,1){8.00}}
\put(100.00,10.00){\line(0,1){16.00}}

\put(76.00,26.00){\line(1,0){24.00}}
\put(84.00,18.00){\line(1,0){8.00}}

\put(108.00,10.00){\line(0,1){16.00}}
\put(116.00,10.00){\line(0,1){8.00}}
\put(124.00,10.00){\line(0,1){8.00}}
\put(132.00,10.00){\line(0,1){16.00}}

\put(108.00,26.00){\line(1,0){24.00}}
\put(116.00,18.00){\line(1,0){8.00}}

\put(170.00,10.00){\line(0,1){16.00}}
\put(178.00,10.00){\line(0,1){8.00}}
\put(186.00,10.00){\line(0,1){8.00}}
\put(194.00,10.00){\line(0,1){8.00}}
\put(202.00,10.00){\line(0,1){8.00}}
\put(210.00,10.00){\line(0,1){8.00}}
\put(218.00,10.00){\line(0,1){8.00}}
\put(226.00,10.00){\line(0,1){16.00}}

\put(170.00,26.00){\line(1,0){56.00}}
\put(178.00,18.00){\line(1,0){8.00}}
\put(194.00,18.00){\line(1,0){8.00}}
\put(210.00,18.00){\line(1,0){8.00}}

\put(264.00,10.00){\line(0,1){32.00}}
\put(272.00,10.00){\line(0,1){24.00}}
\put(280.00,10.00){\line(0,1){16.00}}
\put(288.00,10.00){\line(0,1){8.00}}
\put(296.00,10.00){\line(0,1){8.00}}
\put(304.00,10.00){\line(0,1){16.00}}
\put(312.00,10.00){\line(0,1){24.00}}
\put(320.00,10.00){\line(0,1){32.00}}

\put(264.00,42.00){\line(1,0){56.00}}
\put(272.00,34.00){\line(1,0){40.00}}
\put(280.00,26.00){\line(1,0){24.00}}
\put(288.00,18.00){\line(1,0){8.00}}

\end{picture}
\\
respectively. Thus,
$$
P_{2}(d_1,d_2,d_3)=d_{2}^{\,2}d_{3}^{\,2}+d_{1}d_{2}d_{3}^{\,2}+ d_1d_{2}^{\,2}d_{3}.
$$
If we define multivariate generalizations of Narayana polynomials $N_{k}(t_1, t_2, \ldots , t_p)$ by the formula
$$
P_{k}(d_1,d_2, \ldots , d_{p+1})=d_{1}^{kp}N_{k}(d_2/d_1,d_3/d_1, \ldots, d_{p+1}/d_1),
$$
we obtain
$$
N_{2}(t_1,t_2)=t_1^2t_2^2+t_1t_2^2+t_1^2t_2.
$$ 
These polynomials and their coefficients should
be of some interest from the combinatorial point of view.
}
\end{Example}

Let us write the limit distribution of Theorem 9.1 in the convolution form.
Let $\mu_s=U_{s}\mu$ be the probability measure on the real line defined by 
$$
G_{\mu_s}(z)=sG_{\mu}(z)+\frac{1-s}{z}
$$
for any $s>0$, where $G_{\mu_s}$ and $G_{\mu}$ are Cauchy transforms of $\mu_s$ and $\mu$, respectively.
\begin{Theorem}
If $d_1,d_2, \ldots ,d_{p+1}>0$ and $p>1$, then the limit distribution of Theorem 9.1 takes the form
$$
\mu_1=\varrho_{d_2,d_1}\boxtimes \varrho_{d_3,d_1}\boxtimes \ldots \boxtimes \varrho_{d_{p+1},d_{1}}.
$$
\end{Theorem}
{\it Proof.}
We have
$$
\tau_1(n)\left(\left(S_{1,2}(n)\ldots S_{p,p+1}(n)S_{p,p+1}^{*}(n)\ldots S_{1,2}^{*}(n)\right)^{k}\right)
$$
$$
=
n_2/n_1\cdot \tau_2(n)\left(\left(S_{2,3}(n)\ldots S_{p,p+1}(n)S_{p,p+1}^{*}(n)\ldots S_{2,3}^{*}(n)S_{1,2}^{*}(n)S_{1,2}(n)\right)^{k}\right)
$$
$$
=
n_2/n_1\cdot \tau_2(n)\left(\left(T_{2,3}(n)\ldots T_{p,p+1}(n)T_{p,p+1}(n)\ldots T_{2,3}(n)T_{1,2}(n)T_{1,2}(n)\right)^{k}\right)
$$
which tends to 
$$
d_2/d_1\cdot \Psi_2((\widetilde{\omega}^{2}\widehat{\omega}_{1,2}^{2})^{k}),
$$
as $n\rightarrow \infty$, where $\widetilde{\omega}=\widehat{\omega}_{2,3}\ldots \widehat{\omega}_{p,p+1}$. Now, the pair
$\{ \widetilde{\omega}^{2}, \widehat{\omega}_{1,2}^{2}\}$ is free with respect to $\Psi_2$
since the only common index present in both operators is $2$ and therefore any polynomial in one of them, say $\widehat{\omega}_{1,2}^{2}$, 
which is in the kernel of $\Psi_2$ maps any vector onto a linear combination of vectors which do not begin with $2$ and 
these are in the kernel of the other one, say $\widetilde{\omega}^{2}$. 
Now, the $\Psi_2$-distribution of $\widehat{\omega}_{1,2}^{2}$ is $\varrho_{d_1,d_2}$ by Corollary 9.1. 
Therefore, 
$$
\mu_1=U_{d_2/d_1}(\varrho_{d_1,d_2}\boxtimes \mu_{2})
$$
where $\mu_{2}$ is the $\Psi_2$-distribution of $\widetilde{\omega}^{2}$. If $\widetilde{\omega}$ is a product
of at least two operators, the same procedure is applied to $\mu_{2}$, which gives 
$$
\mu_{2}=U_{d_3/d_2}(\varrho_{d_2,d_3}\boxtimes \mu_{3}),
$$
where $\mu_3$ is the $\Psi_3$-distribution of $(\widehat{\omega}_{3,4}\ldots \widehat{\omega}_{p,p+1})^{2}$.
We continue this inductive procedure and observe that the last step gives 
$$
\mu_{p-1}=U_{d_p/d_{p-1}}(\varrho_{d_{p-1},d_p}\boxtimes \varrho_{d_{p+1},d_p}).
$$
Now, it is easy to show that the S-transform of $\mu'=U_{s}\mu$ takes the form
$$
S_{\mu'}(z)=\frac{1+z}{s+z}S_{\mu}\left(\frac{z}{s}\right).
$$
Therefore,
\begin{eqnarray*}
S_{\mu_{p-1}}(z)&=&(1+z)(d_p/d_{p-1}+z)^{-1}(d_{p-1}+d_{p-1}z)^{-1}(d_{p+1}+d_{p-1}z)^{-1}\\
&=&(d_{p}+d_{p-1}z)^{-1}(d_{p+1}+d_{p-1}z)^{-1}
\end{eqnarray*}
The next step gives
\begin{eqnarray*}
S_{\mu_{p-2}}(z)
&=&
(d_{p-1}+d_{p-2}z)^{-1}(d_{p}+d_{p-2}z)^{-1}(d_{p+1}+d_{p-2}z)^{-1},
\end{eqnarray*}
and continuing in this fashion, we finally get
$$
S_{\mu_1}(z)=(d_{2}+d_{1}z)^{-1}(d_{3}+d_{1}z)^{-1}\ldots (d_{p+1}+d_{1}z)^{-1},
$$
which completes the proof.
\hfill $\blacksquare$\\

Let us finally establish a relation between the limit distributions of Theorem 9.1 and 
free Bessel laws $\pi_{p,\,t}$ of Banica {\it et al} [4] expressed in terms of
the free multiplicative convolution 
$$
\pi_{p,\,t}=\pi^{\,\boxtimes (p-1)}\boxtimes \pi^{\,\boxplus\, t},
$$
where $\pi=\varrho_1$ is the standard Marchenko-Pastur distribution with the 
shape parameter equal to one. Random matrix models for these laws given in [4] were 
based on the multiplication of independent Gaussian random matrices. 

It is easy to observe that one can keep the normalization of Gaussian variables in terms 
of the parameter $n\rightarrow \infty$ in Theorems 8.1 and 9.1 without assuming that $n$ is the dimension of $Y(u,n)$.
In the context of Theorem 9.1, the asymptotic moments are still
$P_{k}(d_{1}, \ldots , d_{p+1})$, except that $d_1, \ldots , d_{p+1}$ are arbitrary non-negative numbers.
In particular, we can set $n_1=\ldots =n_p=n$ and $n_{p+1}$ to be the integer part of $tn$ (for $p=1$,
this reminds the normalization discussed in Remark 9.1). In this case, $d_1=\ldots = d_{p}=1$ and $d_{p+1}=t$, and
the limit laws are free Bessel laws. 

\begin{Corollary}
If $d_1=\ldots =d_p=1$ and $d_{p+1}=t$, the polynomials $P_{k}(d_{1}, \ldots ,d_{p+1})$ are moments of the free Bessel law $\pi_{p,\,t}$.
\end{Corollary}
{\it Proof.}
This is an immediate consequence of Theorem 9.2.
\hfill $\blacksquare$

\section{Appendix}
The definition of the symmetrically matricially free array of units given in [23] should be strenghtened
in order that the symmetrized Gaussian operators be symmetrically matricially free. This requires certain changes to be made in [23] 
which are listed below.
\begin{enumerate}
\item
Condition (2) of Definition 8.1 in [23] should be replaced by condition (2) of Definition 3.4 given in this paper. 
For that purpose, one needs to distinguish even and odd elements.
\item
The proof of symmetric matricial freeness of the array of symmetrized Gaussian operators was stated
in Proposition 8.1 in [23] and the proof was (unfortunately) omitted. Using Definition 3.4, this fact is proved here in Proposition 3.5.
\item
Definition 2.3 in [23] should be replaced by Definition 4.3 given in this paper. In principle, it is possible to modify the old definition and use
conditions on intersections of unordered pairs, but it is much more convenient and simpler to phrase the new definition
using sequences of ordered pairs.
\item
The proof of Theorem 9.1 in [23] needs to be slightly strengthened since the classes
$\mathcal{NC}_{m,q}^{2}(\{p_1,q_1\}, \ldots , \{p_m, q_m\})$ can be smaller, in general, than those considered in [23] since Definition 4.3 is 
stronger than Definition 2.3 in [23].
Therefore, we need to justify that if $\pi(\gamma)\in \mathcal{NC}_{m,q}^{2}\setminus \mathcal{NC}_{m,q}^{2}(\{p_1,q_1\}, \ldots , \{p_m,q_m\})$, then
the corresponding mixed moment of symmetric blocks 
under $\tau_{q}(n)$ vanishes. Clearly, in order that this moment be non-zero, the imaginary block of $\pi(\gamma)$ must be colored by $q$. Then all blocks
of $\pi(\gamma)$ of depth zero must be labeled by pairs $\{p_i,q_i\}$ which contain $q$ in order that 
the corresponding symmetric blocks $T_{p_i,q_i}(n)$ act non-trivially on vectors $e_{j}, j\in N_{q}$. 
Supposing that $q_i=q$ for all such pairs $\{p_i,q_i\}$ assigned to the right legs of $\pi(\gamma)$, 
we set $\mathpzc{v}_{i}=(p_i,q_i)$. 
Then, in turn, all symmetric blocks which correspond to blocks of $\pi(\gamma)$ of depth one whose nearest outer block 
is labeled by given $\{p_i,q_i\}$ must be labeled by pairs $\{p_k,q_k\}$ which contain $p_i$. Supposing that $q_k=p_i$
for all such pairs $\{p_k,q_k\}$ assigned to the right legs of $\pi(\gamma)$, we set $\mathpzc{v}_{k}=(p_k,q_k)$. We continue in this fashion until all symmetric blocks and the corresponding blocks of $\pi(\gamma)$ are taken into account. 
Finally, if $\{i,j\}$ is a block, where $i<j$ and $\mathpzc{v}_{j}=(p_j,q_j)$, then we must have $\mathpzc{v}_i=(q_j,p_j)$ in order that
the action of $T_{p_i,q_i}(n)$ be non-trivial, which follows from an inductive argument starting from the deepest blocks. 
Consequently, in order to get a non-zero contribution 
from the corresponding mixed moment of symmetric blocks, which is proportional to the product of variances as in the proof in [23],
there must exist a tuple $(\mathpzc{v}_{1}, \ldots , \mathpzc{v}_{m})$ such that $\mathpzc{v}_{i}\in \{(p_i,q_i), (q_i,p_i)\}$
to which $\pi(\gamma)$ is adapted. Therefore, the conditions of Definition 4.3 are satisfied. 
\end{enumerate}

\noindent\\[10pt]
{\bf Acknowledgement}\\
I would like to thank the referee for inspiring remarks which enabled me
to generalize the results contained in the original manuscript. The referee's
suggestions were also very helpful in revising several other parts of the text.

This work is partially supported by the Polish Ministry of Science and Higher Education (Wroclaw University of Technology grant 
S20089/I18).

\end{document}